\documentclass[preprint,12pt]{elsarticle}

\usepackage{amssymb, amsmath, amsthm, amscd, mathrsfs}

\usepackage{bigints}
\usepackage{dsfont}
\usepackage{empheq}
\usepackage{cases}
\usepackage{enumitem}

\usepackage{hyperref}
\usepackage{cancel}

\usepackage{algorithm}%
\usepackage{algorithmicx}%

\theoremstyle{plain}
\newtheorem{theorem}{Theorem}[section]

\newtheorem{lemma}[theorem]{Lemma}
\newtheorem{proposition}[theorem]{Proposition}
\theoremstyle{definition}
\newtheorem{definition}[theorem]{Definition}
\newtheorem*{notation}{Notation}
\newtheorem{remark}{Remark}
\numberwithin{equation}{section}

\def \d {\mathrm{d}}
\def \B {\mathbf{B}}
\usepackage{breqn}
\usepackage{xcolor}
\usepackage{algorithm}
\usepackage{algpseudocode}
\usepackage{amsmath}
\allowdisplaybreaks
\journal{***********
}

\begin{document}
		
\begin{frontmatter}



\title{Insensitizing controls of a volume-surface reaction-diffusion equation with dynamic boundary conditions}


\author[1]{Idriss Boutaayamou}
\ead{d.boutaayamou@uiz.ac.ma}

\affiliation[1]{organization={Lab-SIV},
	addressline={Polydisciplinary Faculty-Ouarzazate, Ibnou Zohr University}, 
	city={Ouarzazate},
	postcode={45000}, 
	state={B.P. 638},
	country={Morocco}}

\author[2]{Fouad Et-tahri\corref{cor1}}
\cortext[cor1]{Corresponding author. 
This research did not receive any specific grant from funding agencies in the public, commercial, or not-for-profit sectors.}
\ead{fouad.et-tahri@edu.uiz.ac.ma}

\affiliation[2]{organization={Lab-SIV},
            addressline={Faculty of Sciences-Agadir, Ibnou Zohr University}, 
            city={Agadir},
            state={B.P. 8106},
            country={Morocco}}

\author[3,4]{Lahcen Maniar}
\ead{maniar@uca.ma}

\affiliation[3]{organization={Department of Mathematics},
            addressline={Faculty of Sciences Semlalia, LMDP, UMMISCO (IRD-UPMC), Cadi Ayyad University}, 
            city={Marrakesh},
            postcode={40000}, 
            state={B.P. 2390},
            country={Morocco}}

\affiliation[4]{organization={University Mohammed VI Polytechnic, Vanguard Center},
	city={Benguerir},
	country={Morocco}}

\begin{abstract}
This paper deals with the insensitizing controllability property of the quasilinear parabolic equation with dynamic boundary conditions. This problem can be reformulated as a null controllability problem for a cascade quasilinear system with dynamic boundary conditions. To this end, we approach the problem by first dealing with null controllability in the framework of an inhomogeneous linearized system. Next, we derive new estimates of control and state, allowing us to apply a local inversion theorem to obtain null controllability of the quasilinear system.
\end{abstract}



\begin{keyword}
Reaction-diffusion equations \sep quasilinear parabolic equations \sep volume-surface \sep insensitizing controls \sep null controllability \sep dynamic boundary conditions.
\MSC[2020] 35K55 \sep 35K57 \sep 35K59 \sep 93C20.

\end{keyword}
\end{frontmatter}
\tableofcontents

\section{Introduction}
In this paper, we focus on analyzing the property of the insensitizing controllability of reaction-diffusion equations with dynamic boundary conditions of the surface diffusion type (generalized Wentzell type), which means that controls are resistant to small, unknown perturbations of the initial state. The equation we are dealing with reads as follows:
\begin{equation}\label{equation_quasi-linear}
	\left\{
	\begin{aligned}
		&\psi_{t}-\nabla\cdot\left(\sigma(\psi)\nabla \psi\right) +a(\psi) =f+
		v\mathds{1}_{\omega} & & \text {in}\; \Omega_T, \\
		&\psi_{\Gamma,t}-\nabla_{\Gamma}\cdot\left(\delta(\psi_{\Gamma})\nabla_{\Gamma} \psi_{\Gamma}\right)+\sigma(\psi_{\Gamma})
		\partial_{\nu} \psi + b(\psi_{\Gamma})=f_{\Gamma} & & \text {on}\;\Gamma_T, \\
		& \psi_{\Gamma}= \psi_{|_{\Gamma}}  & & \text {on}\;\Gamma_T, \\
		& (\psi(\cdot,0), \psi_{\Gamma}(\cdot,0))=(\psi_{0}+ \tau \widehat{\psi_0}, \psi_{0,\Gamma}+ \tau_{\Gamma} \widehat{\psi_{0,\Gamma}})  & & \text {in } \Omega\times\Gamma.
	\end{aligned}
	\right.
\end{equation} 
Here $\Omega$ is a bounded domain in $\mathbb{R}^{d}$ ($d\geq 1$) with boundary $\Gamma$ of class $C^{2}$, $\omega\subset\Omega$ an open subset of $\Omega$, for every subset $D$ of $\overline{\Omega}$, we denote $D_T:= D\times (0,T)$ and $\mathds{1}_{D}$ the characteristic function of $D$,
$(\psi, \psi_{\Gamma})$ is the state of the system, $(\psi_{0},\psi_{0,\Gamma})$ is the initial state perturbed by $(\tau \widehat{\psi_0},\tau_{\Gamma} \widehat{\psi_{0,\Gamma}})$,  $\tau, \tau_{\Gamma}\in\mathbb{R}$ are unknown and small enough, $v$ is the function control and $(f,f_{\Gamma})$ is the source term. We denote by $\psi_{|_{\Gamma}}$ the trace of $\psi$ on $\Gamma$, $\partial_{\nu}$ the normal derivative associated to the outward normal $\nu$ of $\Omega$, $\nabla_{\Gamma}$ and $\nabla_{\Gamma}\cdot:=\text{div}_{\Gamma}$ are tangential gradient and tangential divergence, respectively.\\ 
\par 
Given the existence of global solutions of \eqref{equation_quasi-linear} for certain initial state, source and control spaces $X_0$, $X_1$ and $V$ (see Proposition \ref{well-posdness of quasi}), 
we consider the energy function:
\begin{eqnarray}
	\mathcal{J}(\psi,\psi_{\Gamma})=\frac{\theta}{2}\int_{\mathcal{O}_{T}}|\psi(x,t,\tau,\tau_{\Gamma},v)|^{2}\d x \d t + \frac{\theta_\Gamma}{2}\int_{\Sigma_{T}}|\psi_{\Gamma}(x,t,\tau,\tau_{\Gamma},v)|^{2}\d S \d t, \label{energy function}
\end{eqnarray}
where $\theta>0$ and $\theta_\Gamma\geq 0$ are positive real constants, $\mathcal{O}$ and $\Sigma$ are given nonempty open subsets of $\Omega$ and $\Gamma$, respectively, and $(\psi,\psi_{\Gamma})=(\psi(\cdot,t,\tau,\tau_{\Gamma},v), \psi_{\Gamma}(\cdot,t,\tau,\tau_{\Gamma},v))$ is the corresponding (global) solution of the equation \eqref{equation_quasi-linear} associated to $\tau, \tau_{\Gamma}$ and $v$. Now, we introduce the following definition of insensitizing controls for \eqref{equation_quasi-linear}:
\begin{definition} \label{definition of insensitizing controls}
	For given functions $(f,f_{\Gamma})\in X_1$, $(\psi_0,\psi_{0,\Gamma})\in X_0$, a control function $v\in V$ is said to insensitize $\mathcal{J}$ if for all $(\widehat{\psi_0},\widehat{\psi_{0,\Gamma}})\in X_0$ with $\|(\widehat{\psi_0},\widehat{\psi_{0,\Gamma}}\|_{X_0}=1$ the corresponding solution $(\psi,\psi_{\Gamma})$ of \eqref{equation_quasi-linear} satisfies 
	\begin{eqnarray}
		\frac{\partial \mathcal{J}}{\partial\tau}\bigg|_{\tau=\tau_{\Gamma}=0}=\frac{\partial \mathcal{J}}{\partial\tau_\Gamma}\bigg|_{\tau=\tau_{\Gamma}=0}=0. \label{definition of insensitizing}
	\end{eqnarray}
\end{definition}
\noindent \textbf{A brief overview of related literature.}
The issue of insensitizing controls was first introduced by J. L. Lions in \cite{Lions1992SentinellesPL, lions1990quelques}, which subsequently led to numerous studies on the topic, both for hyperbolic and parabolic equations.
The existence of insensitizing controls for parabolic equations with Dirichlet or Neumann boundary conditions has been widely studied by numerous authors (see, for example, \cite{teresa2000insensitizing, bodart2002insensitizing, bodart2004local, bodart2004insensitizing, micu2004example, liu2012insensitizing, boyer2019insensitizing, huaman2023insensitizing}).
However, there are few results on insensitizing controllability for parabolic equations with dynamic boundary conditions. To the best of our knowledge, the paper \cite{zhang2019insensitizing} is the first to address a parabolic equation with semilinearity only in the bulk, while \cite{santos2024insensitizing} discusses an insensitizing control problem involving tangential gradient terms for a specific semilinear equation. As far as we know our paper is the first contribution in the quasilinear framework, both in the bulk and on the surface. Moreover, in the context of equations with dynamic boundary conditions, the functional $\mathcal{J}$ encompasses both the energy within the bulk and the localized surface energy, whereas in the literature, the usual functional only accounts for the local energy in the bulk.\\ \\
\textbf{Structure of this paper.} In Section \ref{Section2} we collect some notation, assumptions, preliminaries
and important tools. Next, Section \ref{Section3} introduces the main result (Theorem \ref{Main result}), reformulated as the null controllability of a quasilinear cascade system, with the proof of this reformulation detailed in \ref{Appendix B}. The discussion in Section \ref{Section4} focuses on the well-posedness and regularity properties required for reformulating the problem for specific linear systems, alongside the proof of the well-posedness for the quasilinear system, as shown in \ref{Appendix A}. Following this, Section \ref{Section5} establishes the null controllability of the inhomogeneous linearized cascade system, accompanied by state estimates. The local null controllability of the quasilinear system is addressed in Section \ref{Section6}, including the proof of main result. Lastly, Section \ref{Secion7} provides concluding remarks and final comments.
\section{Preliminaries, notation and assumptions} \label{Section2}
\subsection{Preliminaries}
Let us first introduce some basic notation. Let $d\in \mathbb N$ and $\Omega\subset\mathbb{R}^{d}$ is a bounded domain with boundary $\Gamma$ of class $C^2$. For any $s\geq 0$ and $p\in [1,\infty]$, the Lebesgue and Sobolev spaces for functions mapping from $\Omega$
to $\mathbb{R}$ are denoted as $L^p(\Omega)$ and $W^{s,p}(\Omega)$. We write $\|\cdot\|_{L^{p}(\Omega)}$ and $\|\cdot\|_{W^{s,p}(\Omega)}$ to denote the standard norms on these spaces. In the case $p=2$, we
use the notation $H^{s}(\Omega)=W^{s,2}(\Omega)$. We use the same notation for Lebesgue and Sobolev spaces on $\Gamma$. For any Banach space $X$, the Bochner
spaces of functions mapping from an interval $I$ into $X$ are denoted by $L^{p}(I;X)$ and $W^{s,p}(I;X)$, and the space
$C(I;X)$ denotes the set of continuous functions mapping from $I$ to $X$. The natural state space for our problems is $$\mathbb{L}^{2}:=L^{2}\left(\overline{\Omega}; \d x \otimes\d S\right),$$
where 
$\d x$ denotes the Lebesgue measure on $\Omega$ and $\d S$ denotes the natural surface measure on $\Gamma$. This space can be identified with  
$$ L^{2}(\Omega, dx)\otimes L^{2}(\Gamma, \d S).$$
It is a Hilbert space endowed with the following inner product
\begin{eqnarray*}
	\langle (y,y_{\Gamma}), (z,z_{\Gamma}) \rangle_{\mathbb{L}^{2}}:=\langle y, z \rangle_{L^{2}(\Omega)}+ \langle y_{\Gamma}, z_{\Gamma} \rangle_{L^{2}(\Gamma)}.
\end{eqnarray*}
We will also denote $\mathbb{L}^{\infty}:=L^{\infty}(\Omega)\otimes L^{\infty}(\Gamma)$.
This is a Banach space equipped with the norm
$$\|(y,y_{\Gamma})\|_{\mathbb{L}^{\infty}}:=\|y\|_{L^{\infty}(\Omega)}+ \|y_{\Gamma}\|_{L^{\infty}(\Gamma)}.$$
We introduce the Sobolev-type spaces needed in the sequel:
\begin{eqnarray*}
	\mathbb{H}^{s}:=\{(y,y_{\Gamma})\in H^{s}(\Omega)\times H^{s}(\Gamma)\;:\; y_{|_{\Gamma}}=y_{\Gamma}\},\quad \text{for}\; s\geq 1.
\end{eqnarray*}
The following continuous embeddings hold when $d\leq 3$ (see \cite[Theorem 4.12]{Adams} and \cite[Theorem 2.2]{lee1987yamabe}):
\begin{eqnarray}
	H^2(M)\hookrightarrow L^\infty(M)\quad \text{and}\quad H^2(M)\hookrightarrow W^{1,4}(M), \quad M=\Omega\;\text{or}\; \Gamma. \label{continuous embeddings}
\end{eqnarray}
\par 
To describe the surface heat diffusion, we need to define certain differential operators on $\Gamma$, defined locally in terms of the standard Riemannian metric on $\Gamma$, see \cite{jost2008riemannian}. In this paper, we will not use the local formulas that define these operators, but rather the relevant properties such as the surface divergence theorem. These operators can be defined by extensions as follows, we refer to \cite{choulli2009introduction}:
\begin{itemize}
	\item We define the tangential gradient $\nabla_{\Gamma} y_{\Gamma}$ for any
	smooth function $y_{\Gamma}$ on $\Gamma$ by
	\begin{eqnarray}
		\nabla_{\Gamma} y_{\Gamma}:=\nabla y-(\partial_{\nu}y)\nu, \label{grad}
	\end{eqnarray}
	where $y$ is an extension of $y_{\Gamma}$ in a neighborhood of $\Gamma$. It can be seen as the projection of the standard Euclidean gradient $\nabla y$ onto the tangent space on $\Gamma$.
	\item We define the tangential divergence $\mbox{div}_{\Gamma}Y_{\Gamma}$ for any
	smooth vector field $Y_{\Gamma}$ on $\Gamma$
	by
	\begin{eqnarray}
		\mbox{div}_{\Gamma} Y_{\Gamma}:=\mbox{div} Y-Y^{\prime}\nu\cdot\nu, \label{div}
	\end{eqnarray}
	where $Y^{\prime}=(\partial_{i}Y_{j})$ and $Y$ is an extension of $Y_{\Gamma}$ in a neighborhood of $\Gamma$. Note that formulas \eqref{grad} and \eqref{div} do not depend on the chosen extension, and $\mbox{div}_{\Gamma}(Y_{\Gamma})$ can be considered as a continuous linear form on $H^{1}(\Gamma)$
	\begin{eqnarray*}
		\mbox{div}_{\Gamma}(Y_{\Gamma})\;:\; H^{1}(\Gamma)\longrightarrow \mathbb{R}, \quad z_{\Gamma}\mapsto -\int_{\Gamma} \langle Y_{\Gamma}, \nabla_{\Gamma}z_{\Gamma} \rangle_{\Gamma} \d S,
	\end{eqnarray*}
	where $\langle \;,\;\rangle_{\Gamma}$ is the Riemannian inner product of tangential vectors on $\Gamma$. In the following, we will denote $\nabla_{\Gamma}\cdot Y_{\Gamma}$ instead of $\mbox{div}_{\Gamma}(Y_{\Gamma})$ and $\cdot$ instead of $\langle \;,\;\rangle_{\Gamma}$.
	\item The Laplace-Beltrami $\Delta_{\Gamma}y_{\Gamma}$ is defined by 
	\begin{eqnarray*}
		\Delta_{\Gamma}y_{\Gamma}=	\mbox{div}_{\Gamma}(\nabla_{\Gamma}y_{\Gamma}) \quad \forall y_{\Gamma}\in H^{2}(\Gamma).
	\end{eqnarray*}
	In
	particular, 
	the Stokes divergence theorem on $\Gamma$ holds, see \cite{taylor1996partial},
	\begin{eqnarray}
		\int_{\Gamma} \Delta_{\Gamma}y_{\Gamma} z_{\Gamma}\d S=-\int_{\Gamma} \nabla_{\Gamma} y_{\Gamma}\cdot\nabla_{\Gamma}z_{\Gamma}\d S \quad \forall y_{\Gamma}\in H^{2}(\Gamma),\;\forall z_{\Gamma}\in H^{1}(\Gamma). \label{Stokes}
	\end{eqnarray}
\end{itemize}
\begin{remark}
	In the one-dimensional case, $\Gamma$ is a manifold of dimension $0$. Consequently, the tangential operators are trivial: $\nabla_{\Gamma}=0$ and $\Delta_{\Gamma}=0$.
\end{remark}

\subsection{Notation}
	Throughout this paper, we adopt the following notations:
	\begin{itemize} 
		\item In order to simplify the presentation of the estimates we introduce the following notation
		\begin{eqnarray*}
			&&\nabla Y :=(\nabla y, \nabla_{\Gamma}y_{\Gamma}), \quad  \Delta Y := (\Delta y, \Delta_{\Gamma}y_{\Gamma}), \quad Y_{t} :=(y_{t}, y_{\Gamma,t}), \quad Y_{tt} := (y_{tt}, y_{\Gamma,tt}),
		\end{eqnarray*}
		where $Y:=(y,y_{\Gamma})$, $\nabla_{\Gamma}$ is the tangential gradient, $\Delta_{\Gamma}$ is the Laplace-Beltrami, $Y_{t}$ and $Y_{tt}$ respectively denote the first and second-order partial derivative of $Y$ with respect to $t$.
		\item 	We introduce the following energy spaces:
		\begin{eqnarray*}
			&&\mathfrak{E}_T:= H^{1}(0,T;\mathbb{L}^{2})\cap L^{2}(0,T;\mathbb{H}^{2})\quad\mbox{and}\quad \mathfrak{F}_T:= H^{1}(0,T;\mathbb{H}^{2})\cap L^{2}(0,T;\mathbb{H}^{4}).
		\end{eqnarray*}
		We recall the usual continuous embedding (see Theorem III.4.10.2 in \cite{amann1995linear}):
		\begin{eqnarray*}
			\mathfrak{E}_T\hookrightarrow C([0,T];\mathbb{H}^1)\quad \text{and}\quad \mathfrak{F}_T\hookrightarrow C([0,T];\mathbb{H}^3).
		\end{eqnarray*}
		In particular, for $d\leq 3$, using \eqref{continuous embeddings}, we obtain the existence of a constant $C>0$ such that
		\begin{eqnarray}
			 \|Y\|_{L^{\infty}(0,T;\mathbb{L}^{\infty})}, \|\nabla Y\|_{L^{\infty}(0,T;\mathbb{L}^{\infty})}\leq C\|Y\|_{\mathfrak{F}_T} \quad \forall Y\in \mathfrak{F}_T. \label{Bounded gradient}
		\end{eqnarray}
		\item For any $Y=(y,y_{\Gamma})\in\mathfrak{E}_T$, we will denote 
		\begin{eqnarray} \label{n1} 
			\begin{cases}
				\textbf{L}_{1}Y:=y_{t}-\sigma(0)\Delta y+a^{\prime}(0)y,\quad \\
				\textbf{L}_{2}Y:=y_{\Gamma,t}-\delta(0)\Delta_{\Gamma}y_{\Gamma}+\sigma(0)\partial_{\nu}y+b^{\prime}(0)y_{\Gamma},\\
				\textbf{L}_{1}^{\star}Y:=-y_{t}-\sigma(0)\Delta y+a^{\prime}(0)y,\\
				\textbf{L}_{2}^{\star}Y:=-y_{\Gamma,t}-\delta(0)\Delta_{\Gamma}y_{\Gamma}+\sigma(0)\partial_{\nu}y+b^{\prime}(0)y_{\Gamma}.
			\end{cases}
		\end{eqnarray}
		\item  For any Banach space 
		$X$, we denote by $B_X(a,r)$
		and $\overline{B}_X(a,r)$ the open and closed balls centered at $a$ with radius 
		$r$, respectively.
		\item We omit the infinitesimals $\d t$, $\d x$ and $\d S$ as they can be deduced by looking at the integration domain.
		\item The symbol $C$ will stand for a generic positive constant depending
		on $\Omega$, $\omega$, $s$, $\lambda$, $T$, $\theta$, $\theta_\Gamma$, $\sigma$, $\delta$, $a$ and $b$, where $s$ and $\lambda$ are known parameters in Carleman estimates.
	\end{itemize}
\subsection{Assumptions}
We make the following general assumptions.
\begin{enumerate}[label=(A\arabic*), ref=(A\arabic*)]
	\item \label{hyp:A1} $\Omega\subset \mathbb{R}^{d}$ is a bounded domain with boundary of class $C^{2}$ and $d\in\{1,2,3\}$.
	\item \label{hyp:A2} $\mathcal{O}$ and $\Sigma$ are given nonempty open subsets of $\Omega$ and $\Gamma$, respectively. 
	\item \label{hyp:A3} $\omega\Subset\Omega$ an open subset such that $\omega\cap \mathcal{O}\neq \varnothing$, where $\omega\Subset\Omega$ denotes $\overline{\omega}\subset\Omega$.
	\item \label{hyp:A4} $(\psi_{0},\psi_{0,\Gamma})=(0,0)$.
	\item \label{hyp:A5} $(\widehat{\psi_0},\widehat{\psi_{0,\Gamma}})\in \mathbb{H}^3$ with $\|(\widehat{\psi_0},\widehat{\psi_{0,\Gamma}})\|_{\mathbb{H}^3}=1$.
	\item  \label{hyp:A6} The parameters $\tau, \tau_{\Gamma}\in\mathbb{R}$ are unknown and small enough.
	\item  \label{hyp:A7} The diffusion coefficients satisfy
	\begin{eqnarray*}
		\sigma,\delta\in C^{3}(\mathbb{R},\mathbb{R}),\quad \sigma(r)\geq \rho \quad\mbox{and}\quad \delta(r)\geq \rho,
	\end{eqnarray*}
	for all $r\in\mathbb{R}$ for some $\rho>0$.
	\item \label{hyp:A8} The reaction terms verify
	\begin{eqnarray*}
		a,b\in C^{2}(\mathbb{R},\mathbb{R}) \quad\mbox{and}\quad a(0)=b(0)=0.
	\end{eqnarray*}
\end{enumerate}
\section{Main result} \label{Section3}
We state the main result of this paper.
\begin{theorem} \label{Main result} Suppose that assumptions \ref{hyp:A1}-\ref{hyp:A8} hold. Then, there exist constants $\varepsilon>0$ and  $C=C(\Omega,\omega,T)>0$ such that for any $F:=(f,f_{\Gamma})\in \mathfrak{E}_T$ verifying
	\begin{eqnarray}
		\|F\|^{2}_{L^{2}(0,T;\mathbb{H}^{2})}+\|e^{C/t}F\|^{2}_{L^{2}(0,T;\mathbb{L}^{2})}+ \|e^{C/t}F_{t}\|^{2}_{L^{2}(0,T;\mathbb{L}^{2})} <\varepsilon, \label{assump source}
	\end{eqnarray}
	one can find a control $v\in H^{1}(0,T;L^{2}(\omega))\cap L^{2}(0,T;H^{2}(\omega))$ that insensitizes the functional $\mathcal{J}$ defined by \eqref{energy function}. Moreover there exist a constant $C^{\prime}>0$ such that 
	\begin{eqnarray*}
		\|v\|^{2}_{H^{1}(0,T;L^{2}(\omega))\cap L^{2}(0,T;H^{2}(\omega))}\leq C^{\prime}\left(\|e^{C/t}F\|^{2}_{L^{2}(0,T;\mathbb{L}^{2})}+ \|e^{C/t}F_{t}\|^{2}_{L^{2}(0,T;\mathbb{L}^{2})}\right).
	\end{eqnarray*}
\end{theorem}
\noindent \textbf{Reduction of the insensitizing problem.} The computation of partial derivatives of the functional defined in \eqref{energy function} and a duality argument, allows us to reformulate our insensitivity problem as a problem of null controllability of a cascade system. Specifically, we present the following cascade system:
\begin{equation} \label{cascade quasi-linear system}
	\left\{
	\begin{aligned}
		&\psi_{t}-\nabla\cdot\left(\sigma(\psi)\nabla \psi\right) +a(\psi) =f+
		v\mathds{1}_{\omega}& & \text {in}\; \Omega_T, \\
		&-h_{t}-\sigma(\psi)\Delta h  +a^{\prime}(\psi)h =\theta\psi\mathds{1}_{\mathcal{O}} & & \text {in}\; \Omega_T, \\
		&\psi_{\Gamma,t}-\nabla_{\Gamma}\cdot\left(\delta(\psi_{\Gamma})\nabla_{\Gamma} \psi_{\Gamma}\right)+\sigma(\psi_{\Gamma})
		\partial_{\nu} \psi + b(\psi_{\Gamma})=f_{\Gamma} & & \text {on}\;\Gamma_T, \\
		&-h_{\Gamma,t}-\delta(\psi_{\Gamma})\nabla_{\Gamma} h_{\Gamma}+\sigma(\psi_{\Gamma})
		\partial_{\nu} h+b^{\prime}(\psi_{\Gamma})h_{\Gamma}=\theta_{\Gamma}\psi_{\Gamma}\mathds{1}_{\Sigma} & & \text {on}\;\Gamma_T, \\
		& \psi_{\Gamma}=\psi_{|_{\Gamma}}, \; h_{\Gamma}=h_{|_{\Gamma}}  & & \text {on}\;\Gamma_T, \\
		& (\psi(\cdot,0),\psi_{\Gamma}(\cdot,0))=(0,0), \; (h(\cdot,T),h_{\Gamma}(\cdot,T))=(0,0) & & \text {in } \Omega\times\Gamma. 
	\end{aligned}
	\right.
\end{equation}
The reformulation is as follows:
\begin{lemma} \label{reformulation of problem}
	Suppose that assumptions \ref{hyp:A1}, \ref{hyp:A2} and \ref{hyp:A4}-\ref{hyp:A8} hold. Let $(f,f_\Gamma)\in L^{2}(0,T;\mathbb{H}^{2})$ and $v\mathds{1}_{\omega}\in  L^{2}(0,T;H^{2}(\Omega))$ such that the system \eqref{equation_quasi-linear} has a unique global solution $\Psi=(\psi,\psi_{\Gamma})\in \mathfrak{F}_{T}$. Then, the following statements are equivalent:
	\begin{itemize}
		\item The solution $(\Psi,H)$ of the
		cascade system \eqref{cascade quasi-linear system} associated with $v$ verifies 
		\begin{eqnarray}
			(h(\cdot,0),h_{\Gamma}(\cdot,0))=(0,0)\; \text{in} \;\Omega\times\Gamma. \label{Null in ident}
		\end{eqnarray}
		\item The control $v$ insensitizes the functional $\mathcal{J}$ in the sense of Definition \ref{definition of insensitizing controls}.
	\end{itemize}
\end{lemma}
	\par The local null controllability of \eqref{cascade quasi-linear system} (in the sense of \eqref{Null in ident}) is reformulated as a surjectivity of the mapping $\varLambda:\mathbb{X}\longrightarrow\mathbb{Y}$:
\begin{eqnarray}
	&&\varLambda(\Psi, H, v):=((\varLambda_{1}(\Psi, H, v), \varLambda_{3}(\Psi, H, v)), (\varLambda_{2}(\Psi, H, v), \varLambda_{4}(\Psi, H, v))), \label{varlamda}
\end{eqnarray}
where
\begin{eqnarray*}
	&&\varLambda_{1}(\Psi, H, v):= \psi_{t}-\nabla\cdot\left(\sigma(\psi)\nabla \psi\right) +a(\psi)-
	v\mathds{1}_{\omega},\\
	&&\varLambda_{2}(\Psi, H, v):= -h_{t}-\sigma(\psi)\Delta h  +a^{\prime}(\psi)h -\theta\psi\mathds{1}_{\mathcal{O}}, \\
	&&\varLambda_{3}(\Psi, H, v):=\psi_{\Gamma,t}-\nabla_{\Gamma}\cdot\left(\delta(\psi_{\Gamma})\nabla_{\Gamma} \psi_{\Gamma}\right)+\sigma(\psi_{\Gamma})
	\partial_{\nu} \psi + b(\psi_{\Gamma}),\\
	&&\varLambda_{4}(\Psi, H, v):= -h_{\Gamma,t}-\delta(\psi_{\Gamma})\Delta_{\Gamma} h_{\Gamma}+\sigma(\psi_{\Gamma})
	\partial_{\nu} h +b^{\prime}(\psi_{\Gamma})h_{\Gamma}-\theta_{\Gamma}\psi_{\Gamma}\mathds{1}_{\Sigma}.       
\end{eqnarray*}
and $\mathbb{X},\mathbb{Y}=\mathbb{Y}_1\times \mathbb{Y}_2$ are appropriate spaces of the state-control (see Section \ref{Section6}). More precisely, Theorem \ref{Main result} is equivalent to 
\begin{eqnarray*}
	\exists\varepsilon>0,\;  \forall F=(f,f_\Gamma)\in B_{\mathbb{Y}_1}(0,\varepsilon),\; \exists (\Psi,v)\in\mathbb{X},\;\mbox{such that}\; \varLambda(\Psi, H, v)=(f, f_\Gamma, 0, 0).
\end{eqnarray*}
To achieve this, we apply the Lyusternik-Graves Inverse Mapping Theorem in infinite
dimensional spaces, whose proof can be referenced in \cite{alekseev1987optimal}.
\begin{theorem}[Lyusternik-Graves’ Theorem]
	\label{Lyusternik} Let $X$ and $Y$ be Banach spaces and let $\varLambda: B_X(0,r)\subset X\rightarrow Y$ be a $C^{1}$ mapping. Let us assume that the derivative $\varLambda^{\prime}(0):X\rightarrow Y$ is surjective and let
	us set $\xi_{0}=\varLambda(0)$. Then, there exist $\varepsilon> 0$, a mapping
	$W: B_Y(\xi_{0},\varepsilon)\subset Y\rightarrow X$ and a constant $C> 0$ satisfying:
	\begin{itemize}
		\item $W(z)\in B_X(0,r)$ and $\varLambda\circ W(z)=z\quad \forall z\in B_Y(\xi_{0},\varepsilon)$,
		\item $\|W(z)\|_{X}\leq C\|z-\xi_{0}\|_{Y}\quad \forall z\in B_Y(\xi_{0},\varepsilon)$.
	\end{itemize}
\end{theorem}
\par 
The most difficult task is therefore the choice of spaces $\mathbb{X}$ and $\mathbb{Y}$ for which mapping $\varLambda$ is well defined. The surjectivity of $\varLambda^{\prime}(0,0, 0)$ is linked to the null controllability of the following inhomogeneous linearized system (around zero) of \eqref{cascade quasi-linear system}:
\begin{equation} \label{linearised cascade system}
	\left\{
	\begin{aligned}
		&\psi_{t}-\sigma(0)\Delta\psi +a^{\prime}(0)\psi =f +v\mathds{1}_{\omega} & & \text {in}\; \Omega_T, \\
		&-h_{t}-\sigma(0)\Delta h  +a^{\prime}(0)h =g+\theta\psi\mathds{1}_{\mathcal{O}} & & \text {in}\; \Omega_T, \\
		&\psi_{\Gamma,t}-\delta(0)\Delta_{\Gamma}\psi_{\Gamma}+\sigma(0)
		\partial_{\nu} \psi + b^{\prime}(0)\psi=f_{\Gamma} & & \text {on}\;\Gamma_T, \\
		&-h_{\Gamma,t}-\delta(0)\Delta_{\Gamma} h_{\Gamma}+\sigma(0)
		\partial_{\nu} h+b^{\prime}(0)h_{\Gamma}=g_{\Gamma}+\theta_{\Gamma}\psi_{\Gamma}\mathds{1}_{\Sigma} & & \text {on}\;\Gamma_T, \\
		& \psi_{\Gamma}= \psi_{|_{\Gamma}},\; h_{\Gamma}=h_{|_{\Gamma}}  & & \text {on}\;\Gamma_T, \\
		& (\psi(\cdot,0),\psi_{\Gamma}(\cdot,0))=(0,0),\; (h(\cdot,T),h_{\Gamma}(\cdot,T))=(0,0) & & \text {in } \Omega\times\Gamma, \\
	\end{aligned}
	\right.
\end{equation}
Given suitable assumptions on $f, f_\Gamma, g$ and $g_{\Gamma}$, the system \eqref{linearised cascade system} is null controllable in the sens of \eqref{Null in ident}.\\

\noindent\textbf{Comments on assumptions \ref{hyp:A1}, \ref{hyp:A3} and \ref{hyp:A4}.} The condition $d\leq 3$ is used to estimate the terms $I_1$ and $I_3$ in the right hand side of \eqref{d2} via the embedding \eqref{continuous embeddings}, which is only valid when $d\leq 3$. Therefore, this question remains open for $d\geq 4$  as mentioned in \cite{huaman2023insensitizing, huaman2023local}.
The main result \ref{Main result} is based on assumption \ref{hyp:A4}, which allows for the establishment of a global Carleman estimate for the adjoint system of \eqref{linearised cascade system}. However, by slightly relaxing the functional $\mathcal{J}$:
\begin{eqnarray}
	\mathcal{J}_{t_0}(\psi,\psi_{\Gamma})=\frac{\theta}{2}\int_{\mathcal{O}\times(t_0, T)}|\psi(x,t,\tau,\tau_{\Gamma},v)|^{2}\d x \d t + \frac{\theta_\Gamma}{2}\int_{\Sigma\times (t_0, T)}|\psi_{\Gamma}(x,t,\tau,\tau_{\Gamma},v)|^{2}\d S \d t, \label{relaxing energy function}
\end{eqnarray}
where $t_0\in (0,T)$. We can find a control $v\in L^{2}(\omega_T)$ that insensitizes the functional $\mathcal{J}_{t_0}$ for sufficiently small initial data and source terms. Indeed, based on \cite{et2025null}, a control can be found that drives the state of \eqref{equation_quasi-linear} to zero at time $t_0$. Consequently, Theorem \ref{Main result} can be applied from $t_0$ onwards. The geometric condition \ref{hyp:A3} is also a technical condition in the proof of the Carleman estimate. However, the case where $\omega\cap \mathcal{O}=\varnothing$ is involved is also of significant interest. For further discussion and remarks on this topic, see \cite{ervedoza2022desensitizing}.


\section{Well-posedness results}\label{Section4}
In this subsection, we will briefly present the well-posedness of the systems that we will need in this paper.
\subsection{Linear systems}
 In this subsection, $\Omega\subset \mathbb{R}^{d}$ ($d\geq 1$) is a bounded domain with boundary $\Gamma$ of class $C^{2}$. Let us first consider the following linear equation:
\begin{equation} \label{linear equation 2}
	\left\{
	\begin{aligned}
		&\psi_{t}-\sigma(0)\Delta\psi +a^{\prime}(0)\psi =f & & \text {in}\; \Omega_T, \\
		&\psi_{\Gamma,t}-\delta(0)\Delta_{\Gamma}\psi_{\Gamma}+\sigma(0)
		\partial_{\nu} \psi + b^{\prime}(0)\psi=f_{\Gamma} & & \text {on}\;\Gamma_T, \\
		& \psi_{\Gamma}= \psi_{|_{\Gamma}}  & & \text {on}\;\Gamma_T, \\
		& (\psi(\cdot,0),\psi_{\Gamma}(\cdot,0))=(\psi_{0},\psi_{0,\Gamma}) & & \text {in } \Omega\times\Gamma, \\
	\end{aligned}
	\right.
\end{equation} 
We are interested in the following categories of solutions for \eqref{linear equation 2},  see\cite{maniar2017null}.
\begin{definition}\label{Well of linear equations}
	Let $F=(f,f_{\Gamma})\in L^{2}(0,T;\mathbb{L}^{2})$, $\Psi_{0}=(\psi_{0},\psi_{0,\Gamma})\in \mathbb{L}^{2}$. 
	\begin{enumerate}[label=(\arabic*)]
		\item A distributional solution (solution by transposition) of \eqref{linear equation 2} is a function $\Psi=(\psi,\psi_{\Gamma})\in L^{2}(0,T;\mathbb{L}^{2})$ such that for any $Z=(z,z_{\Gamma})\in \mathfrak{E}_{T}$ with $Z(\cdot,T)=0$, we have 
		\begin{eqnarray}
			&&\langle \Psi, \mathbf{L}^{*}Z\rangle_{L^{2}(0,T;\mathbb{L}^{2})} =\langle F, Z\rangle_{L^{2}(0,T;\mathbb{L}^{2})}+ \langle \Psi_0, Z(\cdot,0)\rangle_{\mathbb{L}^{2}}.   \label{w1}
		\end{eqnarray}
		where $\textbf{L}^{\star}:=(\textbf{L}^{\star}_1, \textbf{L}^{\star}_2)$.
		\item A strong solution of \eqref{linear equation 2} is a function $(\psi,\psi_{\Gamma})\in \mathfrak{E}_{T}$  fulfilling \eqref{linear equation 2} in $L^{2}(0,T;\mathbb{L}^{2})$.
	\end{enumerate}
\end{definition}
The well-posedness and regularity properties of the solutions to \eqref{linear equation 2} are based on semigroup theory, as studied in detail in \cite[Propositions 2.4 and 2.5]{maniar2017null} and \cite{ait2024internal}.
\begin{proposition} \label{well linear equation}
	Let $F=(f,f_{\Gamma})\in L^{2}(0,T;\mathbb{L}^{2})$. 
	\begin{enumerate}[label=(\arabic*)]
		\item If $\Psi_{0}=(\psi_{0},\psi_{0,\Gamma})\in \mathbb{L}^{2}$. Then, there exists a unique distributional solution $\Psi=(\psi,\psi_{\Gamma})\in  C([0,T];\mathbb{L}^{2})$
		of \eqref{linear equation 2}. Moreover, there is a constant $C>0$ such that 
		\begin{eqnarray}
			&&\|\Psi\|_{C([0,T];\mathbb{L}^{2})}
			\leq C\left( \|\Psi_{0}\|_{\mathbb{L}^{2}} + \|F\|_{L^{2}(0,T;\mathbb{L}^{2})} \right). \label{w2}
		\end{eqnarray}
		\item If $\Psi_{0}=(\psi_{0},\psi_{0,\Gamma})\in \mathbb{H}^{1}$. Then, there exists a unique strong solution $\Psi=(\psi,\psi_{\Gamma})\in \mathfrak{E}_{T}$
		of \eqref{linear equation 2}. Moreover, there is a constant $C>0$ such that 
		\begin{eqnarray}
			&&\|\Psi\|_{\mathfrak{E}_{T}}\leq C\left(\|\Psi_{0}\|_{\mathbb{H}^{1}} + \|F\|_{L^{2}(0,T;\mathbb{L}^{2})}\right). \label{w3}
		\end{eqnarray}
	\end{enumerate}
\end{proposition}
We now present the distributional solution and the strong solution of the linearized cascade system:
\begin{equation} \label{well linearized cascade system}
	\left\{
	\begin{aligned}
		&\psi_{t}-\sigma(0)\Delta\psi +a^{\prime}(0)\psi =f & & \text {in}\; \Omega_T, \\
		&-h_{t}-\sigma(0)\Delta h  +a^{\prime}(0)h =g+\theta\psi\mathds{1}_{\mathcal{O}} & & \text {in}\; \Omega_T, \\
		&\psi_{\Gamma,t}-\delta(0)\Delta_{\Gamma}\psi_{\Gamma}+\sigma(0)
		\partial_{\nu} \psi + b^{\prime}(0)\psi=f_{\Gamma} & & \text {on}\;\Gamma_T, \\
		&-h_{\Gamma,t}-\delta(0)\Delta_{\Gamma} h_{\Gamma}+\sigma(0)
		\partial_{\nu} h+b^{\prime}(0)h_{\Gamma}=g_{\Gamma}+\theta_{\Gamma}\psi_{\Gamma}\mathds{1}_{\Sigma} & & \text {on}\;\Gamma_T, \\
		& \psi_{\Gamma}= \psi_{|_{\Gamma}},\; h_{\Gamma}=h_{|_{\Gamma}}  & & \text {on}\;\Gamma_T, \\
		& (\psi(\cdot,0),\psi_{\Gamma}(\cdot,0))=(0,0),\; (h(\cdot,T),h_{\Gamma}(\cdot,T))=(0,0) & & \text {in } \Omega\times\Gamma, \\
	\end{aligned}
	\right.
\end{equation} 
\begin{definition}
	Let $F=(f,f_{\Gamma}), G=(g,g_{\Gamma})\in L^{2}(0,T;\mathbb{L}^{2})$.
	\begin{enumerate}[label=(\arabic*)]
		\item A distributional solution of \eqref{well linearized cascade system} is a functions $\Psi=(\psi,\psi_{\Gamma}), H=(h, h_{\Gamma})\in L^{2}(0,T;\mathbb{L}^{2})$ such that for any $Z=(z,z_{\Gamma}), W=(w, w_{\Gamma})\in \mathfrak{E}_T$ with $Z(\cdot,T)=Z(\cdot,0)=0$, we have 
		\begin{eqnarray*}
			\langle \Psi, \textbf{L}^{\star}Z -\textbf{B}W \rangle_{L^{2}(0,T;\mathbb{L}^{2})} + \langle H, \textbf{L}W\rangle_{L^{2}(0,T;\mathbb{L}^{2})}&=&\langle F, Z \rangle_{L^{2}(0,T;\mathbb{L}^{2})} + \langle G, W \rangle_{L^{2}(0,T;\mathbb{L}^{2})},
		\end{eqnarray*}
		where $\textbf{L} :=(\textbf{L}_1,\textbf{L}_2)$, $\textbf{L}^{\star}:=(\textbf{L}^{\star}_1, \textbf{L}^{\star}_2)$ and $\textbf{B}W:=(\theta w\mathds{1}_{\mathcal{O}}, \theta_{\Gamma}w_{\Gamma}\mathds{1}_{\Sigma})$.
		\item A strong solution of \eqref{well linearized cascade system} is a functions $\Psi=(\psi,\psi_{\Gamma}), H=(h, h_{\Gamma})\in \mathfrak{E}_{T}$  fulfilling \eqref{well linearized cascade system} in $L^{2}(0,T;\mathbb{L}^{2})$.
	\end{enumerate}
\end{definition}
The well-posedness and regularity properties of the solutions to \eqref{well linearized cascade system} are based on the cascade structure of the system and on Proposition \ref{well linear equation}.
\begin{proposition} 
	Let $F=(f,f_{\Gamma}), G=(g,g_{\Gamma})\in L^{2}(0,T;\mathbb{L}^{2})$. Then, 
	\begin{enumerate}[label=(\arabic*)]
		\item There exists a unique distributional solution $\Psi=(\psi,\psi_{\Gamma}), H=(h,h_{\Gamma})\in  C([0,T];\mathbb{L}^{2})$
		of \eqref{well linearized cascade system}. Moreover, there is a constant $C>0$ such that 
		\begin{eqnarray}
			&&\|\Psi\|_{C([0,T];\mathbb{L}^{2})} +\|H\|_{C([0,T];\mathbb{L}^{2})}
			\leq C\left( \|F\|_{L^{2}(0,T;\mathbb{L}^{2})}+ \|G\|_{L^{2}(0,T;\mathbb{L}^{2})} \right). \label{w2}
		\end{eqnarray}
		\item There exists a unique strong solution $\Psi=(\psi,\psi_{\Gamma}), H=(h,h_{\Gamma})\in \mathfrak{E}_{T}$
		of \eqref{well linearized cascade system}.\\ Moreover, there is a constant $C>0$ such that 
		\begin{eqnarray}
			&&\|\Psi\|_{\mathfrak{E}_{T}}+ \|H\|_{\mathfrak{E}_{T}}\leq C\left( \|F\|_{L^{2}(0,T;\mathbb{L}^{2})}+ \|G\|_{L^{2}(0,T;\mathbb{L}^{2})}\right). \label{energy estimate}
		\end{eqnarray}
	\end{enumerate}
\end{proposition}
\begin{remark}
	The adjoint systems of \eqref{linear equation 2} and \eqref{well linearized cascade system} yield the same results, due to the principal operator of the equation \eqref{linear equation 2} is self-adjoint, see \cite[Proposition 2.1]{maniar2017null}.
\end{remark}
\subsection{Nonlinear system}
Now, we present the well-posedness and regularity properties of solutions to:
\begin{equation}\label{quasi-linear equation well}
	\left\{
	\begin{aligned}
		&\psi_{t}-\nabla\cdot\left(\sigma(\psi)\nabla \psi\right) +a(\psi) =f & & \text {in}\; \Omega_T, \\
		&\psi_{\Gamma,t}-\nabla_{\Gamma}\cdot\left(\delta(\psi_{\Gamma})\nabla_{\Gamma} \psi_{\Gamma}\right)+\sigma(\psi_{\Gamma})
		\partial_{\nu} \psi + b(\psi_{\Gamma})=f_{\Gamma} & & \text {on}\;\Gamma_T, \\
		& \psi_{\Gamma}= \psi_{|_{\Gamma}}  & & \text {on}\;\Gamma_T, \\
		&(\psi(\cdot,0),\psi_{\Gamma}(\cdot,0))=(\psi_{0},\psi_{0,\Gamma}) & & \text {in } \Omega\times\Gamma.
	\end{aligned}
	\right.
\end{equation}
Note that we have proven in \cite{et2025null} that the equation \eqref{quasi-linear equation well} admits local solutions in $\mathfrak{F}_{T^{\prime}}$ for all $T^{\prime}<T_{max}$, and there exists a control and a global solution with regularity different from that of $\mathfrak{F}_{T^{\prime}}$. However, for our current problem, we will need global solutions that belong to $\mathfrak{F}_{T}$ (at least for sufficiently small data). To this end, we will present the Proposition \ref{well-posdness of quasi}, the proof of which will be summarized in \ref{Appendix A}. Existence is guaranteed by Lyusternik-Graves' Theorem \ref{Lyusternik}, while uniqueness is ensured by an energy estimate and Gronwall's inequality.
\begin{proposition} \label{well-posdness of quasi}
	Assume that assumptions \ref{hyp:A1},  \ref{hyp:A7} and \ref{hyp:A8} hold. Then, there exists a constants $\kappa>0$ such that for any $F=(f,f_{\Gamma})\in L^{2}(0,T;\mathbb{H}^{2})$ and $\Psi_0=(\psi_{0},\psi_{0,\Gamma})\in \mathbb{H}^{3}$ verifying
	\begin{eqnarray*}
		\|F\|^{2}_{L^{2}(0,T;\mathbb{H}^{2})}+\|\Psi_{0}\|^{2}_{\mathbb{H}^{3}}<\kappa,
	\end{eqnarray*}
	the equation \eqref{quasi-linear equation well} has a unique global solution $(\psi,\psi_{\Gamma})\in \mathfrak{F}_{T}$. Moreover, there is a constant $C>0$ such that 
	\begin{eqnarray}
		\|\Psi\|_{\mathfrak{F}_{T}}\leq C\left( \|F\|_{L^{2}(0,T;\mathbb{H}^{2})}+ \|\Psi_{0}\|_{\mathbb{H}^{3}}\right). \label{quasi linear energy estimate}
	\end{eqnarray}
\end{proposition}

\section{Carleman estimate for cascade system with dynamic boundary conditions} \label{Section5}
In this section, $\Omega\subset \mathbb{R}^{d}$ ($d\geq 1$) is a bounded domain with boundary $\Gamma$ of class $C^{2}$ and 
assumptions \ref{hyp:A2}, \ref{hyp:A3} are satisfied.
\subsection{Null controllability of \eqref{linearised cascade system}} \label{subsection 2.3}
The aim of this section is to prove null controllability for the system \eqref{linearised cascade system}, we also prove estimates on the state and regularity on the control which require some regularity of the source terms. First, let us recall the definitions of several classical weights, frequently used in this framework, see \cite{fursikov1996controllability}. We consider the following positive weight functions $\alpha$ and $\xi$ which depend on $\Omega$ and $\omega$
$$\alpha(x, t)=\frac{e^{2\lambda m}-e^{\lambda\left(m+ \eta(x)\right)}}{t(T-t)}\quad \mbox{and} \quad \xi(x, t)=\frac{e^{\lambda\left(m+ \eta(x)\right)}}{t(T-t)}.
$$
Here, $\lambda, m >1$ and $\eta=\eta(x)$ is a function in $C^2(\overline{\Omega})$ satisfying
\begin{equation}
	\label{eta}
	\eta>0 \text { in } \Omega, \quad \eta=0 \text { on } \Gamma,\quad  \inf_{\Omega \backslash \omega^{\prime}}\left|\nabla \eta(x)\right| >0\quad\mbox{and}\quad \max_{\Omega}\eta=1,
\end{equation}
where $\omega^{\prime}$ is a nonempty open set of $\omega\cap \mathcal{O}$,
we also consider $\omega^{\prime\prime}$ and $\omega^{\prime\prime\prime}$ two open subsets of $\omega\cap \mathcal{O}$ such that 
$\omega^{\prime}\Subset\omega^{\prime\prime}\Subset\omega^{\prime\prime\prime}\Subset\omega\cap \mathcal{O}$. The following lemma is a Carleman estimate for the adjoint system of \eqref{linearised cascade system}:
\begin{equation} \label{adjoint cascade linear system}
	\left\{
	\begin{aligned}
		&-\phi_{t}-\sigma(0)\Delta\phi +a^{\prime}(0)\phi =f^{1}+
		\theta k\mathds{1}_{\mathcal{O}} & & \text {in}\; \Omega_T, \\
		&k_{t}-\sigma(0)\Delta k +a^{\prime}(0)k =g^{1} & & \text {in}\; \Omega_T, \\
		&-\phi_{\Gamma,t}-\delta(0)\Delta_{\Gamma}\phi_{\Gamma}+\sigma(0)
		\partial_{\nu}\phi + b^{\prime}(0)\phi=f^{1}_{\Gamma}+\theta_{\Gamma}k_{\Gamma}\mathds{1}_{\Sigma} & & \text {on}\;\Gamma_T, \\
		&k_{\Gamma,t}-\delta(0)\Delta_{\Gamma} k_{\Gamma}+\sigma(0)
		\partial_{\nu}k+b^{\prime}(0)k_{\Gamma}=g^{1}_{\Gamma} & & \text {on}\;\Gamma_T, \\
		& \phi_{\Gamma}= \phi_{|_{\Gamma}},\; k_{\Gamma}= k_{|_{\Gamma}}  & & \text {on}\;\Gamma_T, \\
		& (\phi(\cdot,T),\phi_{\Gamma}(\cdot,T))=(0,0),\; (k(\cdot,0),k_{\Gamma}(\cdot,0))=(0, 0) & & \text {in } \Omega\times\Gamma.
	\end{aligned}
	\right.
\end{equation} 
We introduce the following notation:
\begin{eqnarray*}
	&&I(\Phi,s,\lambda, t_1,t_2)
	:=\int_{\Omega\times (t_1,t_2)}e^{-2s\alpha}\left[(s\xi)^{-1}\left(|\phi_{t}|^{2}+|\Delta \phi|^{2}\right)+ \lambda^{2}(s\xi)|\nabla \phi|^{2}+\lambda^{4}(s\xi)^{3}|\phi|^{2}\right]\\
	&&+ \int_{\Gamma\times (t_1,t_2)}e^{-2s\alpha}\left[(s\xi)^{-1}\left(|\phi_{\Gamma,t}|^{2}+|\Delta_{\Gamma} \phi_{\Gamma}|^{2}\right)+ \lambda(s\xi)|\nabla_{\Gamma} \phi_{\Gamma}|^{2}+\lambda^{3}(s\xi)^{3}|\phi_{\Gamma}|^{2}\right]\\
	&&+ \int_{\Gamma\times (t_1,t_2)}e^{-2s\alpha}\lambda(s\xi)|\partial_{\nu}\phi|^{2},
\end{eqnarray*}
where $\Phi:=(\phi,\phi_{\Gamma})\in \mathfrak{E}_T$ and $0\leq t_1\leq t_2\leq T$.
\begin{lemma} \label{Lemma 1}
	There are constants $C_{1}>0$ and $\lambda_1, s_{1}\geq 1$ such that for any $s\geq s_1$, any $\lambda\geq\lambda_1$ and any $\Phi=(\phi,\phi_{\Gamma}), K=(k,k_{\Gamma})\in \mathfrak{E}_T$ solution of \eqref{adjoint cascade linear system}, we have the following estimate
	\begin{eqnarray}
		&&I(\Phi,s,\lambda, 0,T)+I(K,s,\lambda, 0,T) \leq C_{1}\left(s^{7}\lambda^{8}\int_{\omega^{\prime\prime\prime}_{T}}e^{-2s\alpha}\xi^{7}|\phi|^{2}   \right. \nonumber \\
		&& \left.  + \int_{\Omega_{T}}e^{-2s\alpha}\left(s^{3}\lambda^{4}\xi^{3}|f^{1}|^{2}+|g^{1}|^{2}\right)+\int_{\Gamma_T}e^{-2s\alpha}\left(|f^{1}_{\Gamma}|^{2}+|g^{1}_{\Gamma}|^{2}\right) \right). \label{estimate Lemma3.1}
	\end{eqnarray}
	Furthermore, $C_1$ and $\lambda_{1}$ only depend on $\Omega$ and $\omega$, and $s_1$ can be chosen of the form $s_1=C(T+T^{2})$, where $C$ only
	depends on $\Omega$, $\omega$, $\sigma(0)$, $\delta(0)$, $a^{\prime}(0)$ and $b^{\prime}(0)$.
\end{lemma}
\begin{proof}
	Using the Carleman estimate of \cite[Lemma 3.2]{maniar2017null} for $d\geq 2$ and \cite[Lemma 2]{ait2024internal} for $d=1$, we obtain
	\begin{eqnarray}
		I(\Phi,s,\lambda, 0, T) &\leq& C_{0}\left(s^{3}\lambda^{4}\int_{\omega^{\prime\prime}_T}e^{-2s\alpha}\xi^{3}|\phi|^{2} + \theta\int_{\mathcal{O}_T}e^{-2s\alpha}|k|^{2} + \theta_{\Gamma}\int_{\Sigma_T}e^{-2s\alpha}|k_\Gamma|^{2} \right. \nonumber \\
		&& \left.  + \int_{\Omega_T}e^{-2s\alpha}|f^{1}|^{2}  +\int_{\Gamma_T}e^{-2s\alpha}|f^{1}_{\Gamma}|^{2} \right) \label{CarlemanPhi}
	\end{eqnarray}
	and 
	\begin{eqnarray}
		I(K,s,\lambda, 0, T) &\leq& C_{0}\left(s^{3}\lambda^{4}\int_{\omega^{\prime\prime}_T}e^{-2s\alpha}\xi^{3}|k|^{2} + \int_{\Omega_T}e^{-2s\alpha}|g^{1}|^{2}  +\int_{\Gamma_T}e^{-2s\alpha}|g^{1}_{\Gamma}|^{2} \right), \label{CarlemanK}
	\end{eqnarray}
	for any $\lambda\geq\lambda_{0}$ and $s\geq s_{0}$ for some  $C_0>0$ and $\lambda_{0}>0$ only depend on $\Omega$ and $\omega$, and $s_0$ can be chosen of the form $s_0=C(T+T^{2})$, where $C$ only
	depends on $\Omega$, $\omega$, $\sigma(0)$, $\delta(0)$, $a^{\prime}(0)$ and $b^{\prime}(0)$.\\
	Let us introduce a function $\pi\in C^{\infty}_{0}(\omega^{\prime\prime\prime})$ satisfying $0\leq \pi\leq 1$ and $\pi=1$ in $\omega^{\prime\prime}$. Then 
	\begin{eqnarray}
		s^{3}\lambda^{4}\int_{\omega^{\prime\prime}_T}e^{-2s\alpha}\xi^{3}|k|^{2} &\leq& s^{3}\lambda^{4}\int_{\omega^{\prime\prime\prime}_T}e^{-2s\alpha}\xi^{3}\pi|k|^{2} \nonumber\\
		&= & s^{3}\lambda^{4}\theta^{-1}\int_{\omega^{\prime\prime\prime}_T}e^{-2s\alpha}\xi^{3}\pi k(-\phi_{t}-\sigma(0)\Delta\phi +a^{\prime}(0)\phi -f^{1}) \nonumber\\
		&=& I_1+I_2+I_3+I_4. \label{I1234}
	\end{eqnarray}
	Integrating by parts in time, $\Phi(\cdot,T)=K(\cdot,0)=0$, $|\alpha_{t}|,|\xi_{t}|\leq CT\xi^{2}$ and $s\geq CT$ provides  
	\begin{eqnarray*}
		I_1
		&=& 	-2s^{4}\lambda^{4}\theta^{-1}\int_{\omega^{\prime\prime\prime}_T}e^{-2s\alpha}\alpha_{t}\xi^{3}\pi k\phi	+ 3s^{3}\lambda^{4}\theta^{-1}\int_{\omega_T}e^{-2s\alpha}\xi_{t}\xi^{2}\pi k\phi + s^{3}\lambda^{4}\theta^{-1}\int_{\omega^{\prime\prime\prime}_T}e^{-2s\alpha}\xi^{3}\pi k_{t}\phi \\
		&\leq & C\left(s^{5}\lambda^{4}\int_{\omega^{\prime\prime\prime}_T}e^{-2s\alpha}\xi^{5} |k||\phi|  +s^{4}\lambda^{4}\int_{\omega^{\prime\prime\prime}_T}e^{-2s\alpha}\xi^{4} |k||\phi|  + s^{3}\lambda^{4}\int_{\omega^{\prime\prime\prime}_T}e^{-2s\alpha}\xi^{3} |k_{t}||\phi| \right).
	\end{eqnarray*}
	By the Cauchy Schwarz inequality, it follows that
	\begin{eqnarray}
		I_1
		&\leq & \frac{1}{8}s^{3}\lambda^{4}\int_{\omega^{\prime\prime\prime}_T}e^{-2s\alpha}\xi^{3} |k|^{2}+ \frac{1}{2} s^{-1}\int_{\omega^{\prime\prime\prime}_T}e^{-2s\alpha}\xi^{-1} |k_{t}|^{2} +Cs^{7}\lambda^{4}\int_{\omega^{\prime\prime\prime}_T}e^{-2s\alpha}\xi^{7}|\phi|^{2}. \label{Estimate I1}
	\end{eqnarray}
		For $I_2$, we have
	\begin{eqnarray*}
		I_2&=&-\sigma(0)s^{3}\lambda^{4}\theta^{-1}\int_{\omega^{\prime\prime\prime}_T}e^{-2s\alpha}\xi^{3}\pi k\Delta\phi=-\sigma(0)s^{3}\lambda^{4}\theta^{-1}\int_{\omega^{\prime\prime\prime}_T}\Delta(e^{-2s\alpha}\xi^{3}\pi k)\phi.
	\end{eqnarray*}
	A simple computation leads to
	\begin{eqnarray*}
		|\Delta(e^{-2s\alpha}\xi^{3}\pi k)|\leq C\left(s^{2}\lambda^{2} e^{-2s\alpha}\xi^{5} |k|+ s\lambda e^{-2s\alpha}\xi^{4}|\nabla k|+ e^{-2s\alpha}\xi^{3}|\Delta k|\right).
	\end{eqnarray*}
	Then, the Cauchy Schwarz inequality yields
	\begin{eqnarray}
		I_2&\leq &C\left(s^{5}\lambda^{6}\int_{\omega^{\prime\prime\prime}_T}e^{-2s\alpha}\xi^{5} |k||\phi| +s^{4}\lambda^{5}\int_{\omega^{\prime\prime\prime}_T}e^{-2s\alpha}\xi^{4} |\nabla k||\phi|+  s^{3}\lambda^{4}\int_{\omega^{\prime\prime\prime}_T}e^{-2s\alpha}\xi^{3} |\Delta k||\phi|\right) \nonumber \\
		&\leq& \frac{1}{8}s^{3}\lambda^{4}\int_{\omega^{\prime\prime\prime}_T}e^{-2s\alpha}\xi^{3} |k|^{2}+ \frac{1}{2}s\lambda^{2}\int_{\omega_T}e^{-2s\alpha}\xi |\nabla k|^{2} +\frac{1}{2}s^{-1}\int_{\omega^{\prime\prime\prime}_T}e^{-2s\alpha}\xi^{-1} |\Delta k|^{2} \nonumber \\
		&&+Cs^{7}\lambda^{8}\int_{\omega^{\prime\prime\prime}_T}e^{-2s\alpha}\xi^{7} |\phi|^{2}. \label{Estimate I2}
	\end{eqnarray}
	For the remaining terms, one can clearly see that
	\begin{eqnarray}
		I_3&=& a^{\prime}(0)s^{3}\lambda^{4}\theta^{-1}\int_{\omega^{\prime\prime\prime}_T}e^{-2s\alpha}\xi^{3}\pi k\phi \nonumber \\
		&\leq & \frac{1}{8}s^{3}\lambda^{4}\int_{\omega^{\prime\prime\prime}_T}e^{-2s\alpha}\xi^{3} |k|^{2} + Cs^{3}\lambda^{4}\int_{\omega^{\prime\prime\prime}_T}e^{-2s\alpha}\xi^{3}|\phi|^{2} \label{Estimate I3}
	\end{eqnarray}
	and
	\begin{eqnarray}
		I_4&=&-s^{3}\lambda^{4}\theta^{-1}\int_{\omega^{\prime\prime\prime}_T}e^{-2s\alpha}\xi^{3}\pi k f^{1} \nonumber\\
		&\leq & \frac{1}{8}s^{3}\lambda^{4}\int_{\omega^{\prime\prime\prime}_T}e^{-2s\alpha}\xi^{3} |k|^{2} + Cs^{3}\lambda^{4}\int_{\omega^{\prime\prime\prime}_T}e^{-2s\alpha}\xi^{3}|f^{1}|^{2}. \label{Estimate I4}
	\end{eqnarray}
	By summing \eqref{CarlemanPhi} and \eqref{CarlemanK} and then, applying the estimates \eqref{I1234}-\eqref{Estimate I4}, we obtain the required estimate \eqref{estimate Lemma3.1}.
\end{proof}
We will deduce a Carleman estimate similar to \eqref{estimate Lemma3.1} with functions blowing up only at $t=0$.
Define the new weight functions:
\begin{eqnarray*}
	\beta(x, t)=\frac{e^{2\lambda m}-e^{\lambda\left(m+ \eta(x)\right)}}{\ell(t)}\quad \mbox{and} \quad \zeta(x, t)=\frac{e^{\lambda\left(m+ \eta(x)\right)}}{\ell(t)},\quad (x,t)\in\overline{\Omega}\times (0,T),
\end{eqnarray*}
where the function $\ell$ is given by
\begin{eqnarray*}
	\ell(t)=\begin{cases}  t(T-t)
		\quad & \mbox{if}\; t\in [0,T/2],\\
		\frac{T^{2}}{4}   \quad & \mbox{if}\; t\in [T/2,T].
	\end{cases}
\end{eqnarray*}
Note that $\ell\in C^{1}([0,T])$. An estimate with such weights is given in the following result. In the proof, we will use Lemma \ref{Lemma 1} and energy estimate \eqref{energy estimate}. We introduce the following notation:
\begin{eqnarray*}
	&&J(\Phi,s,\lambda, t_1,t_2)
	:=\int_{\Omega\times (t_1,t_2)}e^{-2s\beta}\left[\ell(t)\left(|\phi_{t}|^{2}+|\Delta \phi|^{2}\right)+ \ell^{-1}(t)|\nabla \phi|^{2}+\ell^{-3}(t)|\phi|^{2}\right] \\
	&&+ \int_{\Gamma\times (t_1,t_2)}e^{-2s\beta}\left[\ell(t)\left(|\phi_{\Gamma,t}|^{2}+|\Delta_{\Gamma} \phi_{\Gamma}|^{2}\right)+ \ell^{-1}(t)|\nabla_{\Gamma} \phi_{\Gamma}|^{2}+\ell^{-3}(t)|\phi_{\Gamma}|^{2}\right] \\
	&&+ \int_{\Gamma\times (t_1,t_2)}e^{-2s\beta}\ell^{-1}(t)|\partial_{\nu}\phi|^{2},
\end{eqnarray*}
where $\Phi:=(\phi,\phi_{\Gamma})\in \mathfrak{E}_T$ and $0\leq t_1\leq t_2\leq T$.
\begin{proposition} \label{P0}
	There exist constants $\lambda_1, s_{1}\geq 1$ such that for any $s\geq s_1$, any $\lambda\geq\lambda_1$, there exists a constant $C:=C(s,\lambda,T)>0$ which satisfies the following property: for any $\Phi=(\phi,\phi_{\Gamma}), K=(k,k_\Gamma)\in \mathfrak{E}_T$, we have the following estimate:
	\begin{eqnarray}
		&&J(\Phi,s,\lambda, 0,T)+J(K,s,\lambda, 0,T) \leq C_{1}\left(\int_{\omega^{\prime\prime\prime}_T}e^{-2s\beta}\ell^{-7}(t)|\phi|^{2} \right. \nonumber \\
		&& \left. + \int_{\Omega_T}e^{-2s\beta}\left(\ell^{-3}(t)|f^{1}|^{2}+|g^{1}|^{2}\right)  +\int_{\Gamma_T}e^{-2s\beta}\left(|f^{1}_{\Gamma}|^{2}+|g^{1}_{\Gamma}|^{2}\right)  \right). \label{estimateP0}
	\end{eqnarray}
	with $s_{1}$ and $\lambda_{1}$ as in Lemma \ref{Lemma 1}.
\end{proposition}
\begin{proof} 
	To do so, we split the left-hand side of \eqref{estimateP0} into two parts on $\left(0,T/2\right)$ and $\left(T/2,T\right)$. Firstly, on $\left(0,T/2\right)$. Using Carleman estimate \eqref{estimate Lemma3.1}, $\alpha=\beta$ and $\xi=e^{\lambda(m+\eta(x))}\ell^{-1}$ in $\overline{\Omega}\times \left(0,T/2\right)$, we obtain 
	\begin{eqnarray}
		&&J(\Phi,s,\lambda, 0,T/2)+J(K,s,\lambda, 0,T/2) \leq C\left(\int_{\omega^{\prime\prime\prime}_T}e^{-2s\alpha}\xi^{7}|\phi|^{2} \right. \nonumber\\
		&& \left. + \int_{\Omega_T}e^{-2s\alpha}\left(\xi^{3}|f^{1}|^{2}+|g^{1}|^{2}\right) +\int_{\Gamma_T}e^{-2s\alpha}\left(|f^{1}_{\Gamma}|^{2}+|g^{1}_{\Gamma}|^{2}\right)  \right). 
		\label{estimateP01}
	\end{eqnarray}
	Through a comparison of the weight functions, we obtain
	\begin{eqnarray}
		e^{-2s\alpha}\xi^{j}\leq C e^{-2s\beta}\ell^{-j}\quad \mbox{in}\;\; \overline{\Omega}\times (0,T),\; j=0,1,\label{estimateP02}
	\end{eqnarray}
	where $C>0$ only depending on $s$, $\lambda$ and $T$. 
	Considering \eqref{estimateP01} and \eqref{estimateP02}, we obtain 
	\begin{eqnarray}
		&&J(\Phi,s,\lambda, 0,T/2)+J(K,s,\lambda, 0,T/2) \leq C\left(\int_{\omega^{\prime\prime\prime}_T}e^{-2s\beta}\ell^{-7}(t)|\phi|^{2} \right. \nonumber\\
		&& \left. + \int_{\Omega_T}e^{-2s\beta}\left(\ell^{-3}(t)|f^{1}|^{2}+|g^{1}|^{2}\right) +\int_{\Gamma_T}e^{-2s\beta}\left(|f^{1}_{\Gamma}|^{2}+|g^{1}_{\Gamma}|^{2}\right)  \right). 
		\label{estimateP04}
	\end{eqnarray}
	Secondly, on $\left(T/2,T\right)$.
	Let us define a function $\vartheta\in C^{1}([0,T])$ such that 
		\begin{eqnarray*}
			\vartheta(t)= 0\quad\mbox{in}\quad \left[0,T/4\right] \quad \mbox{and}\quad \vartheta(t)= 1\quad\mbox{in}\quad \left[T/2,T\right].
	\end{eqnarray*}
	Put  $Y:=(y,y_{\Gamma}), Z:=(z,z_\Gamma)$ defined by $Y:= \vartheta\Phi, Z:=\vartheta K$, then $(Y, Z)$ is the strong solution of the system
	\begin{equation*}
		\left\{
		\begin{aligned}
			&\mathbf{L}_1^{\star}(Y)= \vartheta f^{1} +z\mathds{1}_{\mathcal{O}} - \vartheta^{\prime} \phi& & \text {in}\;\Omega_{T}, \\
			&\mathbf{L}_1(Z)= \vartheta g^{1} + \vartheta^{\prime} k& & \text {in}\;\Omega_{T}, \\
			& \mathbf{L}_{2}^{\star}(Y)=\vartheta f^{1}_{\Gamma} +z_{\Gamma}\mathds{1}_{\Sigma}- \vartheta^{\prime}\phi _{\Gamma}& & \text {on}\;\Gamma_{T} \\
			& \mathbf{L}_{2}(Z)=\vartheta g^{1}_{\Gamma} + \vartheta^{\prime}k _{\Gamma}& & \text {on}\;\Gamma_{T} \\
			& y_{\Gamma}= y_{|_{\Gamma}}, \; z_{\Gamma}= z_{|_{\Gamma}}  &&  \text {on}\;\Gamma_{T}, \\
			&(y(\cdot,T),y_{\Gamma}(\cdot,T))=(0,0),\; (z(\cdot,0),z_{\Gamma}(\cdot,0))=(0,0) & & \text {in } \Omega\times\Gamma.
		\end{aligned}
		\right.
	\end{equation*}
	We write the energy estimate \eqref{energy estimate} for $Y$ and, based on the definition $\vartheta$, it follows that
	\begin{eqnarray}
		&&\|\Phi\|^{2}_{L^{2}(T/2, T;\mathbb{L}^{2})} + \|\Phi_{t}\|^{2}_{L^{2}(T/2, T;\mathbb{L}^{2})}+ \|\nabla \Phi\|^{2}_{L^{2}(T/2, T;\mathbb{L}^{2})} +\|\Delta \Phi\|^{2}_{L^{2}(T/2,T;\mathbb{L}^{2})}   \label{estimateP03} \\
		&&  +\|\partial_{\nu} \phi\|^{2}_{L^{2}(T/2,T;L^{2}(\Gamma))}
		+\|K\|^{2}_{L^{2}(T/2,T;\mathbb{L}^{2})} + \|K_{t}\|^{2}_{L^{2}(T/2,T;\mathbb{L}^{2})}\nonumber\\
		&&  + \|\nabla K\|^{2}_{L^{2}(T/2,T;\mathbb{L}^{2})}+\|\Delta K\|^{2}_{L^{2}(T/2,T;\mathbb{L}^{2})}+\|\partial_{\nu} k\|^{2}_{L^{2}(T/2,T;L^{2}(\Gamma))}\nonumber\\
		&&\leq C\left(\|(f^{1},f^{1}_{\Gamma})\|^{2}_{L^{2}(T/4,T;\mathbb{L}^{2})}+ \|\Phi\|^{2}_{L^{2}(T/4,T/2;\mathbb{L}^{2})}+ \|(g^{1},g^{1}_{\Gamma})\|^{2}_{L^{2}(T/4,T;\mathbb{L}^{2})}+ \|K\|^{2}_{L^{2}(T/4,T/2;\mathbb{L}^{2})}  \right). \nonumber
	\end{eqnarray} 
	By the boundedness from above and from below of the weight functions $\beta$ and $\ell$ in $\overline{\Omega}\times\left(T/4, T\right)$ and \eqref{estimateP03}, we can derive that
	\begin{eqnarray*}
		&&J(\Phi,s,\lambda, T/2, T)+J(K,s,\lambda, T/2, T) \leq C\left(  \int_{\Omega_T}e^{-2s\beta}\left(\ell^{-3}(t)|f^{1}|^{2}+|g^{1}|^{2}\right)  
	   \right. \nonumber\\
		&& \left. +\int_{\Gamma_T}e^{-2s\beta}\left(|f^{1}_{\Gamma}|^{2}+|g^{1}_{\Gamma}|^{2}\right)
		+  \int_{T/4}^{T}\int_{\Omega} e^{-2s\beta}\ell^{-3}(t)(|\phi|^{2}+|k|^{2}) \right. \nonumber\\
		&& \left. + \int_{T/4}^{T}\int_{\Gamma} e^{-2s\beta}\ell^{-3}(t)(|\phi|^{2}+|k|^{2})   \right) \nonumber
	\end{eqnarray*}
	which, together with \eqref{estimateP04}, leads to \eqref{estimateP0}.
\end{proof}
Henceforth, we fix $s=s_{1}$, $\lambda=\lambda_{1}$, and introduce the following weights, which will be used in the sequel. 	
\begin{remark}
	Choosing  $m> 1$ large enough for instance $m> \frac{\log(5e^{\lambda}-4)}{\lambda}$, one finds that
	\begin{eqnarray*}
		\hat{\beta}(t)< (5/4) \check{\beta}(t) \quad \forall (x,t)\in \overline{\Omega}\times (0,T),
	\end{eqnarray*}
	where
	\begin{eqnarray*}
		\hat{\beta}(t):=\max_{x\in\overline{\Omega}}\beta (x,t),\quad \check{\beta}(t):=\min_{x\in\overline{\Omega}}\beta (x,t).
	\end{eqnarray*}
\end{remark} 
\begin{notation} Let us introduce the notations:
		\begin{eqnarray*}
			&& \gamma(t):=\hat{\beta}/5, \quad \mu(t):=e^{5s\gamma(t)}\ell^{3/2}, \quad \mu_{0}(t):=e^{4s\gamma(t)}\ell^{3/2}, \\ &&\mu_{1}(t):=\mu_{0}(t)\ell^{2}, \quad \mu_{k}(t):=e^{3s\gamma(t)}\ell^{(2k+9)/2} \quad \forall k\in \{2,\cdots,5\}.
	\end{eqnarray*}
\end{notation} 
\begin{remark} \label{Remark2}
	The following elementary estimates hold:
	\begin{eqnarray*}
		&& \mu_{3}\mu_{1}^{-2}=\mu^{-1}\ell^{2},\quad (\mu_{3}\mu_{1}^{-2})_{t}\leq C \mu^{-1},\quad |\mu_{3,t}|\leq C \mu_{1},\quad \mu_{0}\leq C\mu, \quad \mu\leq C\mu_{5}^{2},\\
		&& \quad \mu_{k}\leq C \mu_{k-1} \quad \forall k\in \{1,\cdots,5\},\quad  |\mu_{k}\mu_{k,t}|\leq C \mu_{k-1}^{2} \quad \forall k\in \{2,\cdots,5\}.
	\end{eqnarray*}
\end{remark}
As a consequence of Proposition \ref{P0}, we we establish the null controllability of \eqref{linearised cascade system} for source terms
which decay sufficiently fast to zero as $t\rightarrow 0^{+}$:
\begin{proposition} \label{P1} Let $F=(f,f_{\Gamma}), G=(g,g_{\Gamma})$ such that $\mu F, \mu G\in L^{2}(0,T;\mathbb{L}^{2})$.
	Then, there exists a control $v$ such that the solution $(\Psi,H)$ of \eqref{linearised cascade system} corresponding to $v$, $F$ and $G$, satisfies
		\begin{eqnarray}
			&&\|\mu_{0}\Psi\|_{L^{2}(0,T;\mathbb{L}^{2})}^{2} + \|\mu_{0}H\|_{L^{2}(0,T;\mathbb{L}^{2})}^{2}+ \|\mu_{1}v\|^{2}_{L^{2}(\omega_{T})}\nonumber\\
			&& \leq  C \left(\|\mu F\|_{L^{2}(0,T;\mathbb{L}^{2})}^{2}+  \|\mu G\|_{L^{2}(0,T;\mathbb{L}^{2})}^{2} \right). \label{c21}
	\end{eqnarray}
	In particular \eqref{linearised cascade system} is null controllable a time $t=0$. Moreover, we can choose $v$ satisfying 
	\begin{eqnarray}
		\mu_{3}v\in H^{1}(0,T;L^{2}(\omega))\cap L^{2}(0,T;H^{2}(\omega)) \label{c22}
	\end{eqnarray}
	and 
	\begin{eqnarray}
		\|\mu_{3}v_{t}\|^{2}_{L^{2}(\omega_{T})}\leq C \left(\|\mu F \|^{2}_{L^{2}(0,T;\mathbb{L}^{2})} +  \|\mu G \|^{2}_{L^{2}(0,T;\mathbb{L}^{2})} \right). \label{c41}
	\end{eqnarray}
\end{proposition}
	\begin{proof} 
	The proof of this result is inspired by the method of Fursikov and Imanuvilov  \cite{fursikov1996controllability}. Let us consider the following space:
	\begin{eqnarray*}
		\mathbf{P}&:=&\left\{(Y,Z)=((y,y_{\Gamma}),(z,z_{\Gamma}))\;:\; y,z\in C^{2}(\overline{\Omega_{T}})\;\;\mbox{and}\;\; y_{|_{\Gamma}}(\cdot,t)=y_{\Gamma}(\cdot,t), \right. \\
		&& \left. \;z_{|_{\Gamma}}(\cdot,t)=z_{\Gamma}(\cdot,t),\; t\in [0,T]\quad\mbox{and}\quad y(x,0)=z(x,T)=0, \; x\in \Omega \right\}.
	\end{eqnarray*}
	We define the bilinear form $\B: \mathbf{P}\times \mathbf{P}\longrightarrow\mathbb{R}$ by
	\begin{eqnarray*}
		&&\B((Y,Z),(\overline{Y},\overline{Z})) :=\int_{\Omega_{T}}\mu_{0}^{-2}(\textbf{L}_{1}^{\star}Y-z\mathds{1}_{\mathcal{O}}) (\textbf{L}^{\star}_{1}\overline{Y}-\overline{z}\mathds{1}_{\mathcal{O}}) \\
		&&+ \int_{\Gamma_{T}}\mu_{0}^{-2}(\textbf{L}_{2}^{\star}Y-z_{\Gamma}\mathds{1}_{\Sigma}) (\textbf{L}^{\star}_{2}\overline{Y}-\overline{z}_{\Gamma}\mathds{1}_{\Sigma})+ \int_{\Omega_{T}}\mu_{0}^{-2}\textbf{L}_{1}Z \textbf{L}_{1}\overline{Z} \\
		&& + \int_{\Gamma_{T}}\mu_{0}^{-2}\textbf{L}_{2}Z \textbf{L}_{2}\overline{Z}+ \int_{\Omega_{T}}\mu_{1}^{-2}\chi y\overline{y} ,
	\end{eqnarray*}
	where $\chi\in C_{0}^{\infty}(\omega)$ (the space of test functions with compact support in $\omega$) is given such that $0\leq \chi\leq 1$ and $\chi_{|_{\omega^{\prime\prime\prime}}}=1$; and
	$\textbf{L}_{1}$, $\textbf{L}_{2}$, $\textbf{L}_{1}^{\star}$ and $\textbf{L}_{2}^{\star}$ are defined in \eqref{n1}. We also define the linear form $\mathbf{F}: \mathbf{P}\longrightarrow\mathbb{R}$ by
	\begin{eqnarray*}
		\mathbf{F}(Y,Z):=\langle F, Y \rangle_{L^{2}(0,T;\mathbb{L}^{2})}+ \langle G, Z \rangle_{L^{2}(0,T;\mathbb{L}^{2})}.
	\end{eqnarray*}
	We claim that $\B$ is an inner
	product in $\mathbf{P}$ and $\mathbf{F}$ is  continuous for the norm $\|\cdot\|_{\mathbf{B}}$ associated with the scalar product $\mathbf{B}$. Indeed, due to Carleman estimate \eqref{estimateP0}, there exists a constant $C:=C(s,\lambda,T)>0$ such that for all $Y=(y,y_{\Gamma}), Z=(z,z_{\Gamma})\in \mathbf{P}$, one has
	\begin{eqnarray}
		&&\int_{\Omega_{T}}\mu^{-2}(|y|^{2}+|z|^{2}) + \int_{\Gamma_{T}}\mu^{-2}(|y_{\Gamma}|^{2}+|z_{\Gamma}|^{2}) \leq C\;\B((Y,Z), (Y,Z)). 
		\label{c9}
	\end{eqnarray}
	In particular, $\B$ is a scalar
	product in $\mathbf{P}$. To ensure the continuity of $\mathbf{F}$, using Cauchy-Schwarz inequality and  \eqref{c9}, we obtain
	\begin{eqnarray}
		|\mathbf{F}(Y,Z)|  
		\leq C \left(\|\mu F \|_{L^{2}(0,T;\mathbb{L}^{2})} + \|\mu G \|_{L^{2}(0,T;\mathbb{L}^{2})} \right) \|(Y,Z)\|_{\B}. \label{c10}
	\end{eqnarray}
	In the sequel, we will denote by $\overline{\mathbf{P}}$ the completion of $\mathbf{P}$ for the norm $\|\cdot\|_{\B}$ and we will still denote $\B$ and $\mathbf{F}$ the corresponding continuous extensions.
	From the Riesz Representation theorem,  there exists a unique
	$(\Phi,K)\in \overline{\mathbf{P}}$
	such that
	\begin{eqnarray}
		\mathbf{B}((\Phi,K), (Y,Z))=\mathbf{F}(Y,Z)\quad \forall (Y,Z)\in \overline{\mathbf{P}}. \label{c11}
	\end{eqnarray}
	Using \eqref{c10} and \eqref{c11}, we obtain
	\begin{eqnarray}
		\|(\Phi,K)\|_{\B}
		\leq   C \left(\|\mu F \|_{L^{2}(0,T;\mathbb{L}^{2})} +  \|\mu G \|_{L^{2}(0,T;\mathbb{L}^{2})} \right). \label{c15}
	\end{eqnarray}
	Let us introduce $(\Psi, H, v)$ with 
		\begin{eqnarray}
			&&\Psi:=\mu_{0}^{-2}(\mathbf{L}_{1}^{*}\Phi-k\mathds{1}_{\mathcal{O}}, \mathbf{L}_{2}^{*}\Phi-k_\Gamma\mathds{1}_{\Sigma}),\; H:=\mu_{0}^{-2}(\mathbf{L}_{1}K, \mathbf{L}_{2}K),\; v:=-\chi\mu_{1}^{-2}\phi |_{\omega_{T}}. \label{c16} 
	\end{eqnarray}
	According to \eqref{c16} and the definition of $\B$ and $\chi^{2}\leq\chi$, we obtain 
		\begin{eqnarray}
			\int_{\Omega_{T}}\mu_{0}^{2}(|\psi|^{2}+|h|^{2}) + \int_{\Gamma_{T}}\mu_{0}^{2}(|\psi_{\Gamma}|^{2}+|h_{\Gamma}|^{2}) + \int_{\omega_{T}}\mu_{1}^{2}|v|^{2} \leq \B((\Phi,K),(\Phi,K)).\nonumber 
	\end{eqnarray} 
	Using this estimate and \eqref{c15}, we deduce estimate \eqref{c21}. As a consequence $\Psi, H\in L^{2}(0,T;\mathbb{L}^{2})$, $v\in L^{2}(\omega_{T})$ and from \eqref{c11}, $(\Psi,H)$ is the unique distributional solution of \eqref{cascade quasi-linear system} associated with the control $v$.
	Let $(W,Z):=-\chi\mu_{1}^{-2}(\Phi,K)$ and  $(\tilde{f},\tilde{g}):=(\mathbf{L}^{\star}_{1}(\mu_{3}W)-\mu_{3}z\mathds{1}_{\mathcal{O}}, \mathbf{L}(\mu_{3}Z))$.  Then 
	\begin{eqnarray}
		\tilde{f} &=&-\chi\mu_{3}\mu_{1}^{-2}(\mathbf{L}^{\star}_{1}\Phi-k\mathds{1}_{\mathcal{O}}) +\chi(\mu_{3}\mu_{1}^{-2})_{t}\phi
		+2\sigma(0)\mu_{3}\mu_{1}^{-2}\nabla\chi\cdot\nabla\phi \nonumber\\
		&& +\sigma(0)\mu_{3}\mu_{1}^{-2}\Delta\chi\;\phi := \tilde{f}_{1}+\tilde{f}_{2}+\tilde{f}_{3}+\tilde{f}_{4} \label{c19}
	\end{eqnarray}
	and
	\begin{eqnarray}
		\tilde{g} &=&-\chi\mu_{3}\mu_{1}^{-2}\mathbf{L}_{1}K-\chi(\mu_{3}\mu_{1}^{-2})_{t}k
		+2\sigma(0)\mu_{3}\mu_{1}^{-2}\nabla\chi\cdot\nabla k \nonumber\\
		&& +\sigma(0)\mu_{3}\mu_{1}^{-2}\Delta\chi\;k := \tilde{g}_{1}+\tilde{g}_{2}+\tilde{g}_{3}+\tilde{g}_{4}. \label{**}
	\end{eqnarray}
	By the definition of $\psi$ and $h$ in \eqref{c16} and Remark \ref{Remark2}, we obtain 
	\begin{eqnarray}
		\begin{cases}
			|\tilde{f}_{1}|\leq C\mu_{0}|\psi|,\quad |\tilde{f}_{2}|\leq C\mu^{-1}|\phi|,\quad |\tilde{f}_{3}|\leq C\mu^{-1}\ell|\nabla\phi|,\quad |\tilde{f}_{4}|\leq C\mu^{-1}|\phi|\\
				|\tilde{g}_{1}|\leq C\mu_{0}|h|,\quad |\tilde{g}_{2}|\leq C\mu^{-1}|k|,\quad |\tilde{g}_{3}|\leq C\mu^{-1}\ell|\nabla k|,\quad |\tilde{g}_{4}|\leq C\mu^{-1}|k|. \label{c17}
		\end{cases}
	\end{eqnarray}
	From the Carleman estimate \eqref{estimateP0} and \eqref{c15}, we have
	\begin{eqnarray}
		&&\int_{\Omega_{T}}\mu^{-2}\left[\ell^{2}|\nabla \phi|^{2} + |\phi|^{2}\right] + \int_{\Omega_{T}}\mu^{-2}\left[\ell^{2}|\nabla k|^{2} + |k|^{2}\right] \nonumber\\
		&&\leq  C \left(\|\mu F \|^{2}_{L^{2}(0,T;\mathbb{L}^{2})} +  \|\mu G \|^{2}_{L^{2}(0,T;\mathbb{L}^{2})} \right). \label{c20}
	\end{eqnarray}
	Taking into account \eqref{c19}-\eqref{c20}, we obtain $\tilde{f}, \tilde{g}\in L^{2}(\Omega_{T})$. Moreover
	\begin{eqnarray}
		 \|\tilde{f}\|^{2}_{L^{2}(\Omega_{T})}+\|\tilde{g}\|^{2}_{L^{2}(\Omega_{T})}\leq  C \left(\|\mu F \|^{2}_{L^{2}(0,T;\mathbb{L}^{2})} +  \|\mu G \|^{2}_{L^{2}(0,T;\mathbb{L}^{2})} \right). \label{21}
	\end{eqnarray} 
	Since $\chi \in C_{0}^{\infty}(\omega)$ and $\omega\Subset \Omega$, then $(\mu_{3}W, \mu_{3}Z)$ is the strong solution of 
	\begin{equation*}
		\left\{
		\begin{aligned}
			&\mathbf{L}^{\star}_{1}(\mu_{3}W)-\mu_{3}z\mathds{1}_{\mathcal{O}} =\tilde{f} & & \text {in}\;\Omega_{T}, \\
			&\mathbf{L}_{1}(\mu_{3}Z)=\tilde{g} & & \text {in}\;\Omega_{T}, \\
			&\mathbf{L}^{\star}_{2}(\mu_{3}W)=0 & & \text {on}\;\Gamma_{T}, \\
			&\mathbf{L}_{2}(\mu_{3}Z)=0 & & \text {on}\;\Gamma_{T}, \\
			& (\mu_{3}w)_{\Gamma}= (\mu_{3}w)_{|_{\Gamma}},\; (\mu_{3}z)_{\Gamma}= (\mu_{3}z)_{|_{\Gamma}} & & \text {on}\;\Gamma_{T}, \\
			&(\mu_{3}w(\cdot,T),\mu_{3}w_{\Gamma}(\cdot,T))=(\mu_{3}z(\cdot,0),\mu_{3}z_{\Gamma}(\cdot,0))=(0,0) & & \text {in } \Omega\times\Gamma.
		\end{aligned}
		\right. 
	\end{equation*}
	Using estimates \eqref{energy estimate} and \eqref{21}, we obtain
	\begin{eqnarray}
		\|\mu_{3}W\|^{2}_{\mathfrak{E}_T}+\|\mu_{3}Z\|^{2}_{\mathfrak{E}_T}\leq C \left(\|\mu F \|^{2}_{L^{2}(0,T;\mathbb{L}^{2})} +  \|\mu G \|^{2}_{L^{2}(0,T;\mathbb{L}^{2})} \right). \label{c54}
	\end{eqnarray}
	Finally, since $\mu_{3}w|_{\omega}=\mu_{3}v$, $\mu_{3}v_{t}=(\mu_{3}v)_{t}-\mu_{3,t}v$ and  $|\mu_{3,t}|\leq C \mu_{1}$, then \eqref{c22} and \eqref{c41} are satisfied.
\end{proof}
\subsection{Additional estimates for the state of \eqref{linearised cascade system}}
	We provide the following additional estimates for the state of \eqref{linearised cascade system}, which will be essential for proving the surjectivity of the mapping $\varLambda$ defined in \eqref{varlamda}.
\begin{proposition} \label{P2}
	Under the assumptions of Proposition \ref{P1} and if $(\Psi, H, v)$ is a state-control for \eqref{linearised cascade system} provided by Proposition \ref{P1} associated with data $F=(f,f_{\Gamma})$ and $G=(g,g_{\Gamma})$ such that $\mu F, \mu G\in L^{2}(0,T;\mathbb{L}^{2})$ and $\mu_{4}F_{t}\in L^{2}(0,T;\mathbb{L}^{2})$. The following estimates are satisfied
	\begin{enumerate}[label=(\arabic*)]
		\item  
		\begin{eqnarray}
			&&\sup_{0\leq t \leq T}\mu_{2}^{2}(t)\|\Psi(\cdot,t)\|_{\mathbb{L}^{2}}^{2}
			+ \|\mu_{2}\nabla\Psi\|^{2}_{L^{2}(0,T;\mathbb{L}^{2})} +\sup_{0\leq t \leq T}\mu_{2}^{2}(t)\|H(\cdot,t)\|_{\mathbb{L}^{2}}^{2}\nonumber\\
			&& 
				+ \|\mu_{2}\nabla H\|^{2}_{L^{2}(0,T;\mathbb{L}^{2})}
			\leq C \left( \|\mu F\|^{2}_{L^{2}(0,T;\mathbb{L}^{2})}+ \|\mu G\|^{2}_{L^{2}(0,T;\mathbb{L}^{2})} \right).  \label{c25}
		\end{eqnarray}
		\item 
		\begin{eqnarray}
			&&\sup_{0\leq t \leq T}\mu_{3}^{2}(t)\|\nabla\Psi(\cdot,t)\|^{2}_{\mathbb{L}^{2}}
			+ \|\mu_{3}\Psi_{t}\|^{2}_{L^{2}(0,T;\mathbb{L}^{2})}+ \|\mu_{3}\Delta\Psi\|^{2}_{L^{2}(0,T;\mathbb{L}^{2})} \nonumber\\
			&& +\sup_{0\leq t \leq T}\mu_{3}^{2}(t)\|\nabla H(\cdot,t)\|^{2}_{\mathbb{L}^{2}}
				+ \|\mu_{3}H_{t}\|^{2}_{L^{2}(0,T;\mathbb{L}^{2})}+ \|\mu_{3}\Delta H\|^{2}_{L^{2}(0,T;\mathbb{L}^{2})} \nonumber\\
			&& \leq C \left( \|\mu F\|^{2}_{L^{2}(0,T;\mathbb{L}^{2})}+ \|\mu G\|^{2}_{L^{2}(0,T;\mathbb{L}^{2})} \right). \label{c26}
		\end{eqnarray}
		\item 
		\begin{eqnarray}
			&&\sup_{0\leq t \leq T}\mu_{4}^{2}(t)\|\Psi_{t}(\cdot,t)\|^{2}_{\mathbb{L}^{2}}
			+ \|\mu_{4}\nabla\Psi_{t}\|^{2}_{L^{2}(0,T;\mathbb{L}^{2})} \nonumber\\
			&& \leq C \left( \|\mu F\|^{2}_{L^{2}(0,T;\mathbb{L}^{2})}+\|\mu G\|^{2}_{L^{2}(0,T;\mathbb{L}^{2})}+\|\mu_{4}F_{t}\|^{2}_{L^{2}(0,T;\mathbb{L}^{2})} \right). 
			\label{c27}
		\end{eqnarray}
		\item 
		\begin{eqnarray}
			&&\sup_{0\leq t \leq T}\mu_{5}^{2}(t)\|\nabla\Psi_{t}(\cdot,t)\|^{2}_{\mathbb{L}^{2}}
			+ \|\mu_{5}\Psi_{tt}\|^{2}_{L^{2}(0,T;\mathbb{L}^{2})}+ \|\mu_{5}\Delta\Psi_{t}\|^{2}_{L^{2}(0,T;\mathbb{L}^{2})}\nonumber\\
			&& +\sup_{0\leq t \leq T}\mu_{5}^{2}(t)\|\Delta\Psi(\cdot,t)\|^{2}_{\mathbb{L}^{2}}
			+ \|\mu_{5}\nabla\Psi_{t}\|^{2}_{L^{2}(0,T;\mathbb{L}^{2})}\nonumber\\
			&&
			\leq C \left( \|\mu F\|^{2}_{L^{2}(0,T;\mathbb{L}^{2})}+\|\mu G\|^{2}_{L^{2}(0,T;\mathbb{L}^{2})}+\|\mu_{4}F_{t}\|^{2}_{L^{2}(0,T;\mathbb{L}^{2})} \right).  \label{c28}
		\end{eqnarray}
	\end{enumerate}
\end{proposition}
\begin{proof} 
	(1)  Firstly,  multiplying the first system of \eqref{linearised cascade system} by $\mu_{2}^{2}\psi$ and  integrating it in $\Omega$, one has
	\begin{eqnarray*}
		&&\frac{1}{2}\frac{\d}{\d t}\int_{\Omega}\mu_{2}^{2}|\psi|^{2} +\sigma(0)\int_{\Omega}\mu_{2}^{2}|\nabla \psi|^{2}= \int_{\Omega}\mu_{2}\mu_{2,t}|\psi|^{2}\nonumber\\
		&&  + \sigma(0)\int_{\Gamma}\mu_{2}^{2}\psi\partial_{\nu}\psi -a^{\prime}(0)\int_{\Omega}\mu_{2}^{2}|\psi|^{2} + \int_{\omega}\mu_{2}^{2} \psi v+ \int_{\Omega}\mu_{2}^{2} \psi f. 
	\end{eqnarray*}
	Secondly,
	multiplying the third system of \eqref{linearised cascade system} by $\mu_{2}^{2}\psi_{\Gamma}$ and integrating it on $\Gamma$, using the Stokes divergence formula \eqref{Stokes}, we obtain
	\begin{eqnarray*}
		&&\frac{1}{2}\frac{\d}{\d t}\int_{\Gamma}\mu_{2}^{2}|\psi_{\Gamma}|^{2} +\delta(0)\int_{\Gamma}\mu_{2}^{2}|\nabla_{\Gamma} \psi_{\Gamma}|^{2}= \int_{\Gamma}\mu_{2}\mu_{2,t}|\psi_{\Gamma}|^{2}  \nonumber\\
		&&- \sigma(0)\int_{\Gamma}\mu_{2}^{2}\psi_{\Gamma}\partial_{\nu}\psi-b^{\prime}(0)\int_{\Gamma}\mu_{2}^{2}|\psi_{\Gamma}|^{2} + \int_{\Gamma}\mu_{2}^{2} \psi_{\Gamma} f_{\Gamma}.
	\end{eqnarray*}
	Next, we add these identities, using the fact that $|\mu_{2}\mu_{2,t}|\leq C \mu_{0}^{2}$, $\mu_{2}\leq C \mu_{0}$
	and Young's inequality, we get
	\begin{eqnarray*}
		&&\frac{1}{2}\frac{\d}{\d t}\|\mu_{2}\Psi\|^{2}_{\mathbb{L}^{2}}
		+ \sigma(0)\int_{\Omega}\mu_{2}^{2}|\nabla \psi|^{2} + \delta(0)\int_{\Gamma}\mu_{2}^{2}|\nabla_{\Gamma} \psi_{\Gamma}|^{2} \nonumber\\	
		&& \leq C\|\mu_{0}\Psi\|_{\mathbb{L}^{2}}^{2}+C \left(\int_{\omega}\mu_{2}^{2}|v|^{2}+ \int_{\Omega}\mu_{2}^{2}|f|^{2} + \int_{\Gamma}\mu_{2}^{2}|f_{\Gamma}|^{2}\right).
	\end{eqnarray*}
	Integrating over $(0,t)$, one has
	\begin{eqnarray}
		&&\sup_{0\leq t \leq T}\mu_{2}^{2}(t)\|\Psi(\cdot,t)\|^{2}_{\mathbb{L}^{2}}
		+ \|\mu_{2}\nabla\Psi\|^{2}_{L^{2}(0,T;\mathbb{L}^{2})}\nonumber\\
		&&\leq C \left(\|\mu_{0}\Psi\|_{L^{2}(0,T;\mathbb{L}^{2})}^{2}+ \|\mu_{2}v\|^{2}_{L^{2}(\omega_{T})} + \|\mu_{2}F\|^{2}_{L^{2}(0,T;\mathbb{L}^{2})} \right). \nonumber
	\end{eqnarray}
	Using the fact that $\mu_{2}\leq C \mu_{1}$, $\mu_{2}\leq C \mu$ and estimate \eqref{c21}, we deduce 
	\begin{eqnarray}
		&&\sup_{0\leq t \leq T}\mu_{2}^{2}(t)\|\Psi(\cdot,t)\|^{2}_{\mathbb{L}^{2}}
		+ \|\mu_{2}\nabla\Psi\|^{2}_{L^{2}(0,T;\mathbb{L}^{2})} \leq C \left(\|\mu F\|^{2}_{L^{2}(0,T;\mathbb{L}^{2})} + \|\mu G\|^{2}_{L^{2}(0,T;\mathbb{L}^{2})} \right). \label{estimate1P34}
	\end{eqnarray}
	A computation similar to the proof of \eqref{estimate1P34}, by multiplying the second system of \eqref{linearised cascade system} by $\mu_{2}^{2}h$ and the fourth system of \eqref{linearised cascade system} by $\mu_{2}^{2}h_{\Gamma}$, leads to the conclusion 
		\begin{eqnarray}
			&&\sup_{0\leq t \leq T}\mu_{2}^{2}(t)\|H(\cdot,t)\|^{2}_{\mathbb{L}^{2}}
			+ \|\mu_{2}\nabla H\|^{2}_{L^{2}(0,T;\mathbb{L}^{2})}\nonumber\\
			&&\leq C \left(\|\mu_{0}H\|_{L^{2}(0,T;\mathbb{L}^{2})}^{2}+ \|\mu_{2}\Psi\|_{L^{2}(0,T;\mathbb{L}^{2})}^{2} + \|\mu_{2}G\|^{2}_{L^{2}(0,T;\mathbb{L}^{2})}  \right). \nonumber
		\end{eqnarray}
		Using the fact that $\mu_{2}\leq C \mu_{0}, \mu_{2}\leq C \mu$ and estimates \eqref{c21} and \eqref{estimate1P34}, we deduce \eqref{c25}.\\
	\par 
	\noindent (2) Firstly, multiplying the first system of \eqref{linearised cascade system} by $\mu_{3}^{2}\psi_{t}$ and  integrating it in $\Omega$, we obtain 
	\begin{eqnarray*}
		&&\int_{\Omega}\mu_{3}^{2}|\psi_{t}|^{2} + \frac{\sigma(0)}{2}\frac{\d}{\d t}\int_{\Omega}\mu_{3}^{2}|\nabla \psi|^{2}=\sigma(0)\int_{\Omega}\mu_{3}\mu_{3,t}|\nabla \psi|^{2} \nonumber\\
		&& +\sigma(0)\int_{\Gamma}\mu_{3}^{2}\psi_{t}\partial_{\nu} \psi -a^{\prime}(0)\int_{\Omega}\mu_{3}^{2}\psi\psi_{t} + \int_{\omega}\mu_{3}^{2}\psi_{t} v+ \int_{\Omega}\mu_{3}^{2}\psi_{t} f. 
	\end{eqnarray*}
	Using Young's inequality and the following elementary estimates
	\begin{eqnarray}
		|\mu_{3}\mu_{3,t}|\leq C \mu_{2}^{2},\quad \quad \mu_{3}\leq C\mu_{0}, \label{c30}
	\end{eqnarray}
	we get 
	\begin{eqnarray}
		&& \frac{1}{2}\int_{\Omega}\mu_{3}^{2}|\psi_{t}|^{2} + \frac{\sigma(0)}{2}\frac{\d}{\d t}\int_{\Omega}\mu_{3}^{2}|\nabla \psi|^{2}
		\leq C\left(\int_{\Omega}\mu_{2}^{2}|\nabla \psi|^{2}  + \int_{\Omega}\mu_{0}^{2}|\psi|^{2}\right) \label{c31} \\
		&& 
		 +\sigma(0)\int_{\Gamma}\mu_{3}^{2}\psi_{t} \partial_{\nu} \psi +C \left(\int_{\omega}\mu_{3}^{2}|v|^{2} + \int_{\Omega}\mu_{3}^{2}|f|^{2} \right). \nonumber
	\end{eqnarray}
	On the other hand, multiplying the third system of \eqref{linearised cascade system} by $\mu_{3}^{2}\psi_{\Gamma,t}$ and integrating it on $\Gamma$, by the Stokes divergence formula \eqref{Stokes}, we find
	\begin{eqnarray*}
		&&\int_{\Gamma}\mu_{3}^{2}|\psi_{\Gamma,t}|^{2} + \frac{\delta(0)}{2}\frac{\d}{\d t}\int_{\Gamma}\mu_{3}^{2}|\nabla_{\Gamma} \psi_{\Gamma}|^{2}=\delta(0)\int_{\Gamma}\mu_{3}\mu_{3,t}|\nabla_{\Gamma} \psi_{\Gamma}|^{2} \nonumber\\
		&&-\sigma(0)\int_{\Gamma}\mu_{3}^{2}\psi_{\Gamma,t}\partial_{\nu} \psi -b^{\prime}(0)\int_{\Gamma}\mu_{3}^{2}\psi_{\Gamma}\psi_{\Gamma,t} + \int_{\Gamma}\mu_{3}^{2}\psi_{\Gamma,t} f_{\Gamma}.
	\end{eqnarray*}
	Using Young's inequality and \eqref{c30}, we have
	\begin{eqnarray}
		&& \frac{1}{2}\int_{\Gamma}\mu_{3}^{2}|\psi_{\Gamma,t}|^{2} + \frac{\delta(0)}{2}\frac{\d}{\d t}\int_{\Gamma}\mu_{3}^{2}|\nabla_{\Gamma} \psi_{\Gamma}|^{2}\leq C\left(\int_{\Gamma}\mu_{2}^{2}|\nabla_{\Gamma} \psi_{\Gamma}|^{2} + \int_{\Gamma}\mu_{0}^{2}| \psi_{\Gamma}|^{2}\right) \nonumber\\
		&&  -\sigma(0)\int_{\Gamma}\mu_{3}^{2}\psi_{\Gamma,t}\partial_{\nu} \psi+ C\int_{\Gamma}\mu_{3}^{2}|f_{\Gamma}|^{2}. \label{c32}
	\end{eqnarray}
	By summing \eqref{c31} and \eqref{c32}, integrating over $(0,t)$, using $\mu_{3}\leq C\mu_{1}$, $\mu_{3}\leq C\mu$ and estimates \eqref{c21} and \eqref{c25}, one has
	\begin{eqnarray}
		&&\|\mu_{3}\Psi_{t}\|^{2}_{L^{2}(0,T;\mathbb{L}^{2})}+ \sup_{0\leq t \leq T}\mu_{3}^{2}(t)\|\nabla\Psi(\cdot,t)\|^{2}_{\mathbb{L}^{2}}\leq C\left(\|\mu F\|^{2}_{L^{2}(0,T;\mathbb{L}^{2})}+ \|\mu G\|^{2}_{L^{2}(0,T;\mathbb{L}^{2})}\right). \label{c35}
	\end{eqnarray}
	\par 
	Now, multiplying the first system of \eqref{linearised cascade system} by $-\mu_{3}^{2}\Delta \psi$ and integrating it in $\Omega$, we obtain
	\begin{eqnarray*}
		&&\frac{1}{2}\frac{\d}{\d t}\int_{\Omega}\mu_{3}^{2}|\nabla \psi|^{2} + \sigma(0)\int_{\Omega}\mu_{3}^{2}|\Delta \psi|^{2} =\int_{\Omega}\mu_{3}\mu_{3,t}|\nabla \psi|^{2}\nonumber\\
		&& +\int_{\Gamma}\mu_{3}^{2}\psi_{t}\partial_{\nu} \psi +a^{\prime}(0)\int_{\Omega}\mu_{3}^{2}\Delta \psi\psi  - \int_{\omega}\mu_{3}^{2}\Delta \psi v- \int_{\Omega}\mu_{3}^{2}\Delta \psi f.
	\end{eqnarray*}
	Using $|\mu_{3}\mu_{3,t}|\leq C \mu_{2}$, 
	Young's inequality and the continuity of the normal derivative from $H^{2}(\Omega)$ to $L^{2}(\Gamma)$, we find
	\begin{eqnarray*}
		&&\frac{1}{2}\frac{\d}{\d t}\int_{\Omega}\mu_{3}^{2}|\nabla \psi|^{2} + \frac{\sigma(0)}{2} \int_{\Omega}\mu_{3}^{2}|\Delta \psi|^{2} \leq C\left( \int_{\Omega}\mu_{2}^{2}|\nabla\psi|^{2}   \right. \nonumber \\
		&& \left. 
		+ \int_{\Gamma}\mu_{3}^{2}|\psi_{\Gamma,t}|^{2} + \int_{\Omega}\mu_{0}^{2}|\psi|^{2} +\int_{\omega}\mu_{1}^{2}|v|^{2}+ \int_{\Omega}\mu^{2}|f|^{2}\right). 
	\end{eqnarray*}
	Integrating over $(0,t)$ and using estimates \eqref{c21}, \eqref{c25} and \eqref{c35}, we obtain 
	\begin{eqnarray}
		&&\sup_{0\leq t \leq T}\mu_{3}^{2}(t)\|\nabla \psi(\cdot,t)\|^{2}_{L^{2}(\Omega)} + \|\mu_{3}\Delta\psi\|^{2}_{L^{2}(\Omega_{T})} \leq C \left( \|\mu F\|^{2}_{L^{2}(0,T;\mathbb{L}^{2})}+\|\mu G\|^{2}_{L^{2}(0,T;\mathbb{L}^{2})} \right). \label{c33}
	\end{eqnarray}  
	On the other hand, multiplying the second system of \eqref{linearised cascade system} by $-\mu_{3}^{2}\Delta_{\Gamma} \psi_{\Gamma}$ and integrating it on $\Gamma$, we get
	\begin{eqnarray*}
		&&\frac{1}{2}\frac{\d}{\d t}\int_{\Gamma}\mu_{3}^{2}|\nabla_{\Gamma} \psi_{\Gamma}|^{2} + \delta(0)\int_{\Gamma}\mu_{3}^{2}|\Delta_{\Gamma} \psi_{\Gamma}|^{2}=\int_{\Gamma}\mu_{3}\mu_{3,t}|\nabla_{\Gamma} \psi_{\Gamma}|^{2}\nonumber\\
		&& +\sigma(0)\int_{\Gamma}\mu_{3}^{2}\Delta_{\Gamma} \psi_{\Gamma}\partial_{\nu} \psi +b^{\prime}(0)\int_{\Gamma}\mu_{3}^{2}\Delta_{\Gamma} \psi_{\Gamma}\psi_{\Gamma}  - \int_{\Gamma}\mu_{3}^{2}\Delta_{\Gamma} \psi_{\Gamma}f_{\Gamma}. 
	\end{eqnarray*}
	Using  
	Young's inequality and the continuity of the normal derivative from $H^{2}(\Omega)$ to $L^{2}(\Gamma)$, one has
	\begin{eqnarray*}
		&&\frac{1}{2}\frac{\d}{\d t}\int_{\Gamma}\mu_{3}^{2}|\nabla_{\Gamma} \psi_{\Gamma}|^{2} + \frac{\delta(0)}{2}\int_{\Gamma}\mu_{3}^{2}|\Delta_{\Gamma} \psi_{\Gamma}|^{2}\leq C\left(\int_{\Gamma}\mu_{2}^{2}|\nabla_{\Gamma} \psi_{\Gamma}|^{2} \right. \nonumber\\
		&& \left. +\int_{\Omega}\mu_{3}^{2}|\Delta\psi|^{2}+ \int_{\Omega}\mu_{3}^{2}|\psi|^{2} +\int_{\Gamma}\mu_{3}^{2} |\psi_{\Gamma}|^{2}  +\int_{\Gamma}\mu_{3}^{2}|f_{\Gamma}|^{2}\right). 
	\end{eqnarray*}
	Integrating over $(0,t)$ and using estimates \eqref{c25} and \eqref{c33}, we get
	\begin{eqnarray}
		&&\sup_{0\leq t \leq T}\mu_{3}^{2}(t)\|\nabla_{\Gamma} \psi_{\Gamma}(\cdot,t)\|^{2}_{L^{2}(\Gamma)} + \|\mu_{3}\Delta_{\Gamma}\psi_{\Gamma}\|^{2}_{L^{2}(\Gamma_{T})} \leq  C \left( \|\mu F\|^{2}_{L^{2}(0,T;\mathbb{L}^{2})}+\|\mu G\|^{2}_{L^{2}(0,T;\mathbb{L}^{2})}\right). \nonumber \\ \label{c34}
	\end{eqnarray}
	Based on \eqref{c35}-\eqref{c34}, we obtain	\begin{eqnarray}
		&&\sup_{0\leq t \leq T}\mu_{3}^{2}(t)\|\nabla\Psi(\cdot,t)\|^{2}_{\mathbb{L}^{2}}
		+ \|\mu_{3}\Psi_{t}\|^{2}_{L^{2}(0,T;\mathbb{L}^{2})}+ \|\mu_{3}\Delta\Psi\|^{2}_{L^{2}(0,T;\mathbb{L}^{2})} \nonumber\\
		&& \leq C \left( \|\mu F\|^{2}_{L^{2}(0,T;\mathbb{L}^{2})}+ \|\mu G\|^{2}_{L^{2}(0,T;\mathbb{L}^{2})} \right). \label{c26forPsi}
	\end{eqnarray}
	Using the same multiplication techniques as in \eqref{c26forPsi}, we arrive at
	\begin{eqnarray}
		&&\sup_{0\leq t \leq T}\mu_{3}^{2}(t)\|\nabla H(\cdot,t)\|^{2}_{\mathbb{L}^{2}}
		+ \|\mu_{3}H_{t}\|^{2}_{L^{2}(0,T;\mathbb{L}^{2})}+ \|\mu_{3}\Delta H\|^{2}_{L^{2}(0,T;\mathbb{L}^{2})} \nonumber\\
		&& \leq C \left( \|\mu F\|^{2}_{L^{2}(0,T;\mathbb{L}^{2})}+ \|\mu G\|^{2}_{L^{2}(0,T;\mathbb{L}^{2})} \right). \label{c26forH}
	\end{eqnarray}
	Finally, from \eqref{c26forPsi} and  \eqref{c26forH}, we deduce \eqref{c26}.\\
	\par 
	\noindent (3) Differentiating with respect to
	time the system \eqref{linearised cascade system}, one has
	\begin{equation} \label{s}
		\left\{
		\begin{aligned}
			&\psi_{tt}-\sigma(0)\Delta\psi_{t} +a^{\prime}(0)\psi_{t} =f_{t} +
			v_{t}\mathds{1}_{\omega} & & \text {in}\; \Omega_T, \\
			&\psi_{\Gamma,tt}-\delta(0)\Delta_{\Gamma}\psi_{\Gamma,t}+\sigma(0)
			\partial_{\nu} \psi_{t} + b^{\prime}(0)\psi_{\Gamma,t}=f_{\Gamma,t} & & \text {on}\;\Gamma_T, \\
			& \psi_{\Gamma,t}=\psi_{t}|_{\Gamma},   & & \text {on}\;\Gamma_T, \\
			& (\psi_{t}(\cdot,0),\psi_{\Gamma,t}(\cdot,0))=(0,0) & & \text {in } \Omega\times\Gamma, 
		\end{aligned}
		\right.
	\end{equation}
	Multiplying the first equation of \eqref{s} by $\mu_{4}^{2}\psi_{t}$ and integrating it in $\Omega$, we obtain
	\begin{eqnarray*}
		&&\frac{1}{2}\frac{\d}{\d t}\int_{\Omega}\mu_{4}^{2}|\psi_{t}|^{2} + \sigma(0)\int_{\Omega}\mu_{4}^{2}|\nabla \psi_{t}|^{2}=\int_{\Omega}\mu_{4}\mu_{4,t}|\psi_{t}|^{2} \nonumber\\
		&&  +\sigma(0)\int_{\Gamma}\mu_{4}^{2}\psi_{t}\partial_{\nu}\psi_{t} -a^{\prime}(0)\int_{\Omega}\mu_{4}^{2}|\psi_{t}|^{2} + \int_{\omega}\mu_{4}^{2}\psi_{t} v_{t}+ \int_{\Omega}\mu_{4}^{2}\psi_{t}f_{t}.
	\end{eqnarray*}
	Using Young's inequality, $|\mu_{4}\mu_{4,t}|\leq C\mu_{3}^{2}$ and $\mu_{4}^{2}\leq C\mu_{3}^{2}$, we get
	\begin{eqnarray}
		&&\frac{1}{2}\frac{\d}{\d t}\int_{\Omega}\mu_{4}^{2}|\psi_{t}|^{2} + \sigma(0)\int_{\Omega}\mu_{4}^{2}|\nabla \psi_{t}|^{2}\leq C\int_{\Omega}\mu_{3}^{2}|\psi_{t}|^{2}\nonumber\\
		&&  +\sigma(0)\int_{\Gamma}\mu_{4}^{2}\psi_{t}\partial_{\nu}\psi_{t} + C\int_{\omega}\mu_{3}^{2}|v_{t}|^{2}+ C\int_{\Omega}\mu_{4}^{2}|f_{t}|^{2}. \label{c36}
	\end{eqnarray}
	Multiplying the second equation of \eqref{s} by $\mu_{4}^{2}\psi_{\Gamma,t}$ and integrating it on $\Gamma$, we get
	\begin{eqnarray} 
		&&\frac{1}{2}\frac{\d}{\d t}\int_{\Gamma}\mu_{4}^{2}|\psi_{\Gamma,t}|^{2} + \delta(0)\int_{\Gamma}\mu_{4}^{2}|\nabla_{\Gamma} \psi_{\Gamma,t}|^{2}=\int_{\Gamma}\mu_{4}\mu_{4,t}|\psi_{\Gamma,t}|^{2} \nonumber\\
		&& -\sigma(0)\int_{\Gamma}\mu_{4}^{2} \psi_{\Gamma,t}\partial_{\nu}\psi_{t}- b^{\prime}(0)\int_{\Gamma}\mu_{4}^{2}|\psi_{\Gamma,t}|^{2} + \int_{\Gamma}\mu_{4}^{2}\psi_{\Gamma,t} f_{\Gamma,t}. \nonumber\\
		&&\leq C\int_{\Gamma}\mu_{3}^{2}|\psi_{\Gamma,t}|^{2}-\sigma(0)\int_{\Gamma}\mu_{4}^{2} \psi_{\Gamma,t}\partial_{\nu}\psi_{t} + C\int_{\Gamma}\mu_{4}^{2}|f_{\Gamma,t}|^{2}. \label{c37}
	\end{eqnarray}
	By summing \eqref{c36} and \eqref{c37} and integrating the estimate obtained  over $(0,t)$, we find 
	\begin{eqnarray}
		&& \sup_{0\leq t \leq T}\mu_{4}^{2}(t)\|\Psi_{t}(\cdot,t)\|^{2}_{\mathbb{L}^{2}}+ \|\mu_{4}\nabla\Psi_{t}\|^{2}_{L^{2}(0,T;\mathbb{L}^{2})} \nonumber\\
		&&\leq C \left(\|\mu_{3}\Psi_{t}\|^{2}_{L^{2}(0,T;\mathbb{L}^{2})}+\|\mu_{3} v_{t}\|^{2}_{L^{2}(\omega_{T})} +\|\mu_{4} F_{t}\|^{2}_{L^{2}(0,T;\mathbb{L}^{2})} \right). \label{c38}
	\end{eqnarray}
	Using estimates \eqref{c38},  \eqref{c26} and \eqref{c41}, we obtain \eqref{c27}.\\
	\par 
	\noindent (4) Firstly, multiplying the first equation of \eqref{s} by $\mu_{5}^{2}\psi_{tt}$ and integrating it in $\Omega$, we obtain
	\begin{eqnarray*}
		&&\int_{\Omega}\mu_{5}^{2}|\psi_{tt}|^{2} + \frac{\sigma(0)}{2}\frac{\d}{\d t}\int_{\Omega}\mu_{5}^{2}|\nabla \psi_{t}|^{2}=\sigma(0)\int_{\Omega}\mu_{5}\mu_{5,t}|\nabla \psi_{t}|^{2}\nonumber\\
		&& +\sigma(0)\int_{\Gamma}\mu_{5}^{2}\psi_{tt}\partial_{\nu} \psi_{t} -a^{\prime}(0)\int_{\Omega}\mu_{5}^{2}\psi_{t}\psi_{tt} + \int_{\omega}\mu_{5}^{2}\psi_{tt} v_{t}+ \int_{\Omega}\mu_{5}^{2}\psi_{tt} f_{t}. 
	\end{eqnarray*}
	Using Young's inequality and the following elementary estimates
	\begin{eqnarray*}
		|\mu_{5}\mu_{5,t}|\leq C \mu_{4}^{2}, \quad \quad \mu_{k+1}\leq C\mu_{k}, \label{c}
	\end{eqnarray*}
	we get 
	\begin{eqnarray}
		&& \frac{1}{2}\int_{\Omega}\mu_{5}^{2}|\psi_{tt}|^{2} + \frac{\sigma(0)}{2}\frac{\d}{\d t}\int_{\Omega}\mu_{5}^{2}|\nabla \psi_{t}|^{2}\leq C\int_{\Omega}\mu_{4}^{2}|\nabla \psi_{t}|^{2}  \nonumber\\
		&& +\sigma(0)\int_{\Gamma}\mu_{5}^{2}\psi_{tt}\partial_{\nu} \psi_{t}+C\int_{\Omega}\mu_{4}^{2}|\psi_{t}|^{2} + C\int_{\omega}\mu_{3}^{2}|v_{t}|^{2}\d x + \int_{\Omega}\mu_{4}^{2}|f_{t}|^{2}. \label{c42}
	\end{eqnarray} 
	Secondly, multiplying the second equation of \eqref{s} by $\mu_{5}^{2}\psi_{\Gamma,tt}$ and integrating over $\Gamma$, we find
	\begin{eqnarray*}
		&&\int_{\Gamma}\mu_{5}^{2}|\psi_{\Gamma,tt}|^{2} + \frac{\delta(0)}{2}\frac{\d}{\d t}\int_{\Gamma}\mu_{5}^{2}|\nabla_{\Gamma} \psi_{\Gamma,t}|^{2}=\delta(0)\int_{\Gamma}\mu_{5}\mu_{5,t}|\nabla_{\Gamma} \psi_{\Gamma,t}|^{2} \nonumber\\
		&&-\sigma(0)\int_{\Gamma}\mu_{5}^{2}\psi_{\Gamma,tt}\partial_{\nu} \psi_{t} -b^{\prime}(0)\int_{\Gamma}\mu_{5}^{2}\psi_{\Gamma,t}\psi_{\Gamma,tt} + \int_{\Gamma}\mu_{5}^{2}\psi_{\Gamma,tt} f_{\Gamma,t}.
	\end{eqnarray*}
	Using Young's inequality, we have
	\begin{eqnarray}
		&& \frac{1}{2}\int_{\Gamma}\mu_{5}^{2}|\psi_{\Gamma,tt}|^{2} + \frac{\delta(0)}{2}\frac{\d}{\d t}\int_{\Gamma}\mu_{5}^{2}|\nabla_{\Gamma} \psi_{\Gamma,t}|^{2}\leq C\int_{\Gamma}\mu_{4}^{2}|\nabla_{\Gamma} \psi_{\Gamma,t}|^{2} \nonumber\\
		&&-\sigma(0)\int_{\Gamma}\mu_{5}^{2}\psi_{\Gamma,tt}\partial_{\nu} \psi_{t} +C\int_{\Gamma}\mu_{4}^{2}|\psi_{\Gamma,t}|^{2} + C\int_{\Gamma}\mu_{4}^{2}|f_{\Gamma,t}|^{2}. \label{c43}
	\end{eqnarray}
	By summing \eqref{c42} and \eqref{c43} and  integrating over $(0,t)$, one has
	\begin{eqnarray}
		&&\|\mu_{5}\Psi_{tt}\|^{2}_{L^{2}(0,T;\mathbb{L}^{2})}+ \sup_{0\leq t \leq T}\mu_{5}^{2}(t)\|\nabla\Psi_{t}(\cdot,t)\|^{2}_{\mathbb{L}^{2}}\leq C \left(\sup_{0\leq t \leq T}\mu_{4}^{2}(t)\|\Psi_{t}(\cdot,t)\|^{2}_{\mathbb{L}^{2}} \right. \nonumber\\
		&&  \left. +\|\mu_{4}\nabla\Psi_{t}\|^{2}_{L^{2}(0,T;\mathbb{L}^{2})} + \|\mu_{3}v_{t}\|^{2}_{L^{2}(\omega_{T})}+\|\mu_{4} F_{t}\|^{2}_{L^{2}(0,T;\mathbb{L}^{2})}  \right). \label{c44}
	\end{eqnarray} 
	Based on \eqref{c44} and estimates \eqref{c27}, \eqref{c41}, we obtain 
	\begin{eqnarray}
		&&\|\mu_{5}\Psi_{tt}\|^{2}_{L^{2}(0,T;\mathbb{L}^{2})}+ \sup_{0\leq t \leq T}\mu_{5}^{2}(t)\|\nabla\Psi_{t}(\cdot,t)\|^{2}_{\mathbb{L}^{2}} \nonumber\\
		&& \leq C\left(\|\mu F\|^{2}_{L^{2}(0,T;\mathbb{L}^{2})}+\|\mu G\|^{2}_{L^{2}(0,T;\mathbb{L}^{2})}+\|\mu_{4} F_{t}\|^{2}_{L^{2}(0,T;\mathbb{L}^{2})}\right). \label{c46}
	\end{eqnarray}
	Multiplying the first equation of \eqref{s} by $-\mu_{5}^{2}\Delta \psi_{t}$ and integrating it in $\Omega$, we obtain
	\begin{eqnarray*}
		&&\frac{1}{2}\frac{\d}{\d t}\int_{\Omega}\mu_{5}^{2}|\nabla \psi_{t}|^{2} + \sigma(0)\int_{\Omega}\mu_{5}^{2}|\Delta \psi_{t}|^{2}=\int_{\Omega}\mu_{5}\mu_{5,t}|\nabla \psi_{t}|^{2}  \nonumber\\
		&& +\int_{\Gamma}\mu_{5}^{2}\psi_{tt}\partial_{\nu} \psi_{t} +a^{\prime}(0)\int_{\Omega}\mu_{5}^{2}\Delta \psi_{t}\psi_{t}  - \int_{\omega}\mu_{5}^{2}\Delta \psi_{t} v_{t}- \int_{\Omega}\mu_{5}^{2}\Delta \psi_{t} f_{t}.
	\end{eqnarray*}
	Using Young's inequality and the continuity of the normal derivative from $H^{2}(\Omega)$ to $L^{2}(\Gamma)$, we find
	\begin{eqnarray*}
		&&\frac{1}{2}\frac{\d}{\d t}\int_{\Omega}\mu_{5}^{2}|\nabla \psi_{t}|^{2} + \frac{\sigma(0)}{2}\int_{\Omega}\mu_{5}^{2}|\Delta \psi_{t}|^{2}\leq C \left(\int_{\Omega}\mu_{4}^{2}|\nabla \psi_{t}|^{2} \right. \nonumber\\
		&& \left.
		+\int_{\Gamma}\mu_{5}^{2}|\psi_{\Gamma,tt}|^{2}+ \int_{\Omega}\mu_{4}^{2}|\psi_{t}|^{2}+ \int_{\omega}\mu_{3}^{2}|v_{t}|^{2} + \int_{\Omega}\mu_{4}^{2}|f_{t}|^{2}\right).
	\end{eqnarray*}
	Integrating over $(0,t)$, we obtain
	\begin{eqnarray*}
		&&\sup_{0\leq t \leq T}\int_{\Omega}\mu_{5}^{2}|\nabla \psi_{t}(\cdot,t)|^{2} + \int_{\Omega_{T}}\mu_{5}^{2}|\Delta \psi_{t}|^{2} \leq C\left(\int_{\Omega_{T}}\mu_{4}^{2}|\nabla \psi_{t}|^{2} \right.\\
		&& \left. + \sup_{0\leq t \leq T}\int_{\Omega}\mu_{4}^{2}|\psi_{t}|^{2}+ \int_{\Gamma_{T}}\mu_{5}^{2}|\psi_{\Gamma,tt}|^{2} 
		+ \int_{\omega_{T}}\mu_{3}^{2}|v_{t}|^{2} + \int_{\Omega_{T}}\mu_{4}^{2}|f_{t}|^{2}  \right). 
	\end{eqnarray*}
	Using estimates \eqref{c27} and  \eqref{c46}, we get
	\begin{eqnarray}
		&& \sup_{0\leq t \leq T}\mu_{5}^{2}(t)\|\nabla \psi_{t}(\cdot,t)\|^{2}_{L^{2}(\Omega)}+ \|\mu_{5}\Delta \psi_{t}\|^{2}_{L^{2}(\Omega_{T})} \nonumber\\
		&&\leq C\left(\|\mu F\|^{2}_{L^{2}(0,T;\mathbb{L}^{2})}+\|\mu G\|^{2}_{L^{2}(0,T;\mathbb{L}^{2})}+\|\mu_{4} F_{t}\|^{2}_{L^{2}(0,T;\mathbb{L}^{2})}\right).  \label{c47}
	\end{eqnarray}
	Now, multiplying the second equation of \eqref{s} by $-\mu_{5}^{2}\Delta_{\Gamma} \psi_{\Gamma,t}$ and integrating it on $\Gamma$, we get
	\begin{eqnarray*}
		&&\frac{1}{2}\frac{\d}{\d t}\int_{\Gamma}\mu_{5}^{2}|\nabla_{\Gamma} \psi_{\Gamma,t}|^{2} + \delta(0)\int_{\Gamma}\mu_{5}^{2}|\Delta_{\Gamma} \psi_{\Gamma,t}|^{2}=\int_{\Gamma}\mu_{5}\mu_{5,t}|\nabla_{\Gamma} \psi_{\Gamma,t}|^{2}  \nonumber\\
		&& +\sigma(0)\int_{\Gamma}\mu_{5}^{2}\Delta_{\Gamma} \psi_{\Gamma,t}\partial_{\nu} \psi_{t} +b^{\prime}(0)\int_{\Gamma}\mu_{5}^{2}\Delta_{\Gamma} \psi_{\Gamma,t}\psi_{\Gamma,t}  - \int_{\Gamma}\mu_{5}^{2}\Delta_{\Gamma} \psi_{\Gamma,t}f_{\Gamma,t}. 
	\end{eqnarray*}
	Similarly as above, one has
	\begin{eqnarray*}
		&&\frac{1}{2}\frac{\d}{\d t}\int_{\Gamma}\mu_{5}^{2}|\nabla_{\Gamma} \psi_{\Gamma,t}|^{2} +\frac{\delta(0)}{2} \int_{\Gamma}\mu_{5}^{2}|\Delta_{\Gamma} \psi_{\Gamma,t}|^{2} \leq C\left(\int_{\Gamma}\mu_{4}^{2}|\nabla_{\Gamma} \psi_{\Gamma,t}|^{2} \right. \nonumber\\
		&& \left. +\int_{\Omega}\mu_{5}^{2}|\psi_{t}|^{2}+ \int_{\Omega}\mu_{5}^{2}|\Delta\psi_{t}|^{2} +\int_{\Gamma}\mu_{5}^{2}|\psi_{\Gamma,t}|^{2}  +\int_{\Gamma}\mu_{4}^{2}|f_{\Gamma,t}|^{2} \right). 
	\end{eqnarray*}
	Hence 
	\begin{eqnarray}
		&& \sup_{0\leq t \leq T}\mu_{5}^{2}(t)\|\nabla_{\Gamma} \psi_{\Gamma,t}(\cdot,t)\|^{2}_{L^{2}(\Gamma)}+ \|\mu_{5}\Delta_{\Gamma} \psi_{\Gamma,t}\|^{2}_{L^{2}(\Gamma_{T})} \nonumber\\
		&&\leq C\left(\|\mu F\|^{2}_{L^{2}(0,T;\mathbb{L}^{2})}+\|\mu G\|^{2}_{L^{2}(0,T;\mathbb{L}^{2})}+\|\mu_{4} F_{t}\|^{2}_{L^{2}(0,T;\mathbb{L}^{2})} \right).  \label{c48}
	\end{eqnarray}
	From \eqref{c47} and \eqref{c48}, we obtain 
	\begin{eqnarray}
		&&\sup_{0\leq t \leq T}\mu_{5}^{2}(t)\|\nabla\Psi_{t}(\cdot,t)|^{2}_{\mathbb{L}^{2}}
		+ \|\mu_{5}\Delta \Psi_{t}\|^{2}_{L^{2}(0,T;\mathbb{L}^{2})} \nonumber\\
		&&\leq C\left(\|\mu F\|^{2}_{L^{2}(0,T;\mathbb{L}^{2})}+\|\mu G\|^{2}_{L^{2}(0,T;\mathbb{L}^{2})}\|\mu_{4} F_{t}\|^{2}_{L^{2}(0,T;\mathbb{L}^{2})}\right).  \label{c49}
	\end{eqnarray}
	Multiplying the first equation of \eqref{linearised cascade system} by $-\mu_{5}^{2}\Delta \psi_{t}$ and integrating it in $\Omega$, we obtain
	\begin{eqnarray*}
		&& \int_{\Omega}\mu_{5}^{2}|\nabla \psi_{t}|^{2}+\frac{\sigma(0)}{2}\frac{\d}{\d t}\int_{\Omega}\mu_{5}^{2}|\Delta \psi|^{2} =\int_{\Gamma}\mu_{5}^{2}\psi_{t}\partial_{\nu}\psi_{t} +\sigma(0)\int_{\Omega}\mu_{5}\mu_{5,t}|\Delta \psi|^{2} \nonumber\\
		&&  +a^{\prime}(0)\int_{\Omega}\mu_{5}^{2}\psi\Delta \psi_{t} - \int_{\omega}\mu_{5}^{2}\Delta \psi_{t} v- \int_{\Omega}\mu_{5}^{2}\Delta \psi_{t}f.
	\end{eqnarray*}
	Using Young's inequality and the trace theorem for the normal derivative, we have 
	\begin{eqnarray}
		&& \int_{\Omega}\mu_{5}^{2}|\nabla \psi_{t}|^{2}+\frac{\sigma(0)}{2}\frac{\d}{\d t}\int_{\Omega}\mu_{5}^{2}|\Delta \psi|^{2} \leq C\left(\int_{\Gamma}\mu_{3}^{2}|\psi_{t}|^{2} + \int_{\Omega}\mu_{3}^{2}|\psi_{t}|^{2} \right.\nonumber\\
		&&  \left. + \int_{\Omega}\mu_{5}^{2}|\Delta\psi_{t}|^{2}+\int_{\Omega}\mu_{3}^{2}|\Delta \psi|^{2}+\int_{\Omega}\mu_{0}^{2}|\psi|^{2}  +\int_{\omega}\mu_{1}^{2}|v|^{2} + \int_{\Omega}\mu^{2}|f|^{2}\right).  \label{cc50}
	\end{eqnarray} 
	Multiplying the third equation of \eqref{linearised cascade system} by $-\mu_{5}^{2}\Delta_{\Gamma} \psi_{\Gamma,t}$ and integrating it in $\Gamma$, we obtain
	\begin{eqnarray*}
		&& \int_{\Gamma}\mu_{5}^{2}|\nabla_{\Gamma} \psi_{\Gamma,t}|^{2}+\frac{\delta(0)}{2}\frac{\d}{\d t}\int_{\Gamma}\mu_{5}^{2}|\Delta_{\Gamma} \psi_{\Gamma}|^{2} =\delta(0)\int_{\Gamma}\mu_{5}\mu_{5,t}|\Delta_{\Gamma} \psi_{\Gamma}|^{2}\nonumber\\
		&& +\sigma(0)\int_{\Gamma}\mu_{5}^{2}\Delta_{\Gamma}\psi_{\Gamma,t}\partial_{\nu}\psi +b^{\prime}(0)\int_{\Gamma}\mu_{5}^{2}\psi_{\Gamma}\Delta_{\Gamma} \psi_{\Gamma,t}  - \int_{\Gamma}\mu_{5}^{2}\Delta_{\Gamma} \psi_{\Gamma,t}f_{\Gamma}.
	\end{eqnarray*}
	Hence, 
	\begin{eqnarray}
		&& \int_{\Gamma}\mu_{5}^{2}|\nabla_{\Gamma} \psi_{\Gamma,t}|^{2}+\frac{\delta(0)}{2}\frac{\d}{\d t}\int_{\Gamma}\mu_{5}^{2}|\Delta_{\Gamma} \psi_{\Gamma}|^{2} \leq C \left(\int_{\Gamma}\mu_{3}^{2}|\Delta_{\Gamma} \psi_{\Gamma}|^{2} \right. \nonumber\\
		&& \left. + \int_{\Gamma}\mu_{5}^{2}|\Delta_{\Gamma}\psi_{\Gamma,t}|^{2}  + \int_{\Omega}\mu_{0}^{2}|\psi|^{2}+  \int_{\Omega}\mu_{3}^{2}|\Delta\psi|^{2}   + \int_{\Gamma}\mu_{0}^{2}|\psi_{\Gamma}|^{2} + \int_{\Gamma}\mu^{2}|f_{\Gamma}|^{2} \right).  \label{cc51}
	\end{eqnarray}
	Summing \eqref{cc50} and \eqref{cc51} and integrating over $(0,t)$, one has
	\begin{eqnarray*}
		&& \|\mu_{5}\nabla\Psi_{t}\|^{2}_{L^{2}(0,T;\mathbb{L}^{2})}+\sup_{0\leq t \leq T}\mu_{5}^{2}(t)\|\Delta\Psi(\cdot,t)\|^{2}_{\mathbb{L}^{2}}\leq 
		C\left(\|\mu_{3}\Psi_{t}\|^{2}_{L^{2}(0,T;\mathbb{L}^{2})} \right.\nonumber \\
		&& \left. + \|\mu_{3}\Delta\Psi\|^{2}_{L^{2}(0,T;\mathbb{L}^{2})}+\|\mu_{5}\Delta\Psi_{t}\|^{2}_{L^{2}(0,T;\mathbb{L}^{2})}+ \|\mu_{0}\Psi\|^{2}_{L^{2}(0,T;\mathbb{L}^{2})}+\|\mu_{1}v\|^{2}_{L^{2}(\omega_{T})}\right. \nonumber\\
		&& \left. + \|\mu F\|^{2}_{L^{2}(0,T;\mathbb{L}^{2})}\right).
	\end{eqnarray*}
	Using \eqref{c26}, \eqref{c21} and \eqref{c49}, we obtain 
	\begin{eqnarray}
		&& \|\mu_{5}\nabla\Psi_{t}\|^{2}_{L^{2}(0,T;\mathbb{L}^{2})}+\sup_{0\leq t \leq T}\mu_{5}^{2}(t)\|\Delta\Psi(\cdot,t)\|^{2}_{\mathbb{L}^{2}}\nonumber\\
		&& \leq C\left(\|\mu F\|^{2}_{L^{2}(0,T;\mathbb{L}^{2})}+\|\mu G\|^{2}_{L^{2}(0,T;\mathbb{L}^{2})}+\|\mu_{4} F_{t}\|^{2}_{L^{2}(0,T;\mathbb{L}^{2})}\right). \label{cc52}
	\end{eqnarray}
	Finally, from \eqref{c46}, \eqref{c49} and \eqref{cc52}, we deduce \eqref{c28}.
\end{proof}

	\section{Null Controllability of the quasilinear system \eqref{cascade quasi-linear system}}\label{Section6}
We will establish the local null controllability of system \eqref{cascade quasi-linear system} using Lyusternik-Graves’ Theorem.

\subsection{Study of the mapping $\varLambda$ given in \eqref{varlamda}} We introduce suitable spaces so that the mapping $\varLambda$ verifies the conditions of Theorem \ref{Lyusternik}.
\begin{eqnarray*}
	&&\mathbb{X}:=\displaystyle\left\{(\Psi, H, v)\;: \quad\mu_{0}\Psi,\, \mu_{0}H,\, \mu_{3}\Delta H,\; \mu_{4}\Psi_{t},\, \mu_{5}\Delta\Psi_{t}\in L^{2}(0,T;\mathbb{L}^{2}),\; \mu_{1}v,\, \mu_{3}v_{t}\in L^{2}(\omega_{T}),  \right.\\
	&& \left.   v\mathds{1}_{\omega_{T}}\in  L^{2}(0,T;H^{2}(\Omega)),\; \mu(\mathbf{L}_{1}\Psi-v\mathds{1}_{\omega},\mathbf{L}_{2}\Psi), \mu_{4}(\mathbf{L}_{1}\Psi-v\mathds{1}_{\omega},\mathbf{L}_{2}\Psi)_{t}\in L^{2}(0,T;\mathbb{L}^{2}), \right.\\ 
	&& \left. \mu(\mathbf{L}_{1}^{*}H-\theta\psi\mathds{1}_{\mathcal{O}}, \mathbf{L}_{2}^{*}H-\theta_{\Gamma}\psi_{\Gamma}\mathds{1}_{\Sigma})\in L^{2}(0,T;\mathbb{L}^{2}),\;  \right.\\
	&&\left. \sup_{0\leq t \leq  T}\mu_{5}^{2}(t)\|\Psi_{t}(\cdot,t)\|^{2}_{\mathbb{H}^{1}}<\infty,\;\sup_{0\leq t \leq  T}\mu_{5}^{2}(t)\|\Psi(\cdot,t)\|^{2}_{\mathbb{H}^{2}}<\infty \right\},\\
	&&\mathbb{Y}_1:=\{F\;:\; \mu F,\, \mu_{4} F_{t}\in L^{2}(0,T;\mathbb{L}^{2}) \},\quad \mathbb{Y}_2:=\{G\;:\; \mu G\in L^{2}(0,T;\mathbb{L}^{2}) \},\nonumber\\
	&&\mathbb{Y}:=\mathbb{Y}_1\times \mathbb{Y}_2,
\end{eqnarray*}
where $\mathbf{L}_{1}$ and $\mathbf{L}_{2}$ are defined in \eqref{n1}. These spaces are naturally equipped with the following norms
\begin{eqnarray*}
	&&\|(\Psi, H, v)\|_{\mathbb{X}}:=\left(\|\mu_{0}\Psi\|^{2}_{L^{2}(0,T;\mathbb{L}^{2})}+\|\mu_{0}H\|^{2}_{L^{2}(0,T;\mathbb{L}^{2})}  +\|\mu_{3}\Delta H\|^{2}_{L^{2}(0,T;\mathbb{L}^{2})} +\|\mu_{4}\Psi_{t}\|^{2}_{L^{2}(0,T;\mathbb{L}^{2})}  \right.\\
	&& \left.  +\|\mu_{5}\Delta\Psi_{t}\|^{2}_{L^{2}(0,T;\mathbb{L}^{2})} + \|\mu_1v\|^{2}_{L^{2}(\omega_{T})}+ \|\mu_3v_t\|^{2}_{L^{2}(\omega_{T})}+ \|v\|^{2}_{L^{2}(0,T;H^{2}(\omega))}, \right.\\
	&& \left. +\|\mu (\mathbf{L}_{1}\Psi-v\mathds{1}_{\omega},\mathbf{L}_{2}\Psi)\|^{2}_{L^{2}(0,T;\mathbb{L}^{2})} +\|\mu_{4} (\mathbf{L}_{1}\Psi-v\mathds{1}_{\omega},\mathbf{L}_{2}\Psi)_{t}\|^{2}_{L^{2}(0,T;\mathbb{L}^{2})} \right. \\
	&& \left.     +\|\mu (\mathbf{L}_{1}^{\star}H-\theta\psi\mathds{1}_{\mathcal{O}},\mathbf{L}_{2}^{\star}H-\theta_{\Gamma}\psi_{\Gamma}\mathds{1}_{\Sigma})\|^{2}_{L^{2}(0,T;\mathbb{L}^{2})}  + \sup_{0\leq t \leq  T}\mu_{5}^{2}(t)\|\Psi_{t}(\cdot,t)\|^{2}_{\mathbb{H}^{1}}  +\sup_{0\leq t \leq  T}\mu_{5}^{2}(t)\|\Psi(\cdot,t)\|^{2}_{\mathbb{H}^{2}} \right)^{1/2},\\
	 &&\|F\|_{\mathbb{Y}}:=\left(\|\mu F\|^{2}_{L^{2}(0,T;\mathbb{L}^{2})}+\|\mu G\|^{2}_{L^{2}(0,T;\mathbb{L}^{2})}+\|\mu_{4}F_{t}\|^{2}_{L^{2}(0,T;\mathbb{L}^{2})}\right)^{1/2}.
\end{eqnarray*}
It is straightforward to show that these are Hilbert spaces.
\begin{remark} \label{Remark 3}
	One can easily observe that there exists a constant $C>0$ such that
	\begin{eqnarray*}
		\int_{0}^{T}\mu_{5}^{2}(t)\|\Psi_{t}(\cdot,t)\|^{2}_{\mathbb{H}^{2}}\leq C\|(\Psi, H, v)\|^{2}_{\mathbb{X}} \quad \forall (\Psi, H, v)\in\mathbb{X}.
	\end{eqnarray*}
\end{remark}
Consider the mapping
$\varLambda:\mathbb{X}\longrightarrow\mathbb{Y}$ given in \eqref{varlamda}.
	To simplify, we consider the linear part of $\varLambda$:
\begin{eqnarray*}
	&&L(\Psi, H, v):=(\mathbf{L}_1\Psi-
	v\mathds{1}_{\omega}, \mathbf{L}_2\Psi, \mathbf{L}^{\star}_1 H -\theta\psi\mathds{1}_{\mathcal{O}},  \mathbf{L}^{\star}_{2}H-\theta_{\Gamma}\psi_{\Gamma}\mathds{1}_{\Sigma} ) \label{d-2}
\end{eqnarray*}
and the nonlinear part of $\varLambda$:
\begin{eqnarray*}
	A(\Psi, H, v)&:=& L(\Psi, H, v)-\varLambda(\Psi, H, v)\\
	&:=& (A_1(\Psi, H, v), A_3(\Psi, H, v), A_2(\Psi, H, v), A_4(\Psi, H, v)),
\end{eqnarray*}
where
\begin{eqnarray*}
	&&A_1(\Psi, H, v):= \nabla\cdot((\sigma(\psi)-\sigma(0))\nabla\psi)-(a(\psi)-a^{\prime}(0)\psi),\\
	&&A_2(\Psi, H, v) := (\sigma(\psi)-\sigma(0))\Delta h -(a^{\prime}(\psi)-a^{\prime}(0))h,\\
	&&A_3(\Psi, H, v):=  \nabla_{\Gamma}\cdot((\delta(\psi_{\Gamma})-\delta(0))\nabla_{\Gamma}\psi_{\Gamma})-(\sigma(\psi_{\Gamma})-\sigma(0))\partial_{\nu}\psi -(b(\psi_{\Gamma})-b^{\prime}(0)\psi_{\Gamma}),\\
	&&A_4(\Psi, H, v) :=  (\delta(\psi_{\Gamma})-\delta(0))\Delta_{\Gamma} h_{\Gamma}-(\sigma(\psi_{\Gamma})-\sigma(0))
	\partial_{\nu} h-(b^{\prime}(\psi_{\Gamma})-b^{\prime}(0))h_{\Gamma}.
\end{eqnarray*}
The following lemma confirms that the mapping 
$\varLambda$ is well-defined.
\begin{lemma} \label{Lemma3}
	Let $r>0$. 
	There is a positive constant $C:=C(r)>0$ such that the mapping $\varLambda : \mathbb{X}\rightarrow\mathbb{Y}$ defined above verifies
	\begin{eqnarray}
		\|\varLambda(\Psi, H, v)\|^{2}_{\mathbb{Y}}\leq C\left(\|(\Psi, H, v)\|^{2}_{\mathbb{X}}+ \|(\Psi, H, v)\|^{4}_{\mathbb{X}}+ \|(\Psi, H, v)\|^{6}_{\mathbb{X}}\right),   \label{d0}
	\end{eqnarray}
	for all $(\Psi, H, v)\in \overline{B}_{\mathbb{X}}(0,r)$. In particular, $\varLambda : \mathbb{X}\rightarrow\mathbb{Y}$ is well defined.
\end{lemma}
\begin{proof} 
	Let $r>0$ and $(\Psi, H, v)\in \overline{B}_{\mathbb{X}}(0,r)$. 
	From the definition of the norm of $\mathbb{X}$ and the continuity of Sobolev embedding $\mathbb{H}^{2}\hookrightarrow \mathbb{L}^{\infty}$, we can easily obtain that 
	\begin{eqnarray}
		\|\Psi\|_{L^{\infty}(0,T;\mathbb{L}^{\infty})}\leq C_{0} \|(\Psi, H, v)\|_{\mathbb{X}}\leq C_{0}r \quad \mbox{for all} \, (\Psi, H, v)\in \overline{B}_{\mathbb{X}}(0,r). \label{constant}
	\end{eqnarray}
	Throughout this proof, we will use the fact that the $i^{\mbox{th}}$ derivative of $\sigma$ and $\delta$, the $j^{\mbox{th}}$ derivative of $a$ and $b$ are Lipschitz-continuous on the  interval $J_{0}:=[-C_{0}r,C_{0}r]$ for $i=0,1,2$ and $j=0,1$. Firstly, according to the definition of the norm of $\mathbb{X}$, the linear part is bounded:
	\begin{eqnarray}
		\|L(\Psi, H, v)\|^{2}_{\mathbb{Y}}
		&=& \|\mu(\mathbf{L}_{1}\Psi-v\mathds{1}_{\omega},\mathbf{L}_{2}\Psi)\|^{2}_{L^{2}(0,T;\mathbb{L}^{2})} + 	\|\mu_{4}(\mathbf{L}_{1}\Psi-v\mathds{1}_{\omega},\mathbf{L}_{2}\Psi)_{t}\|^{2}_{L^{2}(0,T;\mathbb{L}^{2})} \nonumber\\
		&& + \|\mu(\mathbf{L}_{1}^{\star}H-\theta\psi\mathds{1}_{\mathcal{O}},\mathbf{L}_{2}^{\star}H-\theta_{\Gamma}\psi_{\Gamma}\mathds{1}_{\Sigma})\|^{2}_{L^{2}(0,T;\mathbb{L}^{2})} \nonumber\\
		&\leq & \|(\Psi, H, v)\|^{2}. \label{linear part bounded}
	\end{eqnarray}
	\\
	Now, we analyze the nonlinear part $A$.  One has
	\begin{eqnarray}
		&&\|A(\Psi, H, v)\|^{2}_{\mathbb{Y}}
		= \|\mu(A_{1}(\Psi, H, v), A_{3}(\Psi, H, v))\|^{2}_{L^{2}(0,T;\mathbb{L}^{2})} \nonumber\\
		&&+ 	\|\mu_{4}(A_{1}(\Psi, H, v),A_{3}(\Psi, H, v))_{t}\|^{2}_{L^{2}(0,T;\mathbb{L}^{2})}  + \|\mu(A_{2}(\Psi, H, v),A_{4}(\Psi, H, v))\|^{2}_{L^{2}(0,T;\mathbb{L}^{2})}. \label{d-1}
	\end{eqnarray}
	Then, it suffices to show that each term of \eqref{d-1} is bounded as in estimate \eqref{d0}.\\
	Firstly, we claim that
	\begin{eqnarray}
		\|\mu(A_{1}(\Psi, H, v), A_{3}(\Psi, H, v))\|^{2}_{L^{2}(0,T;\mathbb{L}^{2})}\leq C\left(\|(\Psi, H, v)\|^{2}_{\mathbb{X}}+ \|(\Psi, H, v)\|^{4}_{\mathbb{X}}\right). \label{d1}
	\end{eqnarray}
	Using $a$ and $b$ are Lipschitz-continuous on $J_{0}$ and $a(0)=b(0)=0$, we have
	\begin{eqnarray}
		&&\|\mu(\varLambda_{1}(\Psi, H, v),\varLambda_{3}(\Psi,v))\|^{2}_{L^{2}(0,T;\mathbb{L}^{2})}=\|\mu \varLambda_{1}(\Psi, H, v)\|^{2}_{L^{2}(\Omega_{T})}  +  \|\mu \varLambda_{3}(\Psi, H, v)\|^{2}_{L^{2}(\Gamma_{T})} \nonumber\\ 
		&&  \leq  
		C\left(  \int_{\Omega_{T}}\mu^{2}|\nabla\cdot\left((\sigma(\psi)-\sigma(0))\nabla \psi\right)|^{2}  + \int_{\Omega_{T}}\mu^{2}|a(\psi)-a^{\prime}(0)\psi|^{2} \right. \nonumber\\
		&& \left.  + \int_{\Gamma_{T}}\mu^{2}|\nabla_{\Gamma}\cdot\left((\delta(\psi_{\Gamma})-\delta(0))\nabla_{\Gamma} \psi_{\Gamma}\right)|^{2} +\int_{\Gamma_{T}}\mu^{2}|\psi_{\Gamma}|^{2}|\partial_{\nu}\psi|^{2}  \right. \nonumber \\
		&& \left.  + \int_{\Gamma_{T}}\mu^{2}|b(\psi_{\Gamma})-b^{\prime}(0)\psi_{\Gamma}|^{2}  \right) := C\left(\sum_{j=1}^{5}I_{j}\right).  \label{d2}
	\end{eqnarray}
	Let us analyze $I_1$. Since $\sigma$ is Lipschitz-continuous and bounded on $J_{0}$, one has 
	\begin{eqnarray*}
		&&I_{1}= \int_{\Omega_{T}}\mu^{2}|\sigma^{\prime}(\psi)|^{2}|\nabla\psi|^{4}+ \int_{\Omega_{T}}\mu^{2}|(\sigma(\psi)-\sigma(0))\Delta \psi|^{2}\\
		&&\leq  C\left( \int_{\Omega_{T}}\mu^{2}|\nabla\psi|^{4}+ \int_{\Omega_{T}}\mu^{2}|\psi|^{2}|\Delta \psi|^{2} \right)\\
		&&\leq  C\left(\int_{0}^{T}\left[\mu^{2}(t)\int_{\Omega}|\nabla\psi(\cdot,t)|^{4}\right]  + \int_{0}^{T}\left[\mu^{2}(t)\|\psi(\cdot,t)\|^{2}_{L^{\infty}(\Omega)}\int_{\Omega}|\Delta\psi(\cdot,t)|^{2}\right] \right).
	\end{eqnarray*}
	The fact that $H^{2}(\Omega) \hookrightarrow W^{1,4}(\Omega)$ and  $H^{2}(\Omega) \hookrightarrow L^{\infty}(\Omega)$ with continuous embeddings, implies that
	\begin{eqnarray}
		I_{1}&\leq & C\int_{0}^{T}\mu^{2}(t)\|\psi(\cdot,t)\|^{4}_{H^{2}(\Omega)} \nonumber\\
		&\leq & C\sup_{0\leq t \leq T}\mu_{5}^{4}(t)\|\psi(\cdot,t)\|^{4}_{H^{2}(\Omega)} \nonumber\\
		&\leq & C\|(\Psi, H, v)\|^{4}_{\mathbb{X}}. \label{d4}
	\end{eqnarray}
	Similarly to $I_1$, we obtain
	\begin{eqnarray}
		I_{3} &\leq & C\|(\Psi, H, v)\|^{4}_{\mathbb{X}}. \label{d5}
	\end{eqnarray}  
	Using $a(0)=b(0)=0$ and Taylor's Inequality, we obtain  
	\begin{eqnarray}
		I_{2}+I_{5}&\leq & \frac{M^{2}}{4}\int_{\Omega_{T}}\mu^{2}|\psi|^{4}+ \frac{M^{2}_{\Gamma}}{4}\int_{\Gamma_{T}}\mu^{2}|\psi_{\Gamma}|^{4} \nonumber\\
		&\leq & C\left(\int_{\Omega_{T}}\mu_{5}^{4}|\psi|^{4}+ \int_{\Gamma_{T}}\mu_{5}^{4}|\psi_{\Gamma}|^{4} \right)  \nonumber\\
		&\leq & C \sup_{0\leq t \leq  T}\mu_{5}^{4}(t)\|\Psi(\cdot,t)\|^{4}_{L^{2}(0,T;\mathbb{L}^{2})}   \nonumber\\
		&\leq & C\|(\Psi, H, v)\|^{4}_{\mathbb{X}}, \label{d6}
	\end{eqnarray}
	where $M=\displaystyle\sup_{r\in J_{0}}|a^{\prime\prime}(r)|$ and $M_{\Gamma}=\displaystyle\sup_{r\in J_{0}}|b^{\prime\prime}(r)|$. \\
	For the last term $I_{4}$, thanks to the continuity of the Sobolev embedding $H^{2}(\Gamma)\hookrightarrow L^{\infty}(\Gamma)$ and the continuity of the normal derivative from $H^{2}(\Omega)$ to $L^{2}(\Gamma)$, we have
	\begin{eqnarray}
		I_{4}&\leq& C \int_{0}^{T} \mu^{2}(t) \|\psi_{\Gamma}(\cdot,t)\|^{2}_{L^{\infty}(\Gamma)}\|\partial_{\nu}\psi(\cdot,t)\|^{2}_{L^{2}(\Gamma)}  \nonumber\\
		&\leq& C \int_{0}^{T} \mu^{2}(t) \|\psi_{\Gamma}(\cdot,t)\|^{2}_{H^{2}(\Gamma)}\|\psi(\cdot,t)\|^{2}_{H^{2}(\Omega)}  \nonumber\\
		&\leq& C \int_{0}^{T} \mu^{2}(t) \|\Psi(\cdot,t)\|^{4}_{\mathbb{H}^{2}}  \nonumber\\
		&\leq&  C \sup_{0\leq t \leq T}\mu_{5}^{4}(t)\|\Psi(\cdot,t)\|^{4}_{\mathbb{H}^{2}} \nonumber\\
		&\leq & C \|(\Psi, H, v)\|_{\mathbb{X}}^{4}. \label{d7}
	\end{eqnarray}
	From estimates \eqref{d2}-\eqref{d7}, we obtain \eqref{d1}.\\
	Now, we prove that
	\begin{eqnarray}
		\|\mu_{4}(A_{1}(\Psi, H, v),A_{3}(\Psi, H, v))_{t}\|^{2}_{L^{2}(0,T;\mathbb{L}^{2})}\leq C\left(\|(\Psi, H, v)\|^{2}_{\mathbb{X}}+ \|(\Psi, H, v)\|^{4}_{\mathbb{X}}+ \|(\Psi, H, v)\|^{6}_{\mathbb{X}}\right).\nonumber\\ \label{d8}
	\end{eqnarray}
	We have
	\begin{eqnarray}
		&& \|\mu_{4}(\varLambda_{1}(\Psi, H, v),\varLambda_{2}(\Psi, H, v))_{t}\|^{2}_{L^{2}(0,T;\mathbb{L}^{2})}=	\|\mu_{4}(\varLambda_{1}(\Psi, H, v))_{t}\|^{2}_{L^{2}(\Omega_{T})}    \label{d9} \\
		&& + \|\mu_{4}(\varLambda_{2}(\Psi, H, v))_{t}\|^{2}_{L^{2}(\Gamma_{T})} \leq C\left(  \int_{\Omega_{T}}\mu_{4}^{2}\left|\left(\nabla\cdot\left((\sigma(\psi)-\sigma(0))\nabla \psi\right)\right)_{t}\right|^{2}\right. \nonumber \\
		&&  \left.+\int_{\Omega_{T}}\mu_{4}^{2}|(a(\psi)-a^{\prime}(0)\psi)_{t}|^{2} + \int_{\Gamma_{T}}\mu_{4}^{2}\left|\left(\nabla_{\Gamma}\cdot\left((\delta(\psi_{\Gamma})-\delta(0))\nabla_{\Gamma} \psi_{\Gamma}\right)\right)_{t}\right|^{2}  \right. \nonumber\\
		&& \left.    + \int_{\Gamma_{T}}\mu_{4}^{2}\left|\left( \left(\sigma(\psi_{\Gamma})-\sigma(0)\right)\partial_{\nu}\psi\right)_{t} \right|^{2}  + \int_{\Gamma_{T}}\mu_{4}^{2}|(b(\psi_{\Gamma})-b^{\prime}(0)\psi_{\Gamma})_{t}|^{2} \right):= C\left(\sum_{j=1}^{5}J_{j}\right). \nonumber
	\end{eqnarray}
	Let us show the estimate for $J_{1}$. Using  $\sigma^{\prime\prime}$ and $\sigma^{\prime}$ are bounded on $J_{0}$, $\sigma$ is Lipschitz-continuous on $J_0$, we have
	\begin{eqnarray*}
		J_{1} &\leq & C\left(\int_{\Omega_{T}}\mu_{4}^{2}|\sigma^{\prime\prime}(\psi)\psi_{t}|^{2}|\nabla\psi|^{4}+ \int_{\Omega_{T}}\mu_{4}^{2}|\sigma^{\prime}(\psi)\nabla\psi_{t}\cdot\nabla\psi|^{2} \right.\nonumber\\
		&&\left. + \int_{\Omega_{T}}\mu_{4}^{2}|\sigma^{\prime}(\psi)\psi_{t}\Delta\psi|^{2} + \int_{\Omega_{T}}\mu_{4}^{2}|(\sigma(\psi)-\sigma(0))\Delta\psi_{t}|^{2}\right) \nonumber \\
		&\leq & C\left(\int_{\Omega_{T}}\mu_{4}^{2}|\psi_{t}|^{2}|\nabla\psi|^{4}+ \int_{\Omega_{T}}\mu_{4}^{2}|\nabla\psi_{t}\cdot\nabla\psi|^{2} + \int_{\Omega_{T}}\mu_{4}^{2}|\psi_{t}|^{2}|\Delta\psi|^{2}\right.\nonumber\\
		&&\left.  + \int_{\Omega_{T}}\mu_{4}^{2}|\psi|^{2}|\Delta\psi_{t}|^{2}\right):= C\left(\sum_{j=1}^{4}J_{1j}\right).
	\end{eqnarray*}
	Using Sobolev embeddings as above and Remark \ref{Remark 3}, we get
	\begin{eqnarray*}
		J_{11} &\leq & \int_{0}^{T}\left[\mu_{4}^{2}(t)\|\psi_{t}(\cdot,t)\|_{L^{\infty}(\Omega)}^{2}\left(\int_{\Omega}|\nabla\psi(\cdot,t)|^{4}\right)\right]\\
		&\leq & C\int_{0}^{T}\left[\mu_{4}^{2}(t)\|\psi_{t}(\cdot,t)\|_{H^{2}(\Omega)}^{2}\|\psi(\cdot,t)\|_{H^{2}(\Omega)}^{4}\right]\\
		&\leq & C\left(\sup_{0\leq t \leq  T}\mu_{5}^{4}(t)\|\psi(\cdot,t)\|^{4}_{H^{2}(\Omega)} \right)\int_{0}^{T}\mu_{5}^{2}(t)\|\psi_{t}(\cdot,t)\|_{H^{2}(\Omega)}^{2}\\
		&\leq & C\|(\Psi, H, v)\|^{6}_{\mathbb{X}}.
	\end{eqnarray*}
	Applying Cauchy-Schwarz inequality, we obtain 
	\begin{eqnarray*}
		J_{12}&\leq & \int_{0}^{T}\left[\mu_{4}^{2}(t)\left(\int_{\Omega}|\nabla\psi_{t}(\cdot,t)|^{2}|\nabla\psi(\cdot,t)|^{2}\right) \right]\\
		&\leq & \int_{0}^{T}\left[\mu_{4}^{2}(t)\left(\int_{\Omega}|\nabla\psi_{t}(\cdot,t)|^{4}\right)^{1/2}\left(\int_{\Omega}|\nabla\psi(\cdot,t)|^{4}\right)^{1/2} \right]\\
		&\leq & \int_{0}^{T}\mu_{4}^{2}(t)\|\psi_{t}(\cdot,t)\|_{H^{2}(\Omega)}^{2}\|\psi(\cdot,t)\|_{H^{2}(\Omega)}^{2} \\
		&\leq & C\left(\sup_{0\leq t \leq T}\mu_{5}^{2}(t)\|\psi(\cdot,t)\|^{2}_{H^{2}(\Omega)}\right)\int_{0}^{T}\mu_{5}^{2}(t)\|\psi_{t}(\cdot,t)\|_{H^{2}(\Omega)}^{2} \\
		&\leq & C\|(\Psi, H, v)\|^{4}_{\mathbb{X}}.
	\end{eqnarray*}
	Using the fact that $H^{2}(\Omega)\hookrightarrow L^{\infty}(\Omega)$ is continuous and Remark \ref{Remark 3}, one has
	\begin{eqnarray*}
		J_{13} &\leq & \int_{0}^{T}\left[\mu_{4}^{2}(t)\|\psi_{t}(\cdot,t)\|_{L^{\infty}(\Omega)}^{2}\left(\int_{\Omega}|\Delta\psi(\cdot,t)|^{2} \right)\right] \\
		&\leq & C\int_{0}^{T}\mu_{4}^{2}(t)\|\psi_{t}(\cdot,t)\|^{2}_{H^{2}(\Omega)}\|\Delta\psi(\cdot,t)\|^{2}_{L^{2}(\Omega)} \\
		&\leq & C\left(\sup_{0\leq t \leq T}\mu_{5}^{2}(t)\|\Delta\psi(\cdot,t)\|^{2}_{L^{2}(\Omega)}\right)\int_{0}^{T}\mu_{5}^{2}(t)\|\psi_{t}(\cdot,t)\|^{2}_{H^{2}(\Omega)} \\
		&\leq & C\|(\Psi, H, v)\|^{4}_{\mathbb{X}}.
	\end{eqnarray*}
	We also have
	\begin{eqnarray*}
		J_{14} &\leq & \int_{0}^{T}\left[\mu_{4}^{2}(t)\|\psi(\cdot,t)\|_{L^{\infty}(\Omega)}^{2}\left(\int_{\Omega}|\Delta\psi_{t}(\cdot,t)|^{2}\right)\right] \\
		&\leq & C\int_{0}^{T}\mu_{4}^{2}(t)\|\psi(\cdot,t)\|^{2}_{H^{2}(\Omega)}\|\Delta\psi_{t}(\cdot,t)\|^{2}_{L^{2}(\Omega)} \\
		&\leq & C\left(\sup_{0\leq t \leq T}\mu_{5}^{2}(t)\|\psi(\cdot,t)\|^{2}_{H^{2}(\Omega)}\right)\int_{0}^{T}\mu_{5}^{2}(t)\|\Delta\psi_{t}(\cdot,t)\|^{2}_{L^{2}(\Omega)} \\
		&\leq & C\|(\Psi, H, v)\|^{4}_{\mathbb{X}}.
	\end{eqnarray*}
	Consequently,
	\begin{eqnarray}
		J_{1}\leq C\left(\|(\Psi, H, v)\|^{6}_{\mathbb{X}}+\|(\Psi, H, v)\|^{4}_{\mathbb{X}}\right). \label{d11}
	\end{eqnarray}
	Analogously, we obtain 
	\begin{eqnarray}
		J_{3}\leq C\left(\|(\Psi, H, v)\|^{6}_{\mathbb{X}}+\|(\Psi, H, v)\|^{4}_{\mathbb{X}}\right). \label{d12}
	\end{eqnarray}
		Using  $a^{\prime}$ and $b^{\prime}$ are bounded on $J_{0}$, we find
	\begin{eqnarray}
		&& J_{2}+J_{5}\leq C\|\mu_{4}\Psi_{t}\|^{2}_{L^{2}(0,T;\mathbb{L}^{2})}\leq C\|(\Psi, H, v)\|_{\mathbb{X}}^{2}. \label{d13}
	\end{eqnarray}
	It remains to analyze the last term $J_{4}$.
	Using $\sigma^{\prime}$ is bounded on $J_{0}$, $\sigma$ is Lipscitz-continuous on $J_{0}$, $H^{2}(\Gamma)\hookrightarrow L^{\infty}(\Gamma)$ is continuous, the normal derivative is continuous from $H^{2}(\Omega)$ to $L^{2}(\Gamma)$  and Remark \ref{Remark 3}, we find 
	\begin{eqnarray}
		J_{4} &\leq & C\left(\int_{\Gamma_{T}}\mu_{4}^{2}|\psi_{\Gamma,t}|^{2}|\partial_{\nu}\psi|^{2}+ \int_{\Gamma_{T}}\mu_{4}^{2}|\psi_{\Gamma}|^{2}|\partial_{\nu}\psi_{t}|^{2} \right) \nonumber\\
		&\leq& C\left( \int_{0}^{T}\left[\mu_{4}^{2}(t)\|\psi_{\Gamma,t}(\cdot,t)\|^{2}_{L^{\infty}(\Gamma)}\int_{\Gamma}|\partial_{\nu}\psi(\cdot,t)|^{2}\right]  \right. \nonumber\\
		&&\left. + \int_{0}^{T}\left[\mu_{4}^{2}(t)\|\psi_{\Gamma}(\cdot,t)\|^{2}_{L^{\infty}(\Gamma)}\int_{\Gamma}|\partial_{\nu}\psi_{t}(\cdot,t)|^{2}\right] \right) \nonumber\\
		&\leq & C\left( \int_{0}^{T}\mu_{4}^{2}(t)\|\psi_{\Gamma,t}(\cdot,t)\|^{2}_{H^{2}(\Gamma)}\|\psi(\cdot,t)\|^{2}_{H^{2}(\Omega)}  \right. \nonumber\\
		&&\left. + \int_{0}^{T}\mu_{4}^{2}(t)\|\psi_{\Gamma}(\cdot,t)\|^{2}_{H^{2}(\Gamma)}\|\psi_{t}(\cdot,t)\|^{2}_{H^{2}(\Omega)} \right) \nonumber\\
		&\leq& C\left(\sup_{0\leq t \leq  T}\mu_{5}^{2}(t)\|\psi(\cdot,t)\|^{2}_{H^{2}(\Omega)}\int_{0}^{T}\mu_{5}^{2}(t)\|\psi_{\Gamma,t}(\cdot,t)\|^{2}_{H^{2}(\Gamma)} \right. \nonumber\\
		&& \left. +  \sup_{0\leq t \leq  T}\mu_{5}^{2}(t)\|\psi_{\Gamma}(\cdot,t)\|^{2}_{H^{2}(\Gamma)}\int_{0}^{T}\mu_{5}^{2}(t)\|\psi_{t}(\cdot,t)\|^{2}_{H^{2}(\Omega)} \right) \nonumber \\
		&\leq & C\|(\Psi, H, v)\|^{4}_{\mathbb{X}}. \label{d14}
	\end{eqnarray}
	From estimates \eqref{d9}-\eqref{d14}, we obtain \eqref{d8}.\\
	Now, we claim that
	\begin{eqnarray}
		&&\|\mu(\varLambda_{2}(\Psi, H, v),\varLambda_{4}(\Psi, H, v))\|^{2}_{L^{2}(0,T;\mathbb{L}^{2})} \leq C\left(\|(\Psi, H, v)\|^{2}_{\mathbb{X}}+ \|(\Psi, H, v)\|^{4}_{\mathbb{X}}\right). \label{kd1}
	\end{eqnarray}
	We have
	\begin{eqnarray}
		&&\|\mu(\varLambda_{2}(\Psi, H, v),\varLambda_{4}(\Psi, H, v))\|^{2}_{L^{2}(0,T;\mathbb{L}^{2})} =\|\mu \varLambda_{2}(\Psi, H, v)\|^{2}_{L^{2}(\Omega_{T})} +  \|\mu \varLambda_{4}(\Psi, H, v)\|^{2}_{L^{2}(\Gamma_{T})} \nonumber\\ 
		&& \leq  
		C\left( \int_{\Omega_{T}}\mu^{2}|(\sigma(\psi)-\sigma(0))\Delta h|^{2}  + \int_{\Omega_{T}}\mu^{2}|(a^{\prime}(\psi)-a^{\prime}(0))h|^{2}\right. \nonumber\\
		&& \left. +  \int_{\Gamma_{T}}\mu^{2}|(\delta(\psi_{\Gamma})-\delta(0))\Delta_{\Gamma} h_{\Gamma}|^{2} +\int_{\Gamma_{T}}\mu^{2}|(\sigma(\psi_{\Gamma})-\sigma(0))\partial_{\nu}h|^{2} \right. \nonumber\\
		&&\left.    + \int_{\Gamma_{T}}\mu^{2}|(b^{\prime}(\psi_{\Gamma})-b^{\prime}(0))h_{\Gamma}|^{2}   \right) := C\left(\sum_{j=1}^{7}K_{j}\right). \label{kd2}
	\end{eqnarray}
	Using $\sigma$ and $\delta$ are Lipschitz-continuous on $J_{0}$ and $H^{2}(\Omega) \hookrightarrow L^{\infty}(\Omega)$ with continuous embeddings, one has
	\begin{align}
		K_{1}+K_{3}&\leq  C\left(\int_{\Omega_{T}}\mu^{2}|\psi|^{2}|\Delta h|^{2}+ \int_{\Gamma_{T}}\mu^{2}|\psi_{\Gamma}|^{2}|\Delta_{\Gamma}h_{\Gamma}|^{2} \right) \nonumber\\
		&\leq  C\left(\int_{0}^{T}\left[\mu_{5}^{2}(t)\|\psi(\cdot,t)\|^{2}_{L^{\infty}(\Omega)}\int_{\Omega}\mu_{3}^{2}|\Delta h|^{2}\right] \right. \nonumber\\
		& \left. + \int_{0}^{T}\left[\mu_{5}^{2}(t)\|\psi_{\Gamma}(\cdot,t)\|^{2}_{L^{\infty}(\Gamma)}\int_{\Gamma}\mu_{3}^{2}|\Delta_{\Gamma}h_{\Gamma}|^{2} \right]\right)  \nonumber\\
		&\leq  C\left( \sup_{0\leq t \leq T}\mu_{5}^{2}(t)\|\psi(\cdot,t)\|^{2}_{H^{2}(\Omega)}\|\mu_3 \Delta h\|^{2}_{L^{2}(0,T;L^{2}(\Omega))} \right.   \nonumber \\
		&\left. +\sup_{0\leq t \leq T}\mu_{5}^{2}(t)\|\psi(\cdot,t)\|^{2}_{H^{2}(\Gamma)}\|\mu_3 \Delta_{\Gamma}h_\Gamma\|^{2}_{L^{2}(0,T;L^{2}(\Gamma))}\right)  \nonumber \\
		&\leq  C \sup_{0\leq t \leq  T}\mu_{5}^{2}(t)\|\Psi(\cdot,t)\|^{2}_{\mathbb{H}^{2}}\|\mu_3 \Delta H\|^{2}_{L^{2}(0,T;\mathbb{L}^{2})}   \nonumber\\
		&\leq  C\|(\Psi, H, v)\|^{4}_{\mathbb{X}}. \label{kd5}
	\end{align} 		
    Let us now analyze $K_2$ and $K_5$. 
	Using $a^\prime$ and $b^\prime$ are Lipschitz-continuous on $J_{0}$ and $H^{2}(\Omega) \hookrightarrow L^{\infty}(\Omega)$ with continuous embeddings, we obtain 
	\begin{eqnarray}
		K_{2}+K_{5}&\leq & C\left(\int_{\Omega_{T}}\mu^{2}|\psi|^{2}|h|^{2}+ \int_{\Gamma_{T}}\mu^{2}|\psi_{\Gamma}|^{2}|h_{\Gamma}|^{2} \right) \nonumber\\
		&\leq & C\left(\int_{0}^{T}\left[\mu_{5}^{2}(t)\|\psi(\cdot,t)\|^{2}_{L^{\infty}(\Omega)}\int_{\Omega}\mu_{0}^{2}|h|^{2}\right] \right. \nonumber\\
		&& \left. + \int_{0}^{T}\left[\mu_{5}^{2}(t)\|\psi_{\Gamma}(\cdot,t)\|^{2}_{L^{\infty}(\Gamma)}\int_{\Gamma}\mu_{0}^{2}|h_{\Gamma}|^{2}\right]\right)  \nonumber\\
		&\leq & C\left( \sup_{0\leq t \leq T}\mu_{5}^{2}(t)\|\psi(\cdot,t)\|^{2}_{H^{2}(\Omega)}\|\mu_0 h\|^{2}_{L^{2}(0,T;L^{2}(\Omega))} \right.   \nonumber \\
		&&\left. +\sup_{0\leq t \leq T}\mu_{5}^{2}(t)\|\psi(\cdot,t)\|^{2}_{H^{2}(\Gamma)}\|\mu_0 h_\Gamma\|^{2}_{L^{2}(0,T;L^{2}(\Gamma))}\right)  \nonumber \\
		&\leq & C \sup_{0\leq t \leq  T}\mu_{5}^{2}(t)\|\Psi(\cdot,t)\|^{2}_{\mathbb{H}^{2}}\|\mu_0 H\|^{2}_{L^{2}(0,T;\mathbb{L}^{2})}   \nonumber\\
		&\leq & C\|(\Psi, H, v)\|^{4}_{\mathbb{X}}. \label{kd4}
	\end{eqnarray}
	For the last term $K_{4}$, thanks to the continuity of the Sobolev embedding $H^{2}(\Gamma)\hookrightarrow L^{\infty}(\Gamma)$ and the continuity of the normal derivative from $H^{2}(\Omega)$ to $L^{2}(\Gamma)$, we have
	\begin{eqnarray}
		K_{4}&\leq& C \int_{0}^{T} \mu^{2}(t) \|\psi_{\Gamma}(\cdot,t)\|^{2}_{L^{\infty}(\Gamma)}\|\partial_{\nu}h(\cdot,t)\|^{2}_{L^{2}(\Gamma)} \nonumber\\
		&\leq& C \int_{0}^{T} \mu^{2}(t) \|\psi_{\Gamma}(\cdot,t)\|^{2}_{H^{2}(\Gamma)}\|h(\cdot,t)\|^{2}_{H^{2}(\Omega)}  \nonumber\\
		&\leq&  C \sup_{0\leq t \leq T}\mu_{5}^{2}(t)\|\Psi(\cdot,t)\|^{2}_{\mathbb{H}^{2}}\|\mu_3 H\|^{2}_{L^{2}(0,T;\mathbb{H}^{2})} \nonumber\\
		&\leq & C \|(\Psi,v)\|_{\mathbb{X}}^{4}. \label{kd6}
	\end{eqnarray}
	From estimates \eqref{kd2}-\eqref{kd6}, we obtain \eqref{kd1}.\\
Finally, by $\varLambda=L-A$ and estimates \eqref{linear part bounded}-\eqref{d1}, \eqref{d8} and \eqref{kd1}, we obtain \eqref{d0}.
\end{proof}

	\begin{lemma} \label{Lemma4}
	The mapping $\varLambda : \mathbb{X}\rightarrow\mathbb{Y}$ is continuously differentiable.
\end{lemma}
\begin{proof}
	Recall that $\varLambda=L-A$ and $L$ is linear continuous from $\mathbb{X}$ to $\mathbb{Y}$, then $L$ is continuously differentiable and $DL(\Psi, H, v)=L$, for any $(\Psi, H, v)\in\mathbb{X}$.
	Let us show that $A$ is continuously differentiable. We will start by proving that $A$ is Gateaux-differentiable at any $(\Psi, H, v)\in\mathbb{X}$ and determine the $G$-derivative $A^{\prime}(\Psi, H, v)$. Let $\varepsilon\in (-1,1)\setminus\{0\}$ and $(\Psi, H, v), (\Phi, K, u) \in\mathbb{X}$.\\ A simple computation yields
	\begin{eqnarray*}
		&&\frac{A_{1}((\Psi, H, v)+\varepsilon(\Phi, K, u))-A_{1}(\Psi, H, v)}{\varepsilon}
		= \frac{\sigma(\psi+\varepsilon\phi)-\sigma(\psi)}{\varepsilon}\Delta\psi \\
		&& + \left[\frac{\sigma^{\prime}(\psi+\varepsilon\phi)-\sigma^{\prime}(\psi)}{\varepsilon}\right]|\nabla\psi|^{2} +(\sigma(\psi+\varepsilon\phi)-\sigma(0))\Delta\phi\\
		&& +2\sigma^{\prime}(\psi+\varepsilon\phi)\nabla\psi\cdot\nabla\phi-\frac{a(\psi+\varepsilon\phi)-a(\psi)}{\varepsilon}+\varepsilon\sigma^{\prime}(\psi+\varepsilon\phi)|\nabla\phi|^{2}+a^{\prime}(0)\phi :=\sum_{j=1}^{7}A_{1,j}^\varepsilon,
	\end{eqnarray*}
	\begin{eqnarray*}
		&&\frac{A_{3}((\Psi, H, v)+\varepsilon(\Phi, K, u))-A_{3}(\Psi, H, v)}{\varepsilon}
		= \frac{\delta(\psi_{\Gamma}+\varepsilon\phi_{\Gamma})-\delta(\psi_{\Gamma})}{\varepsilon}\Delta_{\Gamma}\psi_{\Gamma}\nonumber\\
		&&+ \left[\frac{\delta^{\prime}(\psi_{\Gamma}+\varepsilon\phi_{\Gamma})-\delta^{\prime}(\psi_{\Gamma})}{\varepsilon}\right]|\nabla_{\Gamma}\psi_{\Gamma}|^{2} +(\delta(\psi_{\Gamma}+\varepsilon\phi_{\Gamma})-\delta(0))\Delta_{\Gamma}\phi_{\Gamma}\nonumber\\
		&&+2\delta^{\prime}(\psi_{\Gamma}+\varepsilon\phi_{\Gamma})\nabla_{\Gamma}\psi_{\Gamma}\cdot\nabla_{\Gamma}\phi_{\Gamma} -\frac{\sigma(\psi_{\Gamma}+\varepsilon\phi_{\Gamma})-\sigma(\psi_{\Gamma})}{\varepsilon}\partial_{\nu}\psi\\
		&&-(\sigma(\psi_{\Gamma}+\varepsilon\phi_{\Gamma})-\sigma(0))\partial_{\nu}\phi-\frac{b(\psi_{\Gamma}+\varepsilon\phi_{\Gamma})-b(\psi_{\Gamma})}{\varepsilon}  \nonumber\\
		&& +\varepsilon\delta^{\prime}(\psi_{\Gamma}+\varepsilon\phi_{\Gamma})|\nabla_{\Gamma}\phi_{\Gamma}|^{2} +b^{\prime}(0)\phi_{\Gamma}:=\sum_{j=1}^{9}A_{3,j}^{\varepsilon},
	\end{eqnarray*}
	\begin{eqnarray*}
		&&\frac{A_{2}((\Psi, H, v)+\varepsilon(\Phi, K, u))-A_{2}(\Psi, H, v)}{\varepsilon}
		= \frac{\sigma(\psi+\varepsilon\phi)-\sigma(\psi)}{\varepsilon}\Delta h \\
		&&  +(\sigma(\psi+\varepsilon\phi)-\sigma(0))\Delta k-\frac{a^{\prime}(\psi+\varepsilon\phi)-a^{\prime}(\psi)}{\varepsilon}h\nonumber\\
		&&-(a^{\prime}(\psi+\varepsilon\phi)-a^{\prime}(0))k:=\sum_{j=1}^{4}A_{2,j}^\varepsilon
	\end{eqnarray*}
	and 
	\begin{align*}
		&\frac{A_{4}((\Psi, H, v)+\varepsilon(\Phi, K, u))-A_{4}(\Psi, H, v)}{\varepsilon}
		= \frac{\delta(\psi_{\Gamma}+\varepsilon\phi_{\Gamma})-\delta(\psi_{\Gamma})}{\varepsilon}\Delta_{\Gamma} h_{\Gamma} \\
		& +(\delta(\psi_{\Gamma}+\varepsilon\phi_{\Gamma})-\delta(0))\Delta_{\Gamma} k_{\Gamma}-\frac{\sigma(\psi_{\Gamma}+\varepsilon\phi_{\Gamma})-\sigma(\psi_{\Gamma})}{\varepsilon}\partial_{\nu}h\\
		&-(\sigma(\psi_{\Gamma}+\varepsilon\phi_{\Gamma})-\sigma(0))\partial_{\nu}k-\frac{b^{\prime}(\psi_{\Gamma}+\varepsilon\phi_{\Gamma})-b^{\prime}(\psi_{\Gamma})}{\varepsilon}h_{\Gamma} \nonumber\\
		&-(b^{\prime}(\psi_{\Gamma}+\varepsilon\phi_{\Gamma})-b^{\prime}(0))k_{\Gamma}:=\sum_{j=1}^{6}A_{4,j}^{\varepsilon},
	\end{align*}
	We consider the linear mapping $\mathbf{D}A(\Psi, H, v): \mathbb{X}\rightarrow \mathbb{Y}$ defined by
	\begin{eqnarray}
		&&\mathbf{D}A_1(\Psi, H, v)(\Phi, K, u):= \sigma^{\prime}(\psi)\phi\Delta\psi+\sigma^{\prime\prime}(\psi)\phi|\nabla\psi|^{2}+(\sigma(\psi)-\sigma(0))\Delta\phi   \nonumber\\
		&& +2\sigma^{\prime}(\psi)\nabla\psi\cdot\nabla\phi-a^{\prime}(\psi)\phi +a^{\prime}(0)\phi:=\sum_{j=1}^{7}A_{1,j}, \label{diff1}\\
		&&\mathbf{D}A_3(\Psi, H, v)(\Phi, K, u) := \delta^{\prime}(\psi_{\Gamma})\phi_{\Gamma}\Delta_{\Gamma}\psi_{\Gamma}+\delta^{\prime\prime}(\psi_{\Gamma})\phi_{\Gamma}|\nabla_{\Gamma}\psi_{\Gamma}|^{2} \nonumber\\
		&& +(\delta(\psi_{\Gamma})-\delta(0))\Delta_{\Gamma}\phi_{\Gamma}+2\delta^{\prime}(\psi_{\Gamma})\nabla_{\Gamma}\psi_{\Gamma}\cdot\nabla_{\Gamma}\phi_{\Gamma} -\sigma^{\prime}(\psi_{\Gamma})\phi_{\Gamma}\partial_{\nu}\psi-(\sigma(\psi_{\Gamma})-\sigma(0))\partial_{\nu}\phi\nonumber\\
		&& -b^{\prime}(\psi_{\Gamma})\phi_{\Gamma}+b^{\prime}(0)\phi_{\Gamma}:=\sum_{j=1}^{9}A_{3,j}, \label{diff2}\\
		&&\mathbf{D}A_2(\Psi, H, v)(\Phi, K, u) := \sigma^{\prime}(\psi)\phi\Delta h+(\sigma(\psi)-\sigma(0))\Delta k-a^{\prime\prime}(\psi)\phi h\nonumber\\
		&&-(a^{\prime}(\psi)-a^{\prime}(0))k:=\sum_{j=1}^{4}A_{2,j}, \label{diff3}\\
		&&\mathbf{D}A_4(\Psi, H, v)(\Phi, K, u) := \delta^{\prime}(\psi_{\Gamma})\phi_{\Gamma}\Delta_{\Gamma} h_{\Gamma}+(\delta(\psi_{\Gamma})-\delta(0))\Delta_{\Gamma} k_{\Gamma}-\sigma^{\prime}(\psi_{\Gamma})\phi_{\Gamma}\partial_{\nu}h \nonumber\\
		&&-(\sigma(\psi_{\Gamma})-\sigma(0))\partial_{\nu}k-b^{\prime\prime}(\psi_{\Gamma})\phi_{\Gamma} h_{\Gamma}-(b^{\prime}(\psi_{\Gamma})-b^{\prime}(0))k_{\Gamma}:=\sum_{j=1}^{6}A_{4,j}, \label{diff4}
	\end{eqnarray}
	were we have denoted $A_{1,6}=A_{3,8}=0$ (Since $A^{\varepsilon}_{1,6}\longrightarrow 0$ and $A^{\varepsilon}_{3,8}\longrightarrow 0$ as $\varepsilon\longrightarrow 0$).
	Using the dominated convergence theorem, we can show that 
	\begin{eqnarray*}
		\frac{A((\Psi, H, v)+\varepsilon(\Phi, K, u))-A(\Psi, H, v)}{\varepsilon}\longrightarrow \mathbf{D}A(\Psi, H, u)(\Phi, K, v)\; \mbox{in}\; \mathbb{Y}\; \mbox{as}\; \varepsilon\longrightarrow 0.
	\end{eqnarray*}
	Indeed, using $A^\varepsilon_{1,7}=A_{1,7}$ and $A^\varepsilon_{3,9}=A_{3,9}=0$, then
	\begin{align*}
		&\left|\left|\frac{A((\Psi, H, v)+\varepsilon(\Phi, K, u))-A(\Psi, H, v)}{\varepsilon}- \mathbf{D}A(\Psi, H, v)(\Phi, K, u)\right|\right|^{2}_{\mathbb{Y}} \\
		&\leq  C\left(\sum_{j=1}^{6}\left[\|\mu(A_{1,j}^{\varepsilon}-A_{1,j})\|^{2}_{L^{2}(\Omega_{T})}+\|\mu_{4}(A_{1,j}^{\varepsilon}-A_{1,j})_{t}\|^{2}_{L^{2}(\Omega_{T})} \right]\right. \\
		& \left. +\sum_{j=1}^{8}\left[\|\mu(A_{3,j}^{\varepsilon}-A_{3,j})\|^{2}_{L^{2}(\Gamma_{T})}+\|\mu_{4}(A_{3,j}^{\varepsilon}-A_{3,j})_{t}\|^{2}_{L^{2}(\Gamma_{T})}   \right] \right.\\
		& \left. +\sum_{j=1}^{4}\|\mu(A_{2,j}^{\varepsilon}-A_{2,j})\|^{2}_{L^{2}(\Omega_{T})} + \sum_{j=1}^{6}\|\mu(A_{4,j}^{\varepsilon}-A_{4,j})\|^{2}_{L^{2}(\Gamma_{T})} \right).
	\end{align*}
	Let us start with the term  $\|\mu(A_{1,1}^{\varepsilon}-A_{1,1})\|^{2}_{L^{2}(\Omega_{T})}$.
	\begin{eqnarray*}
		\|\mu(A_{1,1}^{\varepsilon}-A_{1,1})\|^{2}_{L^{2}(\Omega_{T})}=\int_{\Omega_{T}}\mu^{2}\left|\frac{\sigma(\psi+\varepsilon\phi)-\sigma(\psi)-\varepsilon\phi\sigma^{\prime}(\psi)}{\varepsilon}\right|^{2}|\Delta\psi|^{2}.
	\end{eqnarray*}
	Since $\sigma$ is of class $C^{2}$, by  Taylor-Lagrange inequality, we have 
	$$\mu^{2}\left|\frac{\sigma(\psi+\varepsilon\phi)-\sigma(\psi)-\sigma^{\prime}(\psi)\varepsilon\phi}{\varepsilon}\right|^{2}|\Delta\psi|^{2}\leq \frac{M^{2}}{4} \mu^{2}|\phi|^{4}|\Delta\psi|^{2},$$
	where $M:=\displaystyle\sup_{r\in [-\rho_0, \rho_0]}|\sigma^{\prime\prime}(r)|$ and $\rho_0=C_0(\|(\Psi, H,v)\|_{\mathbb{X}}+\|(\Phi, K,u)\|_{\mathbb{X}})$. Moreover
	\begin{eqnarray*}
		\int_{\Omega_{T}}\mu^{2}|\phi|^{4}|\Delta\psi|^{2}&\leq& C\sup_{0\leq t \leq  T}\mu_{5}^{4}(t)\|\phi(\cdot,t)\|^{4}_{H^{2}(\Omega)}\int_{\Omega_{T}}\mu_{5}^{2}|\Delta\psi|^{2}\\
		&\leq & C\|(\Phi, K, u)\|_{\mathbb{X}}^{4}\|(\Psi, H, v)\|_{\mathbb{X}}^{2}.
	\end{eqnarray*}
	Then, using the dominated convergence theorem, we obtain 
	$$\|\mu(A_{1,1}^{\varepsilon}-A_{1,1})\|^{2}_{L^{2}(\Omega_{T})}\longrightarrow 0\quad \text{as}\quad \varepsilon\longrightarrow 0.$$
	For the term $\|\mu_{4}(A_{1,1}^{\varepsilon}-A_{1,1})_{t}\|^{2}_{L^{2}(\Omega_{T})}$, one has
	\begin{eqnarray*}
		&&\|\mu_{4}(A_{1,1}^{\varepsilon}-A_{1,1})_{t}\|^{2}_{L^{2}(\Omega_{T})}
		\leq C\left(\int_{\Omega_{T}}\mu_{4}^{2}\left|\frac{\sigma^{\prime}(\psi+\varepsilon\phi)-\sigma^{\prime}(\psi)-\varepsilon\phi\sigma^{\prime\prime}(\psi)}{\varepsilon}\right|^{2}|\psi_{t}|^{2}|\Delta\psi|^{2}\right. \\
		&& \left. + \int_{\Omega_{T}}\mu_{4}^{2}|\sigma^{\prime}(\psi+\varepsilon\phi)-\sigma^{\prime}(\psi)|^{2}|\phi_{t}|^{2}|\Delta\psi|^{2}  +\int_{\Omega_{T}}\mu_{4}^{2}\left|\frac{\sigma(\psi+\varepsilon\phi)-\sigma(\psi)-\varepsilon\phi\sigma^{\prime}(\psi)}{\varepsilon}\right|^{2}|\Delta\psi_{t}|^{2}\right).
	\end{eqnarray*}
	Using the dominated convergence theorem as above, we obtain 
	$$\|\mu_{4}(A_{1,1}^{\varepsilon}-A_{1,1})_{t}\|^{2}_{L^{2}(\Omega_{T})}\longrightarrow 0 \quad \text{as} \quad \varepsilon\longrightarrow 0.$$
	The same applies to other terms in fraction form whose denominator is $\varepsilon$ and the remaining terms are easy to study.
	Thus concluding that $A$ is Gateaux-differentiable at $(\Psi, H, v)$ and $A^{\prime}(\Psi, H, v)=\mathbf{D}A(\Psi, H, v)$.
	To conclude, it is sufficient to show $A^{\prime}:\mathbb{X}\rightarrow\mathcal{L}(\mathbb{X},\mathbb{Y})$ is continuous. Let $(\Psi^{n}, H^{n}, v^{n}), (\Psi, H, v), (\Phi, K, u)\in \mathbb{X}$ such that $(\Psi^{n}, H^{n}, v^{n})$ converges to $(\Psi, H, v)$ in $\mathbb{X}$ and we will prove that 
	\begin{eqnarray}
		\|(	A^{\prime}(\Psi^{n}, H^{n}, v^{n})(\Phi, K, u)-A^{\prime}(\Psi, H, v)(\Phi, K, u))\|_{\mathbb{Y}}\leq \varepsilon_{n} \|(\Phi, K, u)\|_{\mathbb{X}}, \label{bd0}
	\end{eqnarray}
	for some $(\varepsilon_{n})$ converging to $0$.\\
	Firstly, let $r>0$, such that $(\Psi, H, v)\in \overline{B}_{\mathbb{X}}(0,r)$, according to $(\Psi^{n}, H^{n}, v^{n})$ converges to $(\Psi, H, v)$ in $\mathbb{X}$ and \eqref{constant} , we can assume that 
	\begin{eqnarray*}
		\|\Psi^{n}\|_{L^{\infty}(0,T;\mathbb{L}^{\infty})}, \|\Psi\|_{L^{\infty}(0,T;\mathbb{L}^{\infty})}\leq C_{0}r \quad \mbox{for $n$ large enough}. 
	\end{eqnarray*}
	We consider
	$$B_{i,j}^{n}:=A_{i,j}(\Psi^{n}, H^{n}, v^{n})(\Phi, K, u)-A_{i,j}(\Psi, H, v)(\Phi, K, u).$$
	where $A_{i,j}$ are defined in \eqref{diff1}-\eqref{diff4}. Recall that $A_{1,6}=A_{3,8}=0$, and note that $B_{1,6}=B_{3,8}=0$ as well as $B_{1,7}=B_{3,9}=0$. Then
	\begin{eqnarray}
		&&\left|\left| A^{\prime}(\Psi^{n}, H^{n}, v^{n})(\Phi, K, u)-A^{\prime}(\Psi, H, v)(\Phi, K, u)\right|\right|^{2}_{\mathbb{Y}} \label{bd1} \\
		&&\leq  C\left(\sum_{j=1}^{5}\left[\|\mu B_{1,j}^{n}\|^{2}_{L^{2}(\Omega_{T})}+\|\mu_{4}(B_{1,j}^{n})_{t}\|^{2}_{L^{2}(\Omega_{T})} \right] \right. \nonumber \\
		&& \left. +\sum_{j=1}^{7}\left[\|\mu B_{3,j}^{n}\|^{2}_{L^{2}(\Gamma_{T})}+\|\mu_{4}(B_{3,j}^{n})_{t}\|^{2}_{L^{2}(\Gamma_{T})}   \right]  +\sum_{j=1}^{4}\|\mu B_{2,j}^{n}\|^{2}_{L^{2}(\Omega_{T})} + \sum_{j=1}^{6}\|\mu B_{4,j}^{n}\|^{2}_{L^{2}(\Gamma_{T})} \right). \nonumber
	\end{eqnarray}
	Let us start with the term $\|\mu B_{1,1}^{n}\|^{2}_{L^{2}(\Omega_{T})}$. Using $\sigma^{\prime}$ is Lipschitz-continuous and bounded on $J_0=[-C_0 r, C_0 r]$, we obtain 
	\begin{eqnarray*}
		\|\mu B_{1,1}^{n}\|^{2}_{L^{2}(\Omega_{T})} 
		&\leq & C\left(\int_{\Omega_{T}}\mu^{2}|(\psi^{n}-\psi)\phi \Delta\psi^{n}|^{2} + \int_{\Omega_{T}}\mu^{2}|\phi \Delta(\psi^{n}-\psi)|^{2}\right).
	\end{eqnarray*}
	Now using the fact that $H^{2}(\Omega)\hookrightarrow L^{\infty}(\Omega)$ is continuous we get
	\begin{eqnarray}
		\|\mu B_{1,1}^{n}\|^{2}_{L^{2}(\Omega_{T})}
		&\leq&  C\left(\int_{0}^{T}\left[\mu^{2}(t)\|(\psi^{n}-\psi)(\cdot,t)\|^{2}_{H^{2}(\Omega)}\|\phi(\cdot,t)\|^{2}_{H^{2}(\Omega)}\int_{\Omega}|\Delta\psi^{n}(x,t)|^{2}\right] \right. \nonumber\\
		&& \left. +\int_{0}^{T}\left[\mu^{2}(t)\|\phi(\cdot,t)\|^{2}_{H^{2}(\Omega)}\int_{\Omega}|\Delta(\psi^{n}-\psi)(x,t)|^{2}\right] \right) \nonumber\\
		&\leq&  C\sup_{0\leq t \leq T}\mu_{5}^{2}(t)\|(\psi^{n}-\psi)(\cdot,t)\|^{2}_{H^{2}(\Omega)}\sup_{0\leq t \leq T}\mu_{5}^{2}(t)\|\phi(\cdot,t)\|^{2}_{H^{2}(\Omega)} \nonumber\\
		&& \times\left(\int_{\Omega_{T}}\mu_{5}^{2}|\Delta\psi^{n}(x,t)|^{2}+1\right) \nonumber\\
		&\leq & C\varepsilon_{11}^{n}\|(\Phi, K, u)\|_{\mathbb{X}}^{2}, \label{bd2}
	\end{eqnarray}
	where $$\varepsilon_{11}^{n}:=\|(\Psi^{n}, H^{n}, v^{n})-(\Psi, H, v)\|_{\mathbb{X}}^{2}(\|(\Psi^{n}, H^{n}, v^{n})\|_{\mathbb{X}}^{2}+1).$$
	In the same way, we obtain 
		\begin{eqnarray}
			\|\mu B_{3,1}^{n}\|^{2}_{L^{2}(\Gamma_{T})}+\|\mu B_{3,5}^{n}\|^{2}_{L^{2}(\Gamma_{T})} +\|\mu B_{2,1}^{n}\|^{2}_{L^{2}(\Omega_{T})}+\|\mu B_{4,1}^{n}\|^{2}_{L^{2}(\Gamma_{T})} \leq C\varepsilon_{11}^{n}\|(\Phi, K, u)\|_{\mathbb{X}}^{2}, \label{bd3}
	\end{eqnarray}
were we have used the continuity of the normal derivative from $H^{2}(\Omega)$ to $L^{2}(\Gamma)$ for the term $\|\mu B_{3,5}^{n}\|^{2}_{L^{2}(\Gamma_{T})}$.\\
	Using $\sigma^{\prime}$ and $\sigma^{\prime\prime}$ are Lipschitz-continuous and bounded on $J_0$, we obtain 
	\begin{eqnarray}
		&&	\|\mu (B_{1,1}^{n})_{t}\|^{2}_{L^{2}(\Omega_{T})}
		\leq  C\left(\int_{\Omega_{T}}\mu^{2}|(\psi^{n}-\psi)\psi^{n}_{t}\phi \Delta\psi^{n}|^{2}  \right. \nonumber\\
		&&\left. + \int_{\Omega_{T}}\mu^{2}|(\psi^{n}_{t}-\psi_{t})\phi \Delta\psi^{n}|^{2} + \int_{\Omega_{T}}\mu^{2}|(\psi^{n}-\psi)\phi_{t}\Delta\psi^{n}|^{2} \right. \nonumber\\
		&& \left. +\int_{\Omega_{T}}\mu^{2}|(\psi^{n}-\psi)\phi\Delta\psi^{n}_{t}|^{2} +\int_{\Omega_{T}}\mu^{2}|\psi_{t}\phi\Delta(\psi^{n}-\psi)|^{2}  \right. \nonumber\\
		&& \left. +\int_{\Omega_{T}}\mu^{2}|\phi_{t}\Delta(\psi^{n}-\psi)|^{2} +\int_{\Omega_{T}}\mu^{2}|\phi\Delta(\psi^{n}_{t}-\psi_{t})|^{2}\right). \label{B11t}
	\end{eqnarray}
Let us analyze the first term in the right-hand side of \eqref{B11t}. Using the fact that the embeddings $H^{2}(\Omega)\hookrightarrow L^{\infty}(\Omega)$ is continuous and Remark \ref{Remark 3}, we obtain
	\begin{eqnarray*}
		&&\int_{\Omega_{T}}\mu^{2}|(\psi^{n}-\psi)\psi^{n}_{t}\phi \Delta\psi^{n}|^{2} \\
		&&\leq  \int_{0}^{T}\left[\mu^{2}(t)\|(\psi^{n}-\psi)(\cdot,t)\|^{2}_{H^{2}(\Omega)}\|\phi(\cdot,t)\|^{2}_{H^{2}(\Omega)}\|\psi^{n}_{t}(\cdot,t)\|^{2}_{H^{2}(\Omega)}\int_{\Omega}|\Delta\psi^{n}(x,t)|^{2} \right] \\
		&&\leq C\sup_{0\leq t \leq T}\mu_{5}^{2}(t)\|(\psi^{n}-\psi)(\cdot,t)\|^{2}_{H^{2}(\Omega)}\sup_{0\leq t \leq T}\mu_{5}^{2}(t)\|\phi(\cdot,t)\|^{2}_{H^{2}(\Omega)}\\
		&& \times\left(\int_{0}^{T}\mu_{5}^{2}(t)\|\psi^{n}_{t}(\cdot,t)\|^{2}_{H^{2}(\Omega)} \right)\sup_{0\leq t \leq T}\mu_{5}^{2}(t)\|\psi^{n}(\cdot,t)\|^{2}_{H^{2}(\Omega)}\\
		&& \leq C\|(\Psi^{n}, H^{n}, v^{n})-(\Psi, H, v)\|_{\mathbb{X}}^{2}\|(\Phi, K, u)\|_{\mathbb{X}}^{2}\|(\Psi^{n}, H^{n}, v^{n})\|_{\mathbb{X}}^{4}.
	\end{eqnarray*}
	The same applies to the other terms in the right-hand side of \eqref{B11t}, consequently 
	\begin{eqnarray}
		\|\mu (B_{1,1}^{n})_{t}\|^{2}_{L^{2}(\Omega_{T})} \leq  \tau^{n}_{1,1}\|(\Phi, K, u)\|_{\mathbb{X}}^{2}, \label{bd4}
	\end{eqnarray}
	where
	\begin{eqnarray*}
		\tau^{n}_{1,1}:=C\|(\Psi^{n}, H^{n}, v^{n})-(\Psi, H, v)\|_{\mathbb{X}}^{2}(\|(\Psi^{n}, H^{n}, v^{n})\|_{\mathbb{X}}^{4}+\|(\Psi^{n}, H^{n}, v^{n})\|_{\mathbb{X}}^{2}+\|(\Psi, H, v)\|_{\mathbb{X}}^{2}+1).		
	\end{eqnarray*}
	In a similar way to $\|\mu (B_{1,1}^{n})_{t}\|^{2}_{L^{2}(\Omega_{T})}$, we have 
		\begin{eqnarray}
			\|\mu (B_{3,1}^{n})_{t}\|^{2}_{L^{2}(\Omega_{T})} + \|\mu (B_{3,5}^{n})_{t}\|^{2}_{L^{2}(\Omega_{T})}\leq C\tau_{1,1}^{n}\|(\Phi, K, u)\|_{\mathbb{X}}^{2}. \label{bd5}
	\end{eqnarray}
In the computation of $\|\mu B_{1,2}^{n}\|^{2}_{L^{2}(\Omega_{T})}$, we use $\sigma^{\prime\prime}$ is Lipschitz-continuous and bounded on $J_0$ and  Cauchy Schwarz inequality. 
		\begin{eqnarray}
		&&\|\mu B_{1,2}^{n}\|^{2}_{L^{2}(\Omega_{T})} \leq  C\left(\int_{\Omega_{T}}\mu^{2}|(\sigma^{\prime\prime}(\psi^{n})-\sigma^{\prime\prime}(\psi))\phi |\nabla\psi^{n}|^{2}|^{2}  + \int_{\Omega_{T}}\mu^{2}|\sigma^{\prime\prime}(\psi)\phi (|\nabla\psi^{n}|^{2}-|\nabla\psi|^{2})|^{2}\right)   \nonumber\\
		&&\leq  C\left(\int_{\Omega_{T}}\mu^{2}|(\psi^{n}-\psi)\phi |\nabla\psi^{n}|^{2}|^{2} + \int_{\Omega_{T}}\mu^{2}|\phi (|\nabla\psi^{n}|^{2}-|\nabla\psi|^{2})|^{2} \right)  \nonumber\\
		&&\leq  C\left(\int_{0}^{T}\left[\mu^{2}(t)\|(\psi^{n}-\psi)(\cdot,t)\|^{2}_{H^{2}(\Omega)}\|\phi(\cdot,t)\|^{2}_{H^{2}(\Omega)}\int_{\Omega}|\nabla\psi^{n}|^{4}\right] \right.  \nonumber\\
		&& \left. +\int_{0}^{T}\left[\mu^{2}(t)\|\phi(\cdot,t)\|^{2}_{H^{2}(\Omega)}\left(\int_{\Omega}|\nabla(\psi^{n}-\psi)|^{4} \right)^{1/2} \left(\int_{\Omega}| \nabla(\psi^{n}+\psi)|^{4} \right)^{1/2} \right] \right) \nonumber
	\end{eqnarray} 
	By the embeddings $H^{2}(\Omega)\hookrightarrow W^{1,4}(\Omega)$ is continuous,  it follows that
	\begin{eqnarray}
	&&\|\mu B_{1,2}^{n}\|^{2}_{L^{2}(\Omega_{T})} \leq C\left(\int_{0}^{T}\mu^{2}(t)\|(\psi^{n}-\psi)(\cdot,t)\|^{2}_{H^{2}(\Omega)}\|\phi(\cdot,t)\|^{2}_{H^{2}(\Omega)}\|\psi^{n}(\cdot,t)\|^{4}_{H^{2}(\Omega)}  \right.   \nonumber\\
		&& \left. +\int_{0}^{T}\mu^{2}(t)\|\phi(\cdot,t)\|^{2}_{H^{2}(\Omega)}\|(\psi^{n}-\psi)(\cdot,t)\|^{2}_{H^{2}(\Omega)}\|(\psi^{n}+\psi)(\cdot,t)\|^{2}_{H^{2}(\Omega)}  \right)  \nonumber\\
		&&\leq C\sup_{0\leq t \leq T}\mu_{5}^{2}(t)\|(\psi^{n}-\psi)(\cdot,t)\|^{2}_{H^{2}(\Omega)}\sup_{0\leq t \leq T}\mu_{5}^{2}(t)\|\phi(\cdot,t)\|^{2}_{H^{2}(\Omega)}  \nonumber\\
		&& \times \left(\sup_{0\leq t \leq T}\mu_{5}^{4}(t)\|\psi^{n}(\cdot,t)\|^{4}_{H^{2}(\Omega)}+ \sup_{0\leq t \leq T}\mu_{5}^{2}(t)\|(\psi^{n}+\psi)(\cdot,t)\|^{2}_{H^{2}(\Omega)} \right) \nonumber\\
		&&\leq  C\varepsilon_{1,2}^{n}\|(\Phi, K, u)\|_{\mathbb{X}}^{2}, \label{bd6}
	\end{eqnarray} 
	where 
	\begin{eqnarray*}
		\varepsilon_{1,2}^{n}:=\|(\Psi^{n}, H^{n}, v^{n})-(\Psi, H, v)\|_{\mathbb{X}}^{2} \left(\|(\Psi^{n}, H^{n}, v^{n})\|_{\mathbb{X}}^{4}+ \|(\Psi^{n}, H^{n}, v^{n})+(\Psi, H, v)\|_{\mathbb{X}}^{2}\right).
	\end{eqnarray*}
	In the same way as for $\|\mu B_{1,2}^{n}\|^{2}_{L^{2}(\Omega_{T})}$, we obtain
		\begin{eqnarray}
			\|\mu B_{3,2}^{n}\|^{2}_{L^{2}(\Gamma_{T})} \leq C\varepsilon_{1,2}^{n}\|(\Phi, K, u)\|_{\mathbb{X}}^{2}. \label{bd7}
	\end{eqnarray}
	 Now we will compute $\|\mu (B_{1,2}^{n})_{t}\|^{2}_{L^{2}(\Omega_{T})}$.
	\begin{equation*}
		\begin{aligned}
		&\|\mu (B_{1,2}^{n})_{t}\|^{2}_{L^{2}(\Omega_{T})} \leq  C\left(\int_{\Omega_{T}}\mu^{2}|(\sigma^{\prime\prime\prime}(\psi^{n})-\sigma^{\prime\prime\prime}(\psi))\psi_{t}^{n}\phi |\nabla\psi^{n}|^{2}|^{2} \right.\\
		& \left. + \int_{\Omega_{T}}\mu^{2}|\sigma^{\prime\prime\prime}(\psi)(\psi_{t}^{n}-\psi_{t})\phi |\nabla\psi^{n}|^{2}|^{2} + \int_{\Omega_{T}}\mu^{2}|(\sigma^{\prime\prime}(\psi^{n})-\sigma^{\prime\prime}(\psi))\phi_{t} |\nabla\psi^{n}|^{2}|^{2}\right.\\
		&\left. + \int_{\Omega_{T}}\mu^{2}|(\sigma^{\prime\prime}(\psi^{n})-\sigma^{\prime\prime}(\psi))\phi \nabla\psi^{n}\cdot\nabla\psi^{n}_{t}|^{2}  \right.\\
		& \left. + \int_{\Omega_{T}}\mu^{2}|\sigma^{\prime\prime\prime}(\psi)\psi_{t}\phi (|\nabla\psi^{n}|^{2}-|\nabla\psi|^{2})|^{2}+ \int_{\Omega_{T}}\mu^{2}|\sigma^{\prime\prime}(\psi)\phi_{t} (|\nabla\psi^{n}|^{2}-|\nabla\psi|^{2})|^{2}  \right.\\
		& \left. + \int_{\Omega_{T}}\mu^{2}|\sigma^{\prime\prime}(\psi)\phi \nabla(\psi^{n}-\psi)\cdot\nabla\psi^{n}_{t}|^{2}+ \int_{\Omega_{T}}\mu^{2}|\sigma^{\prime\prime}(\psi)\phi \nabla\psi\cdot\nabla(\psi^{n}-\psi)|^{2}\right).
	\end{aligned}
	\end{equation*}
		As above, we arrive at 
		\begin{eqnarray}
			&&\|\mu (B_{1,2}^{n})_{t}\|^{2}_{L^{2}(\Omega_{T})}  \leq  C\left(\int_{\Omega_{T}}\mu^{2}|(\psi^{n}-\psi)\psi_{t}^{n}\phi |\nabla\psi^{n}|^{2}|^{2} \right.   \nonumber\\
		&& \left. + \int_{\Omega_{T}}\mu^{2}|(\psi_{t}^{n}-\psi_{t})\phi |\nabla\psi^{n}|^{2}|^{2} + \int_{\Omega_{T}}\mu^{2}|(\psi^{n}-\psi)\phi_{t} |\nabla\psi^{n}|^{2}|^{2} \right.  \nonumber\\
		&&\left. + \int_{\Omega_{T}}\mu^{2}|(\psi^{n}-\psi)\phi \nabla\psi^{n}\cdot\nabla\psi^{n}_{t}|^{2} + \int_{\Omega_{T}}\mu^{2}|\psi_{t}\phi (|\nabla\psi^{n}|^{2}-|\nabla\psi|^{2})|^{2} \right.  \nonumber\\
		&& \left. + \int_{\Omega_{T}}\mu^{2}|\phi_{t} (|\nabla\psi^{n}|^{2}-|\nabla\psi|^{2})|^{2} + \int_{\Omega_{T}}\mu^{2}|\phi \nabla(\psi^{n}-\psi)\cdot\nabla\psi^{n}_{t}|^{2} \right.  \nonumber\\
		&& \left. + \int_{\Omega_{T}}\mu^{2}|\phi \nabla\psi\cdot\nabla(\psi^{n}-\psi)|^{2} \right)  \nonumber\\
		&&\leq C\tau^{n}_{1,2}\|(\Phi, K, u)\|_{\mathbb{X}}^{2}, \label{bd8}
	\end{eqnarray} 
	where 
	\begin{eqnarray*}
		\tau^{n}_{1,2}&:=& \|(\Psi^{n}, H^{n}, v^{n})-(\Psi, H, v)\|_{\mathbb{X}}^{2}\left(\|(\Psi^{n}, H^{n}, v^{n})\|_{\mathbb{X}}^{6}+ \|(\Psi^{n}, H^{n}, v^{n})\|_{\mathbb{X}}^{4} \right. \\
		&& \left. + \|(\Psi^{n}, H^{n}, v^{n})\|_{\mathbb{X}}^{2} \|(\Psi, H, v)\|_{\mathbb{X}}^{2}  +\|(\Psi^{n}, H^{n}, v^{n})\|_{\mathbb{X}}^{2} \|(\Psi^{n}, H^{n}, v^{n})+(\Psi, H, v)\|_{\mathbb{X}}^{2} \right.\\
		&& \left. + \|(\Psi^{n}, H^{n}, v^{n})\|_{\mathbb{X}}^{2}+ \|(\Psi, H, v)\|_{\mathbb{X}}^{2}\right).
	\end{eqnarray*}
	In the same way as for $\|\mu (B_{1,2}^{n})_{t}\|^{2}_{L^{2}(\Gamma_{T})}$, we have
		\begin{eqnarray}
			\|\mu (B_{3,2}^{n})_{t}\|^{2}_{L^{2}(\Gamma_{T})} \leq C\tau_{1,2}^{n}\|(\Phi, K, u)\|_{\mathbb{X}}^{2}. \label{bd9}
	\end{eqnarray}
	For terms $\|\mu B_{1,3}^{n}\|^{2}_{L^{2}(\Omega_{T})}$ and $\|\mu (B_{1,3}^{n})_{t}\|^{2}_{L^{2}(\Omega_{T})}$, it is obvious that 
	\begin{eqnarray}
		&&\|\mu B_{1,3}^{n}\|^{2}_{L^{2}(\Omega_{T})} \leq  C\varepsilon_{1,3}^{n}\|(\Phi, K, u)\|_{\mathbb{X}}^{2}, \label{bd10}\\
		&&\|\mu (B_{1,3}^{n})_{t}\|^{2}_{L^{2}(\Omega_{T})} \leq  C\tau_{1,3}^{n}\|(\Phi, K, u)\|_{\mathbb{X}}^{2}, \label{bd11}
	\end{eqnarray}
	where 
	\begin{eqnarray*}
		&&\varepsilon_{1,3}^{n}:=\|(\Psi^{n}, H^{n}, v^{n})-(\Psi, H, v)\|_{\mathbb{X}}^{2},\\
		&&\tau_{1,3}^{n}:=\|(\Psi^{n}, H^{n}, v^{n})-(\Psi, H, v)\|_{\mathbb{X}}^{2}\left( \|(\Psi^{n}, H^{n}, v^{n})\|_{\mathbb{X}}^{2}+1\right).
	\end{eqnarray*}
	In the same way, we obtain
		\begin{eqnarray}
			&&\|\mu B_{3,3}^{n}\|^{2}_{L^{2}(\Gamma_{T})}+\|\mu B_{3,6}^{n}\|^{2}_{L^{2}(\Gamma_{T})}+\|\mu B_{2,2}^{n}\|^{2}_{L^{2}(\Omega_{T})}+\|\mu B_{2,4}^{n}\|^{2}_{L^{2}(\Gamma_{T})}  \nonumber \\
			&&+ \|\mu B_{4,2}^{n}\|^{2}_{L^{2}(\Gamma_{T})} + \|\mu B_{4,4}^{n}\|^{2}_{L^{2}(\Gamma_{T})}+ \|\mu B_{4,6}^{n}\|^{2}_{L^{2}(\Gamma_{T})}\leq C\varepsilon_{1,3}^{n}\|(\Phi, K, u)\|_{\mathbb{X}}^{2}, \label{bd12}\\
			&& \|\mu (B_{3,3}^{n})_{t}\|^{2}_{L^{2}(\Gamma_{T})} + \|\mu (B_{3,6}^{n})_{t}\|^{2}_{L^{2}(\Gamma_{T})}\leq C\tau_{1,3}^{n}\|(\Phi, K, u)\|_{\mathbb{X}}^{2},\label{bd13}
	\end{eqnarray}
were we have used the continuity of the normal derivative from $H^{2}(\Omega)$ to $L^{2}(\Gamma)$ for the terms $\|\mu B_{3,6}^{n}\|^{2}_{L^{2}(\Gamma_{T})}, \|\mu B_{4,4}^{n}\|^{2}_{L^{2}(\Gamma_{T})}$ and $\|\mu (B_{3,6}^{n})_{t}\|^{2}_{L^{2}(\Gamma_{T})}$.\\
	In the Computation of $\|\mu B_{1,4}^{n}\|^{2}_{L^{2}(\Omega_{T})}$, using $\sigma^{\prime}$ is Lipschitz-continuous and bounded on $\mathbb{R}$, we have 
	\begin{eqnarray*}
		\|\mu B_{1,4}^{n}\|^{2}_{L^{2}(\Omega_{T})} 
		&\leq & C\left(\int_{\Omega_{T}}\mu^{2}|(\psi^{n}-\psi)\nabla\psi^{n}\cdot\nabla\phi|^{2}  + \int_{\Omega_{T}}\mu^{2}|\nabla(\psi^{n}-\psi)\cdot\nabla\phi|^{2}\right).
	\end{eqnarray*}
	By the fact that $H^{2}(\Omega) \hookrightarrow L^{\infty}(\Omega)$ and $H^{2}(\Omega) \hookrightarrow W^{1,4}(\Omega)$ are continuous and Cauchy Schwarz inequality, it follows that 
	\begin{eqnarray}
		&&\|\mu B_{1,4}^{n}\|^{2}_{L^{2}(\Omega_{T})}\leq  C\left(\int_{0}^{T}\left[\mu^{2}(t)\|(\psi^{n}-\psi)(\cdot,t)\|^{2}_{H^{2}(\Omega)}\left(\int_{\Omega}|\nabla\psi^{n}|^{4} \right)^{1/2} \left(\int_{\Omega}|\nabla\phi|^{4} \right)^{1/2}\right]  \right.  \nonumber\\
		&& \left. + \int_{0}^{T}\left[\mu^{2}(t)\left(\int_{\Omega}|\nabla(\psi^{n}-\psi)|^{4}\d x \right)^{1/2}\left(\int_{\Omega}|\nabla\phi|^{4}\d x \right)^{1/2} \right] \right)  \nonumber\\
		&&\leq  C\left(\int_{0}^{T}\mu^{2}(t)\|(\psi^{n}-\psi)(\cdot,t)\|^{2}_{H^{2}(\Omega)}\|\psi^{n}(\cdot,t)\|^{2}_{H^{2}(\Omega)}\|\phi(\cdot,t)\|^{2}_{H^{2}(\Omega)}  \right.  \nonumber\\
		&& \left. + \int_{0}^{T}\mu^{2}(t)\|\psi^{n}(\cdot,t)-\psi(\cdot,t)\|^{2}_{H^{2}(\Omega)}\|\phi(\cdot,t)\|^{2}_{H^{2}(\Omega)}  \right) \nonumber\\
		&&\leq  C\sup_{0\leq t \leq T}\mu_{5}^{2}(t)\|(\psi^{n}-\psi)(\cdot,t)\|^{2}_{H^{2}(\Omega)}\sup_{0\leq t \leq T}\mu_{5}^{2}(t)\|\phi(\cdot,t)\|^{2}_{H^{2}(\Omega)}   \nonumber\\
		&& \times\left(\sup_{0\leq t \leq T}\mu_{5}^{2}(t)\|\psi^{n}(\cdot,t)\|^{2}_{H^{2}(\Omega)} +1\right)  \nonumber\\
		&&\leq C\varepsilon_{14}^{n}\|(\Phi, K, u)\|_{\mathbb{X}}^{2}, \label{bd14}
	\end{eqnarray} 
	where $$\varepsilon_{14}^{n}:=\|(\Psi^{n}, H^{n}, v^{n})-(\Psi, H, v)\|_{\mathbb{X}}^{2}(\|(\Psi^{n}, H^{n}, v^{n})\|_{\mathbb{X}}^{2}+1).$$
	In the same way as for $\|\mu B_{1,4}^{n}\|^{2}_{L^{2}(\Omega_{T})}$, we obtain
		\begin{eqnarray}
			\|\mu B_{3,4}^{n}\|^{2}_{L^{2}(\Gamma_{T})} \leq C\varepsilon_{14}^{n}\|(\Phi, K, u)\|_{\mathbb{X}}^{2}. \label{bd15}
	\end{eqnarray}
   Now we compute $\|\mu (B_{1,4}^{n})_{t}\|^{2}_{L^{2}(\Omega_{T})}$. Since $\sigma^{\prime\prime}$ and $\sigma^{\prime}$ are Lipschitz-continuous and bounded on $J_0$, then 
	\begin{eqnarray*}
		&&\|\mu (B_{1,4}^{n})_{t}\|^{2}_{L^{2}(\Omega_{T})} 
		\leq  C\left(\int_{\Omega_{T}}\mu^{2}|(\psi^{n}-\psi)\psi^{n}_{t}\nabla\psi^{n}\cdot\nabla\phi|^{2} \right. \\
		&& \left. + \int_{\Omega_{T}}\mu^{2}|(\psi^{n}_{t}-\psi_{t})\nabla\psi^{n}\cdot\nabla\phi|^{2} + \int_{\Omega_{T}}\mu^{2}|(\psi^{n}-\psi)\nabla\psi^{n}_{t}\cdot\nabla\phi|^{2} \right. \\
		&&\left. + \int_{\Omega_{T}}\mu^{2}|(\psi^{n}-\psi)\nabla\psi^{n}\cdot\nabla\phi_{t}|^{2} + \int_{\Omega_{T}}\mu^{2}|\psi_{t}\nabla(\psi^{n}-\psi)\cdot\nabla\phi|^{2} \right. \\
		&&\left.  + \int_{\Omega_{T}}\mu^{2}|\nabla(\psi^{n}_{t}-\psi_{t})\cdot\nabla\phi|^{2}+  \int_{\Omega_{T}}\mu^{2}|\nabla(\psi^{n}-\psi)\cdot\nabla\phi_{t}|^{2} \right).
	\end{eqnarray*}
	Now using the fact that $H^{2}(\Omega) \hookrightarrow L^{\infty}(\Omega)$ and $H^{2}(\Omega) \hookrightarrow W^{1,4}(\Omega)$ are continuous and Cauchy Schwarz inequality we get
	\begin{eqnarray}
		\|\mu (B_{1,4}^{n})_{t}\|^{2}_{L^{2}(\Omega_{T})}\leq C\tau_{14}^{n}\|(\Phi, K, u)\|_{\mathbb{X}}^{2}, \label{bd16}
	\end{eqnarray}
	where 
	\begin{eqnarray*}
		\tau_{14}^{n}&:=&\|(\Psi^{n}, H^{n}, v^{n})-(\Psi, H, v)\|_{\mathbb{X}}^{2}\left(\|(\Psi^{n}, H^{n}, v^{n})\|_{\mathbb{X}}^{4}+\|(\Psi^{n}, H^{n}, v^{n})\|_{\mathbb{X}}^{2}+\|(\Psi, H, v)\|_{\mathbb{X}}^{2}+1\right).
	\end{eqnarray*}
	In the same way, we obtain
		\begin{eqnarray}
			\|\mu (B_{3,4}^{n})_{t}\|^{2}_{L^{2}(\Gamma_{T})} \leq C\tau_{14}^{n}\|(\Phi, K, u)\|_{\mathbb{X}}^{2}. \label{bd17}
	\end{eqnarray}
	Let us compute of $\|\mu B_{1,5}^{n}\|^{2}_{L^{2}(\Omega_{T})}$ and $\|\mu (B_{1,5}^{n})_{t}\|^{2}_{L^{2}(\Omega_{T})}$. Using $a^{\prime}$ and $a^{\prime\prime}$ are Lipschitz-continuous and bounded on $J_0$, we obtain 
	\begin{eqnarray}
		&&\|\mu B_{1,5}^{n}\|^{2}_{L^{2}(\Omega_{T})}\leq C\varepsilon_{15}^{n}\|(\Phi, K, u)\|_{\mathbb{X}}^{2}, \label{bd18}\\
		&&\|\mu (B_{1,5}^{n})_{t}\|^{2}_{L^{2}(\Omega_{T})}\leq C\tau_{15}^{n}\|(\Phi, K, u)\|_{\mathbb{X}}^{2}, \label{bd19}
	\end{eqnarray}
	where 
	\begin{eqnarray*}
		&&\varepsilon_{15}^{n}:=\|(\Psi^{n}, H^{n}, v^{n})-(\Psi, H, v)\|_{\mathbb{X}}^{2},\\
		&&\tau_{15}^{n}:=\|(\Psi^{n}, H^{n}, v^{n})-(\Psi, H, v)\|_{\mathbb{X}}^{2}(\|(\Psi^{n}, H^{n}, v^{n})\|_{\mathbb{X}}^{2}+1).
	\end{eqnarray*}
	In the same way, we obtain
		\begin{eqnarray}
			&&\|\mu B_{3,7}^{n}\|^{2}_{L^{2}(\Gamma_{T})} \leq C\varepsilon_{15}^{n}\|(\Phi, K, u)\|_{\mathbb{X}}^{2}, \label{bd20}\\
			&&\|\mu (B_{3,7}^{n})_{t}\|^{2}_{L^{2}(\Gamma_{T})} \leq C\tau_{15}^{n}\|(\Phi, K, u)\|_{\mathbb{X}}^{2}. \label{bd21}
	\end{eqnarray}
	Now we compute $\|\mu B_{2,3}^{n}\|^{2}_{L^{2}(\Omega_{T})}$. Using $a^{\prime\prime}$ is Lipschitz-continuous and bounded on $J_0$, we obtain 
	\begin{eqnarray}
		&&\|\mu B_{2,3}^{n}\|^{2}_{L^{2}(\Omega_{T})}\leq C\left(\int_{\Omega_{T}}\mu^{2}|(a^{\prime\prime}(\psi^n)-a^{\prime\prime}(\psi))\phi h^{n}|^{2}   + \int_{\Omega_{T}}\mu^{2}|a^{\prime\prime}(\psi)\phi(h^{n}-h)|^{2}\right) \nonumber\\
		&&\leq C\left(\int_{\Omega_{T}}\mu^{2}|(\psi^{n}-\psi)\phi h^{n}|^{2} + \int_{\Omega_{T}}\mu^{2}|\phi(h^{n}-h)|^{2}\right) \nonumber\\
		&&\leq C\left( \int_{0}^{T}\left[\mu^{2}(t)\|(\psi^{n}-\psi)(\cdot,t)\|^{2}_{H^{2}(\Omega)}\|\phi(\cdot,t)\|^{2}_{H^{2}(\Omega)} \int_{\Omega}|h^{n}(x,t)|^{2} \right]  \right. \nonumber\\
		&& \left. + \int_{0}^{T}\left[\mu^{2}(t) \|\phi(\cdot,t)\|^{2}_{H^{2}(\Omega)} \int_{\Omega}|h^{n}(x,t)-h(x,t)|^{2}\d x \right] \right) \nonumber\\
		&&\leq C\left(\sup_{0\leq t \leq T}\mu_{5}^{2}(t)\|(\psi^{n}-\psi)(\cdot,t)\|^{2}_{H^{2}(\Omega)}\sup_{0\leq t \leq T}\mu_{5}^{2}(t)\|\phi(\cdot,t)\|^{2}_{H^{2}(\Omega)} \right. \nonumber\\
		&& \left. \times\left(\int_{\Omega_T}\mu_{0}^{2}|h^{n}|^{2} \right) + \sup_{0\leq t \leq T}\mu_{5}^{2}(t)\|\phi(\cdot,t)\|^{2}_{H^{2}(\Omega)}\int_{\Omega_T}\mu_{0}^{2}|h^{n}-h|^{2} \right) \nonumber\\
		&& \leq C\varepsilon_{2,3}^{n}\|(\Phi, K, u)\|_{\mathbb{X}}^{2}. \label{bd22}
	\end{eqnarray}
	where $\varepsilon_{2,3}^{n}:=\|(\Psi^{n}, H^{n}, v^{n})-(\Psi, H, v)\|_{\mathbb{X}}^{2}\left(\|(\Psi^{n}, H^{n}, v^{n})\|_{\mathbb{X}}^{2}+1\right)$.\\
Similar arguments lead to
		\begin{eqnarray}
			&& \|\mu B_{4,3}^{n}\|^{2}_{L^{2}(\Gamma_{T})} + \|\mu B_{4,5}^{n}\|^{2}_{L^{2}(\Gamma_{T})}\leq C\varepsilon_{2,3}^{n}\|(\Phi, K, u)\|_{\mathbb{X}}^{2}, \label{bd23}
	\end{eqnarray}
	were we have used the continuity of the normal derivative from $H^{2}(\Omega)$ to $L^{2}(\Gamma)$ for the term $\|\mu B_{4,3}^{n}\|^{2}_{L^{2}(\Gamma_{T})}$.\\
	Finally, from \eqref{bd1}-\eqref{bd3} and \eqref{bd4}-\eqref{bd23} , we obtain \eqref{bd0}.
\end{proof}

\subsection{Proof of Theorem \ref{Main result}}
Due to Lemmas \ref{Lemma3} and \ref{Lemma4}, the mapping $\Lambda$ is well-defined and continuously differentiable. Moreover 
\begin{eqnarray*}
	\Lambda^{\prime}(0,0, 0)(\Psi, H, v)=(\mathbf{L}_{1}\Psi-v\mathds{1}_{\omega}, \mathbf{L}_{2}\Psi, \mathbf{L}^{*}_{1}H-\theta\psi\mathds{1}_{\mathcal{O}}, \mathbf{L}^{*}_{2}H-\theta_{\Gamma}\psi_{\Gamma}\mathds{1}_{\Sigma}).
\end{eqnarray*}
Null controllability results of inhomogenuous linearized system \eqref{linearised cascade system} obtained in Propositions \ref{P1} and \ref{P2}, show that $\Lambda^{\prime}(0,0,0): \mathbb{X}\longrightarrow\mathbb{Y}$ is surjective. We conclude that Lyusternik-Graves' Theorem \ref{Lyusternik} can be applied to the operator $\Lambda$, in particular there exist $C=C(\Omega,\omega, T)$ and $\varepsilon>0$ such that, under the condition \eqref{assump source}, $\|F\|_{\mathbb{Y}_1}$ becomes sufficiently small and consequently, there exists $(\Psi, H, v)\in \mathbb{X}$ such that $\Lambda(\Psi, H, v)=(F, 0)$. 
As a result, $(\Psi, H)$ is the solution of the system \eqref{cascade quasi-linear system} associated to the control $v$, and 
the exponential growth of the weight $\mu_{0}$ as $t\rightarrow 0^{+}$ and $\mu_{0}H\in L^{2}(0,T,\mathbb{L}^{2})$, ensures that
\begin{eqnarray*}
	(h(\cdot,0),h_{\Gamma}(\cdot,0))=(0,0)\quad\mbox{in}\quad \Omega\times\Gamma.
\end{eqnarray*}
Moreover, there exists a constant $C>0$ such that
\begin{eqnarray*}
\|v\|_{L^{2}(0,T;H^{2}(\omega))}\leq\|(\Psi, H,v)\|_{\mathbb{X}}\leq C\|F\|_{\mathbb{Y}_1}.
\end{eqnarray*}
Therefore, by Proposition \ref{well-posdness of quasi}, we obtain $\Psi \in \mathfrak{F}_{T}$ which allows us to apply Lemma \ref{reformulation of problem}, which in turn establishes the proof of Theorem \eqref{Main result}. \qed

\section{Conclusion and final comments} \label{Secion7}
In this work, we have studied the insensitizing controls of a quasilinear reaction-diffusion equation of volume-surface type subject to dynamic boundary conditions. The strategy we have adopted relies on the null controllability of the inhomogeneous linearized cascad system and a local inversion theorem. Although we have considered reaction and diffusion coefficients that depend only on the state variable, our strategy can be adapted to more general diffusion and reaction coefficients depending on both the state and its gradients.

To the best of the authors knowledge, the insensitizing controls of quasilinear parabolic equations with general boundary conditions has not been considered before in the literature. Most of the previous works have dealt only with Dirichlet and Neumann boundary conditions, semi- linear equations with dynamic boundary conditions, as far as we know.

From a numerical perspective, it would be of much interest to investigate the approximation of null controls numerically based on the theoretical results we have obtained. This will eventually be done in a forthcoming paper.

\appendix
\section{Proof of Proposition \ref{well-posdness of quasi}} \label{Appendix A}
This paragraph is devoted to the proof of Proposition \ref{well-posdness of quasi}.\\
\begin{proof}[Proof of Proposition~{\upshape\ref{well-posdness of quasi}}]
Consider the mapping
	$\Upsilon:\mathfrak{F}_{T}\longrightarrow L^{2}(0,T;\mathbb{H}^{2})\times \mathbb{H}^{3}$ given by
\begin{eqnarray*}
	&&\Upsilon(\Psi):=(\Upsilon_1(\Psi), \Upsilon_2(\Psi), \Psi(\cdot, 0)), 
\end{eqnarray*}
where 
\begin{eqnarray*}
	&&\Upsilon_1(\Psi):=\psi_{t}-\nabla\cdot\left(\sigma(\psi)\nabla \psi\right) +a(\psi),\\
	&&\Upsilon_2(\Psi):=\psi_{\Gamma,t}-\nabla_{\Gamma}\cdot\left(\delta(\psi_{\Gamma})\nabla_{\Gamma} \psi_{\Gamma}\right)+\sigma(\psi_{\Gamma})
	\partial_{\nu} \psi + b(\psi_{\Gamma}).
\end{eqnarray*}
Firstly, using $\mathfrak{F}_T\hookrightarrow C([0,T],\mathbb{H}^{3})$ and \eqref{Bounded gradient}, we obtain $\|\Psi(\cdot, 0)\|_{\mathbb{H}^3}\leq \|\Psi\|_{\mathfrak{F}_T}$ and there exists $C>0$ such that $\|\Psi\|_{L^{\infty}(0,T;\mathbb{L}^{\infty})} \leq C\|\Psi\|_{\mathfrak{F}_T}$ and $ \|\nabla\Psi\|_{L^{\infty}(0,T;\mathbb{L}^{\infty})} \leq C\|\Psi\|_{\mathfrak{F}_T}$.\\
Now, by applying ideas similar to those in Lemmas \ref{Lemma3} and \ref{Lemma4}, we can establish that 
$\Upsilon$ is well-defined, continuously differentiable, and that
\begin{eqnarray*}
	\Upsilon^{\prime}(0)(\Psi)=(\mathbf{L}_1\Psi, \mathbf{L}_2\Psi, \Psi_0).
\end{eqnarray*}
Next, using semigroup theory through the form method, as illustrated in \cite[Proposition 2.4]{maniar2017null} or \cite[Theorem 4.2]{arendt2014maximal}, we can prove that $\Upsilon^{\prime}(0)$ is suejctive. Thus, we are able to apply Lyusternik-Graves' Theorem \ref{Lyusternik}, which guarantees the existence of solutions to \eqref{equation_quasi-linear} in $\mathfrak{F}_T$ that satisfy the estimate \eqref{quasi linear energy estimate} for sufficiently small data $(F,\Psi_{0})\in L^{2}(0,T;\mathbb{H}^{2}) \times\mathbb{H}^{3}$. The uniqueness of solutions to \eqref{quasi-linear equation well} in $\mathfrak{F}_T$ follows directly from an energy estimate and Grönwall's lemma.
Indeed, we put $Y=\Psi^{2}- \Psi^1$, where $\Psi^{1}, \Psi^2 \in\mathfrak{F}_T$ are both solutions of \eqref{quasi-linear equation well}. Obviously, $Y=(y,y_{\Gamma})$ satisfies the following equation:
\begin{equation} \label{sA1}
	\left\{
	\begin{aligned}
		&y_{t}-\nabla\cdot\left[\sigma(\psi^{2})\nabla \psi^2-\sigma(\psi^1)\nabla \psi^1\right] +(a(\psi^2)-a(\psi^1)) =0 & & \text {in}\; \Omega_T, \\
		&y_{\Gamma,t}-\nabla_{\Gamma}\cdot\left[\delta(\psi^{2}_{\Gamma})\nabla_{\Gamma} \psi^{2}_{\Gamma}-\delta(\psi^1_{\Gamma})\nabla_{\Gamma} \psi^1_{\Gamma}\right]+(\sigma(\psi^{2}_{\Gamma})\partial_{\nu}\psi^{2}-\sigma(\psi^1_{\Gamma})\partial_{\nu}\psi^1) \\
		& +(b(\psi^{2}_{\Gamma})-b(\psi^1_{\Gamma}))=0 & & \text {on}\;\Gamma_T, \\
		& y_{\Gamma}= y|_{\Gamma}  & & \text {on}\;\Gamma_T, \\
		& (y(\cdot,0), y_{\Gamma}(\cdot,0))=(0, 0) & & \text {in } \Omega\times\Gamma.
	\end{aligned}
	\right.
\end{equation}  
Multiplying the first equation of \eqref{sA1} by $y$ and integrating it in $\Omega$, one has
\begin{eqnarray*}
	&&\frac{1}{2}\frac{\d}{\d t}\int_{\Omega}|y|^{2}-\int_{\Gamma}(\sigma(\psi^{2}_{\Gamma})\partial_{\nu}\psi^{2}-\sigma(\psi^{1}_{\Gamma})\partial_{\nu}\psi^{1})y+\int_{\Omega}\sigma(\psi^{2})|\nabla y|^{2}  \\
	&&=-\int_{\Omega}(\sigma(\psi^{2})-\sigma(\psi^1))\nabla\psi^1\cdot\nabla y - \int_{\Omega}(a(\psi^{2})-a(\psi^{1}))y. 
\end{eqnarray*}
Multiplying the second equation of \eqref{sA1} by $y_{\Gamma}$ and integrating it on $\Gamma$, using the Stokes divergence formula \eqref{Stokes}, we obtain
\begin{eqnarray*}
	&&\frac{1}{2}\frac{\d}{\d t}\int_{\Gamma}|y_{\Gamma}|^{2}\d S+\int_{\Gamma}(\sigma(\psi^{2}_{\Gamma})\partial_{\nu}\psi^{2}-\sigma(\psi^{1}_{\Gamma})\partial_{\nu}\psi^{1})y+\int_{\Gamma}\delta(\psi^{2}_{\Gamma})|\nabla_{\Gamma} y_{\Gamma}|^{2}  \\
	&&=-\int_{\Gamma}(\delta(\psi^{2}_{\Gamma})-\delta(\psi^{1}_{\Gamma}))\nabla_{\Gamma}\psi^{1}_{\Gamma}\cdot\nabla_{\Gamma} y_{\Gamma} - \int_{\Gamma}(b(\psi^{2}_{\Gamma})-b(\psi^{1}_{\Gamma}))y_{\Gamma}.
\end{eqnarray*}
Next, we add these identities, using the fact that $a, b, \sigma$ and $\delta$ are Lipschitz- continuous on $[-C_0, C_0]$, $\sigma, \delta \geq\rho$, $\nabla\Psi^1\in L^{\infty}(0,T;\mathbb{L}^{\infty})$ (see \eqref{Bounded gradient})
and Young's inequality, we get
\begin{eqnarray*}
	&&\frac{1}{2}\frac{\d}{\d t}\|Y\|^{2}_{\mathbb{L}^{2}}
	+ \rho\|\nabla Y\|^{2}_{\mathbb{L}^{2}} \leq \frac{\rho}{2}\|\nabla Y\|^{2}_{\mathbb{L}^{2}}+C \|Y\|^{2}_{\mathbb{L}^{2}},
\end{eqnarray*}
where $C>0$ depends on $\sigma, \delta, a$ and $b$.
By Gronwall's
inequality, we find that
\begin{eqnarray}
	&&\|Y\|_{C([0,T];\mathbb{L}^{2})}
	\leq e^{CT}\|Y(\cdot, 0)\|_{\mathbb{L}^{2}}=0. \nonumber
\end{eqnarray}
This shows $\Psi^1=\Psi^2$.
\end{proof}
\section{Proof of Lemma \ref{reformulation of problem}} \label{Appendix B}
This paragraph is devoted to the proof of Lemma \ref{reformulation of problem}.\\
\begin{proof}[Proof of Lemma~{\upshape\ref{reformulation of problem}}]
For a fixed  $\widehat{\Psi_0}=(\widehat{\psi_0},\widehat{\psi_{0,\Gamma}})\in \mathbb{H}^{3}$ and $v\mathds{1}_{\omega}\in L^{2}(0,T;H^{2}(\Omega))$, we denote by $\Psi=(\psi,\psi_{\Gamma})\in\mathfrak{F}_T$ the solution of \eqref{equation_quasi-linear} associated with $\tau=\tau_{\Gamma}=0$ and for any $\tau\in (-1,1)$, $\Psi^{\tau}=(\psi^{\tau},\psi^{\tau}_{\Gamma})\in\mathfrak{F}_T$ denote the solution of \eqref{equation_quasi-linear} associated with $\tau_{\Gamma}=0$. Using \eqref{quasi linear energy estimate} and \eqref{Bounded gradient}, we obtain
	\begin{eqnarray*}
		\|\Psi^{\tau}\|_{L^{\infty}(0,T;\mathbb{L}^{\infty})}\leq C_0, \quad \|\Psi\|_{L^{\infty}(0,T;\mathbb{L}^{\infty})}\leq C_0 \quad\mbox{and}\quad \|\nabla\Psi\|_{L^{\infty}(0,T;\mathbb{L}^{\infty})}\leq C_0,
	\end{eqnarray*}
	for $C_0>0$, independent of $\tau$. According to the definition of the partial derivative, we have
	\begin{eqnarray}
		\frac{\partial \mathcal{J}}{\partial\tau}\bigg|_{\tau=\tau_{\Gamma}=0}=\lim_{\tau\rightarrow 0}\left( \frac{\theta}{2}\int_{\mathcal{O}_{T}}(\psi^{\tau}
		+\psi)\frac{\psi^{\tau}-\psi}{\tau}+\frac{\theta_{\Gamma}}{2}\int_{\Sigma_{T}}(\psi^{\tau}_{\Gamma}
		+\psi_{\Gamma})\frac{\psi^{\tau}_{\Gamma}-\psi_{\Gamma}}{\tau}\right). \label{partial derivative of J}
	\end{eqnarray} 
	We claim that 
	\begin{eqnarray}
		\Psi^{\tau}\longrightarrow \Psi \; \mbox{in}\; C([0,T];\mathbb{L}^{2}),\; \mbox{as} \; \tau\longrightarrow 0. \label{Limit1}
	\end{eqnarray}
	Indeed, we put $Y^{\tau}=\Psi^{\tau}- \Psi$. Obviously, $Y^{\tau}=(y^{\tau},y^{\tau}_{\Gamma})$ satisfies the following equation:
	\begin{equation*} 
		\left\{
		\begin{aligned}
			&y^{\tau}_{t}-\nabla\cdot\left[\sigma(\psi^{\tau})\nabla \psi^\tau-\sigma(\psi)\nabla \psi\right] +(a(\psi^\tau)-a(\psi)) =0 & & \text {in}\; \Omega_T, \\
			&y^{\tau}_{\Gamma,t}-\nabla_{\Gamma}\cdot\left[\delta(\psi^{\tau}_{\Gamma})\nabla_{\Gamma} \psi^{\tau}_{\Gamma}-\delta(\psi_{\Gamma})\nabla_{\Gamma} \psi_{\Gamma}\right]+(\sigma(\psi^{\tau}_{\Gamma})\partial_{\nu}\psi^{\tau}-\sigma(\psi_{\Gamma})\partial_{\nu}\psi) \\
			& +(b(\psi^{\tau}_{\Gamma})-b(\psi_{\Gamma}))=0 & & \text {on}\;\Gamma_T, \\
			& y^{\tau}_{\Gamma}= y^{\tau}|_{\Gamma}  & & \text {on}\;\Gamma_T, \\
			& (y^\tau(\cdot,0), y^{\tau}_{\Gamma}(\cdot,0))=(\tau\widehat{\psi_0}, 0) & & \text {in } \Omega\times\Gamma.
		\end{aligned}
		\right.
	\end{equation*}  
	To avoid repetition, by applying the same argument used in the uniqueness of solutions in \ref{Appendix A}, we arrive at 
	\begin{eqnarray}
		&&\|Y^{\tau}\|_{C([0,T];\mathbb{L}^{2})}
		\leq e^{CT}\|Y^{\tau}(\cdot,0)\|_{\mathbb{L}^{2}}= e^{CT}\tau\|\widehat{\psi_{0}}\|_{L^{2}(\Omega)}. \nonumber
	\end{eqnarray}
	This shows \eqref{Limit1}. Next, we claim that 
	\begin{eqnarray}
		\frac{\Psi^{\tau}- \Psi}{\tau}\longrightarrow Z \; \mbox{in}\; C([0,T];\mathbb{L}^{2}),\; \mbox{as} \; \tau\longrightarrow 0, \label{Limit2}
	\end{eqnarray}
	where $Z=(z,z_\Gamma)$ satisfies the following linear parabolic equation:
	\begin{equation} \label{s2 in proof lemma 2.1}
		\left\{
		\begin{aligned}
			&z_{t}-\nabla\cdot\left[\sigma(\psi)\nabla z +\sigma^{\prime}(\psi)z\nabla\psi\right] +a^{\prime}(\psi)z =0 & & \text {in}\; \Omega_T, \\
			&z_{\Gamma,t}-\nabla_{\Gamma}\cdot\left[\delta(\psi_{\Gamma})\nabla_{\Gamma} z_{\Gamma}+\delta^{\prime}(\psi_{\Gamma})z_{\Gamma}\nabla_{\Gamma}\psi_{\Gamma}\right]+\sigma(\psi_{\Gamma})
			\partial_{\nu} z\\
			&+\sigma^{\prime}(\psi_{\Gamma})\partial_{\nu}\psi z_{\Gamma} + b^{\prime}(\psi_{\Gamma})z_{\Gamma}=0 & & \text {on}\;\Gamma_T, \\
			& z_{\Gamma}= z_{|_{\Gamma}}  & & \text {on}\;\Gamma_T, \\
			& (z(\cdot,0), z_{\Gamma}(\cdot,0)) =(\widehat{\psi_{0}}, 0) & & \text {in } \Omega\times\Gamma. 
		\end{aligned}
		\right.
	\end{equation}
	To show this, set $W^{\tau}=\frac{\Psi^{\tau}- \Psi}{\tau}- Z$. It is easy to check that $W^{\tau}=(w^{\tau},w^{\tau}_{\Gamma})$ satisfies the following: 
	\small{
		\begin{equation} \label{s3 in proof lemma 2.1}
			\left\{
			\begin{aligned}
				&w^{\tau}_{t}-\nabla\cdot\left[\sigma(\psi^{\tau})\nabla w^{\tau}+D^{\tau}_{\sigma}(\psi,z)\nabla\psi+\left(\sigma(\psi^{\tau})-\sigma(\psi)\right)\nabla z\right] +D^{\tau}_{a}(\psi,z) =0 & & \text {in}\; \Omega_T, \\
				&w^{\tau}_{\Gamma,t}-\nabla_{\Gamma}\cdot\left[\delta(\psi^{\tau}_{\Gamma})\nabla_{\Gamma} w^{\tau}_{\Gamma}+D^{\tau}_{\delta}(\psi_{\Gamma},z_{\Gamma})\nabla_{\Gamma}\psi_{\Gamma}+\left(\delta(\psi^{\tau}_{\Gamma})-\delta(\psi_{\Gamma})\right)\nabla_{\Gamma} z_{\Gamma}\right]\\
				&+\sigma(\psi^{\tau}_{\Gamma})\partial_{\nu}w^{\tau}+D^{\tau}_{\sigma}(\psi_{\Gamma},z_{\Gamma})\partial_{\nu}\psi+D^{\tau}_{b}(\psi_{\Gamma},z_{\Gamma}) =0 & & \text {on}\; \Gamma_T, \\
				& w^{\tau}_{\Gamma}= w^{\tau}|_{\Gamma}  & & \text {on}\;\Gamma_T, \\
				& (w^\tau(\cdot,0), w^{\tau}_{\Gamma}(\cdot,0))=(0,0) & & \text {in } \Omega\times\Gamma.
			\end{aligned}
			\right.
	\end{equation}}
	where we have used the notation: $D^{\tau}_{\sigma}(\psi,z):=\frac{\sigma(\psi^{\tau})-\sigma(\psi)}{\tau}-\sigma^{\prime}(\psi)z.$\\
	Similarly for $D^{\tau}_{a}(\psi,z), D^{\tau}_{\delta}(\psi_{\Gamma},z_{\Gamma})$ and $D^{\tau}_{b}(\psi_{\Gamma},z_{\Gamma})$. Multiplying the first equation of \eqref{s3 in proof lemma 2.1} by $w^{\tau}$ and integrating it in $\Omega$, one has
	\begin{eqnarray*}
		&&\frac{1}{2}\frac{\d}{\d t}\int_{\Omega}|w^{\tau}|^{2}-\int_{\Gamma}\left[\sigma(\psi^{\tau})\partial_{\nu} w^{\tau}+D_{\sigma}^{\tau}(\psi,z)\partial_{\nu}\psi\right]w^{\tau} +\int_{\Omega}\sigma(\psi^{\tau})|\nabla w^{\tau}|^{2} =\int_{\Gamma}\left(\sigma(\psi^{\tau})-\sigma(\psi)\right)\partial_{\nu} z w^{\tau}\\
		&&-\int_{\Omega}\left[D_{\sigma}^{\tau}(\psi,z)\nabla\psi\cdot \nabla w^{\tau}+\left(\sigma(\psi^{\tau})-\sigma(\psi)\right)\nabla z\cdot \nabla w^{\tau}\right] -\int_{\Omega}D_{a}^{\tau}(\psi,z)w^{\tau}. 
	\end{eqnarray*}
	Multiplying the second equation of \eqref{s3 in proof lemma 2.1} by $w_{\Gamma}^{\tau}$ and integrating it on $\Gamma$, using the Stokes divergence formula \eqref{Stokes}, we obtain
	
	\begin{eqnarray*}
		&&\frac{1}{2}\frac{\d}{\d t}\int_{\Gamma}|w_{\Gamma}^{\tau}|^{2}  +\int_{\Gamma}\left(\sigma(\psi^{\tau}_{\Gamma})\partial_{\nu}w^{\tau}+D_{\sigma}^{\tau}(\psi_\Gamma, z_{\Gamma})\partial_{\nu}\psi\right)w_{\Gamma}^{\tau} +\int_{\Gamma}\delta(\psi^{\tau}_{\Gamma})|\nabla_{\Gamma} w^{\tau}_{\Gamma}|^{2}\\
		&&=-\int_{\Gamma}\left[D_{\delta}^{\tau}(\psi_{\Gamma},z_{\Gamma})\nabla_{\Gamma}\psi_{\Gamma}\cdot\nabla_{\Gamma} w^{\tau}_{\Gamma}  +\left(\delta(\psi^{\tau}_{\Gamma})-\delta(\psi_{\Gamma})\right)\nabla_{\Gamma} z_{\Gamma}\cdot\nabla_{\Gamma} w^{\tau}_{\Gamma}\right] -\int_{\Gamma}D_{b}^{\tau}(\psi_{\Gamma}, z_{\Gamma})w^{\tau}_{\Gamma}.
	\end{eqnarray*}
	Next, we add these identities, 
	using the fact that $\sigma$ and $\delta$ are Lipschitz-continuous on $[-C_0, C_0]$, $\sigma, \delta \geq\rho$, $|\sigma(\psi^{\tau})-\sigma(\psi)|\leq 2\displaystyle\max_{r\in [-C_0, C_0]}|\sigma(r)|$, $\partial_{\nu}z\in C([0,T]; L^{2}(\Gamma))$, $\|\nabla\Psi\|_{L^{\infty}(0,T;\mathbb{L}^{\infty})}\leq C_0$
	and Young's inequality, we obtain
	\begin{eqnarray*}
		&&\frac{1}{2}\frac{\d}{\d t}\|W^{\tau}\|^{2}_{\mathbb{L}^{2}}
		+ \rho\|\nabla W^{\tau}\|^{2}_{\mathbb{L}^{2}} \leq \frac{\rho}{2}\|\nabla W^{\tau}\|^{2}_{\mathbb{L}^{2}}+C \|W^{\tau}\|^{2}_{\mathbb{L}^{2}}+C\|\Psi^{\tau}-\Psi\|^{2}_{\mathbb{L}^{2}}\\
		&&+C\left(\int_{\Omega}\left|D_{\sigma}^{\tau}(\psi,z)\right|^{2} +\int_{\Omega}\left|D_{a}^{\tau}(\psi,z)\right|^{2}+\int_{\Gamma}\left|D_{\delta}^{\tau}(\psi_{\Gamma},z_{\Gamma})\right|^{2} +\int_{\Gamma}\left|D_{b}^{\tau}(\psi_{\Gamma},z_{\Gamma})\right|^{2}\right),
	\end{eqnarray*}
	where $C>0$ depends on $f, f_\Gamma, v, \psi_0, \psi_{0,\Gamma}, \sigma, \delta, a$ and $b$. The last four terms are similar to treat: 
	\begin{eqnarray*}
		D_{\sigma}^{\tau}(\psi,z)&=&\int_{0}^{1}\sigma^{\prime}((1-s)\psi^{\tau}+s\psi)ds 	\left(\frac{\psi^{\tau}-\psi}{\tau}\right)-\sigma^{\prime}(\psi)z\\
		&=& \int_{0}^{1}(\sigma^{\prime}((1-s)\psi^{\tau}+s\psi)-\sigma^{\prime}(\psi))ds z + \int_{0}^{1}\sigma^{\prime}((1-s)\psi^{\tau}+s\psi)ds w^{\tau}.
	\end{eqnarray*}
	For the first term, using $\sigma^{\prime}$ is lipschitz-continuous and $z$ is bounded, and for the second, using $\sigma^{\prime}$ is bounded, we obtain
	\begin{eqnarray*}
		\left|D_{\sigma}^{\tau}(\psi,z)\right|^{2} \leq C\left(|\psi^{\tau}-\psi|^{2} + |w^{\tau}|^{2}\right).
	\end{eqnarray*}
	Consequently, by Young's inequality, we obtain
	\begin{eqnarray*}
		\frac{1}{2}\frac{\d}{\d t}\|W^{\tau}\|^{2}_{\mathbb{L}^{2}}
		+ \rho\|\nabla W^{\tau}\|^{2}_{\mathbb{L}^{2}} \leq \frac{\rho}{2}\|\nabla W^{\tau}\|^{2}_{\mathbb{L}^{2}}+C \|W^{\tau}\|^{2}_{\mathbb{L}^{2}}+C\|\Psi^{\tau}-\Psi\|^{2}_{\mathbb{L}^{2}}
	\end{eqnarray*}
	Using Gronwall's
	inequality, we obtain 
	\begin{eqnarray}
		&&\|W^{\tau}\|_{C([0,T];\mathbb{L}^{2})}
		\leq e^{CT}\|\Psi^{\tau}-\Psi\|_{L^{2}(0,T;\mathbb{L}^{2})}. \nonumber
	\end{eqnarray}
	This yields \eqref{Limit2}. Now, using \eqref{partial derivative of J}, \eqref{Limit1} and \eqref{Limit2}, we find that 
	\begin{eqnarray*}
		\frac{\partial \mathcal{J}}{\partial\tau}\bigg|_{\tau=\tau_{\Gamma}=0}=\theta\int_{\mathcal{O}_{T}}\psi z + \theta_{\Gamma}\int_{\Sigma_{T}}\psi_{\Gamma}z_{\Gamma}.
	\end{eqnarray*}
	By a duality argument, we obtain
	\begin{eqnarray*}
		\frac{\partial \mathcal{J}}{\partial\tau}\bigg|_{\tau=\tau_{\Gamma}=0}
		=\langle Z(\cdot,0), H(\cdot,0) \rangle_{\mathbb{L}^{2}}-\langle(Z(\cdot,T), H(\cdot,T) \rangle_{\mathbb{L}^{2}}= \int_{\Omega}\widehat{\psi_{0}}h(\cdot,0).
	\end{eqnarray*}
	where $H=(h,h_\Gamma)$ satisfies the following adjoint linear parabolic equation of \eqref{s2 in proof lemma 2.1}:
	\begin{equation}
		\left\{
		\begin{aligned}
			&-h_{t}-\sigma(\psi)\Delta h +a^{\prime}(\psi)h =\theta\psi\mathds{1}_{\mathcal{O}} & & \text {in}\; \Omega_T, \\
			&-h_{\Gamma,t}-\delta(\psi_{\Gamma})\Delta_{\Gamma} h_{\Gamma}+\sigma(\psi_{\Gamma})
			\partial_{\nu}h+b^{\prime}(\psi_{\Gamma})h_{\Gamma}=\theta_{\Gamma}\psi_{\Gamma}\mathds{1}_{\Sigma} & & \text {on}\;\Gamma_T, \\
			& h_{\Gamma}= h_{|_{\Gamma}}  & & \text {on}\;\Gamma_T, \\
			& (h(\cdot,T),h_{\Gamma}(\cdot,T))=(0,0) & & \text {in } \Omega\times\Gamma. 
		\end{aligned}
		\right.
	\end{equation}
	Similarly, we can deduce 
	\begin{eqnarray*}
		\frac{\partial \Phi}{\partial\tau_\Gamma}\bigg|_{\tau=\tau_{\Gamma}=0}&=& \int_{\Gamma}\widehat{\psi_{0,\Gamma}}h_{\Gamma}(\cdot,0).
	\end{eqnarray*}
	Consequently, \eqref{definition of insensitizing} is satisfied for all $(\widehat{\psi_0},\widehat{\psi_{0,\Gamma}})\in \mathbb{H}^{3}$ with $\|(\widehat{\psi_0}, \widehat{\psi_{0,\Gamma}})\|_{\mathbb{H}^{3}}=1$ if, and only if $(h(\cdot,0),h_{\Gamma}(\cdot,0))=(0,0)$ in $\Omega\times\Gamma$.
\end{proof}


\end{document}